


\input amstex
\expandafter\ifx\csname mathdefs.tex\endcsname\relax
  \expandafter\gdef\csname mathdefs.tex\endcsname{}
\else \message{Hey!  Apparently you were trying to
  \string twice.   This does not make sense.} 
\errmessage{Please edit your file (probably \jobname.tex) and remove
any duplicate ``\string\input'' lines} \fi




\catcode`\X=12\catcode`\@=11

\def\n@wcount{\alloc@0\count\countdef\insc@unt}
\def\n@wwrite{\alloc@7\write\chardef\sixt@@n}
\def\n@wread{\alloc@6\read\chardef\sixt@@n}
\def\r@s@t{\relax}\def\v@idline{\par}\def\@mputate#1/{#1}
\def\l@c@l#1X{\firstpart.#1}\def\gl@b@l#1X{#1}\def\t@d@l#1X{{}}

\def\crossrefs#1{\ifx\all#1\let\tr@ce=\all\else\def\tr@ce{#1,}\fi
   \n@wwrite\cit@tionsout\openout\cit@tionsout=\jobname.cit 
   \write\cit@tionsout{\tr@ce}\expandafter\setfl@gs\tr@ce,}
\def\setfl@gs#1,{\def\@{#1}\ifx\@\empty\let\next=\relax
   \else\let\next=\setfl@gs\expandafter\xdef
   \csname#1tr@cetrue\endcsname{}\fi\next}
\def\m@ketag#1#2{\expandafter\n@wcount\csname#2tagno\endcsname
     \csname#2tagno\endcsname=0\let\tail=\all\xdef\all{\tail#2,}
   \ifx#1\l@c@l\let\tail=\r@s@t\xdef\r@s@t{\csname#2tagno\endcsname=0\tail}\fi
   \expandafter\gdef\csname#2cite\endcsname##1{\expandafter
     \ifx\csname#2tag##1\endcsname\relax?\else\csname#2tag##1\endcsname\fi
     \expandafter\ifx\csname#2tr@cetrue\endcsname\relax\else
     \write\cit@tionsout{#2tag ##1 cited on page \folio.}\fi}
   \expandafter\gdef\csname#2page\endcsname##1{\expandafter
     \ifx\csname#2page##1\endcsname\relax?\else\csname#2page##1\endcsname\fi
     \expandafter\ifx\csname#2tr@cetrue\endcsname\relax\else
     \write\cit@tionsout{#2tag ##1 cited on page \folio.}\fi}
   \expandafter\gdef\csname#2tag\endcsname##1{\expandafter
      \ifx\csname#2check##1\endcsname\relax
      \expandafter\xdef\csname#2check##1\endcsname{}%
      \else\immediate\write16{Warning: #2tag ##1 used more than once.}\fi
      \multit@g{#1}{#2}##1/X%
      \write\t@gsout{#2tag ##1 assigned number \csname#2tag##1\endcsname\space
      on page \number\count0.}%
   \csname#2tag##1\endcsname}}
\def\multit@g#1#2#3/#4X{\def\t@mp{#4}\ifx\t@mp\empty%
      \global\advance\csname#2tagno\endcsname by 1 
      \expandafter\xdef\csname#2tag#3\endcsname
      {#1\number\csname#2tagno\endcsnameX}%
   \else\expandafter\ifx\csname#2last#3\endcsname\relax
      \expandafter\n@wcount\csname#2last#3\endcsname
      \global\advance\csname#2tagno\endcsname by 1 
      \expandafter\xdef\csname#2tag#3\endcsname
      {#1\number\csname#2tagno\endcsnameX}
      \write\t@gsout{#2tag #3 assigned number \csname#2tag#3\endcsname\space
      on page \number\count0.}\fi
   \global\advance\csname#2last#3\endcsname by 1
   \def\t@mp{\expandafter\xdef\csname#2tag#3/}%
   \expandafter\t@mp\@mputate#4\endcsname
   {\csname#2tag#3\endcsname\lastpart{\csname#2last#3\endcsname}}\fi}
\def\t@gs#1{\def\all{}\m@ketag#1e\m@ketag#1s\m@ketag\t@d@l p
   \m@ketag\gl@b@l r \n@wread\t@gsin
   \openin\t@gsin=\jobname.tgs \re@der \closein\t@gsin
   \n@wwrite\t@gsout\openout\t@gsout=\jobname.tgs }
\outer\def\localtags{\t@gs\l@c@l}
\outer\def\globaltags{\t@gs\gl@b@l}
\outer\def\newlocaltag#1{\m@ketag\l@c@l{#1}}
\outer\def\newglobaltag#1{\m@ketag\gl@b@l{#1}}

\newif\ifpr@ 
\def\m@kecs #1tag #2 assigned number #3 on page #4.%
   {\expandafter\gdef\csname#1tag#2\endcsname{#3}
   \expandafter\gdef\csname#1page#2\endcsname{#4}
   \ifpr@\expandafter\xdef\csname#1check#2\endcsname{}\fi}
\def\re@der{\ifeof\t@gsin\let\next=\relax\else
   \read\t@gsin to\t@gline\ifx\t@gline\v@idline\else
   \expandafter\m@kecs \t@gline\fi\let \next=\re@der\fi\next}
\def\pretags#1{\pr@true\pret@gs#1,,}
\def\pret@gs#1,{\def\@{#1}\ifx\@\empty\let\n@xtfile=\relax
   \else\let\n@xtfile=\pret@gs \openin\t@gsin=#1.tgs \message{#1} \re@der 
   \closein\t@gsin\fi \n@xtfile}

\newcount\sectno\sectno=0\newcount\subsectno\subsectno=0
\newif\ifultr@local \def\ultralocal{\ultr@localtrue}
\def\firstpart{\number\sectno}
\def\lastpart#1{\ifcase#1 \or a\or b\or c\or d\or e\or f\or g\or h\or 
   i\or k\or l\or m\or n\or o\or p\or q\or r\or s\or t\or u\or v\or w\or 
   x\or y\or z \fi}

\def\resetall{\global\advance\sectno by 1\subsectno=0
   \gdef\firstpart{\number\sectno}\r@s@t}
\def\resetsub{\global\advance\subsectno by 1
   \gdef\firstpart{\number\sectno.\number\subsectno}\r@s@t}
\def\newsection#1\par{\resetall\vskip0pt plus.3\vsize\penalty-250
   \vskip0pt plus-.3\vsize\bigskip\bigskip
   \message{#1}\leftline{\bf#1}\nobreak\bigskip}
\def\subsection#1\par{\ifultr@local\resetsub\fi
   \vskip0pt plus.2\vsize\penalty-250\vskip0pt plus-.2\vsize
   \bigskip\smallskip\message{#1}\leftline{\bf#1}\nobreak\medskip}

\def\t@gsoff#1,{\def\@{#1}\ifx\@\empty\let\next=\relax\else\let\next=\t@gsoff
   \def\@@{p}\ifx\@\@@\else
   \expandafter\gdef\csname#1cite\endcsname##1{\zeigen{##1}}
   \expandafter\gdef\csname#1page\endcsname##1{?}
   \expandafter\gdef\csname#1tag\endcsname##1{\zeigen{##1}}\fi\fi\next}
\def\verbatimtags{\ifx\all\relax\else\expandafter\t@gsoff\all,\fi}
\def\zeigen#1{\hbox{$\langle$}#1\hbox{$\rangle$}}

\def\(#1){\edef\dot@g{\ifmmode\ifinner(\hbox{\noexpand\etag{#1}})
   \else\noexpand\eqno(\hbox{\noexpand\etag{#1}})\fi
   \else(\noexpand\ecite{#1})\fi}\dot@g}

\newif\ifbr@ck
\def\eat#1{}
\def\[#1]{\br@cktrue[\br@cket#1'X]}
\def\br@cket#1'#2X{\def\temp{#2}\ifx\temp\empty\let\next\eat
   \else\let\next\br@cket\fi
   \ifbr@ck\br@ckfalse\br@ck@t#1,X\else\br@cktrue#1\fi\next#2X}
\def\br@ck@t#1,#2X{\def\temp{#2}\ifx\temp\empty\let\neext\eat
   \else\let\neext\br@ck@t\def\temp{,}\fi
   \def\teemp{#1}\ifx\teemp\empty\else\rcite{#1}\fi\temp\neext#2X}
\def\resetbr@cket{\gdef\[##1]{[\rtag{##1}]}}
\def\references{\resetbr@cket\newsection References\par}

\newtoks\symb@ls\newtoks\s@mb@ls\newtoks\p@gelist\n@wcount\ftn@mber
    \ftn@mber=1\newif\ifftn@mbers\ftn@mbersfalse\newif\ifbyp@ge\byp@gefalse
\def\defm@rk{\ifftn@mbers\n@mberm@rk\else\symb@lm@rk\fi}
\def\n@mberm@rk{\xdef\m@rk{{\the\ftn@mber}}%
    \global\advance\ftn@mber by 1 }
\def\rot@te#1{\let\temp=#1\global#1=\expandafter\r@t@te\the\temp,X}
\def\r@t@te#1,#2X{{#2#1}\xdef\m@rk{{#1}}}
\def\b@@st#1{{$^{#1}$}}\def\str@p#1{#1}
\def\symb@lm@rk{\ifbyp@ge\rot@te\p@gelist\ifnum\expandafter\str@p\m@rk=1 
    \s@mb@ls=\symb@ls\fi\write\f@nsout{\number\count0}\fi \rot@te\s@mb@ls}
\def\byp@ge{\byp@getrue\n@wwrite\f@nsin\openin\f@nsin=\jobname.fns 
    \n@wcount\currentp@ge\currentp@ge=0\p@gelist={0}
    \re@dfns\closein\f@nsin\rot@te\p@gelist
    \n@wread\f@nsout\openout\f@nsout=\jobname.fns }
\def\m@kelist#1X#2{{#1,#2}}
\def\re@dfns{\ifeof\f@nsin\let\next=\relax\else\read\f@nsin to \f@nline
    \ifx\f@nline\v@idline\else\let\t@mplist=\p@gelist
    \ifnum\currentp@ge=\f@nline
    \global\p@gelist=\expandafter\m@kelist\the\t@mplistX0
    \else\currentp@ge=\f@nline
    \global\p@gelist=\expandafter\m@kelist\the\t@mplistX1\fi\fi
    \let\next=\re@dfns\fi\next}
\def\symbols#1{\symb@ls={#1}\s@mb@ls=\symb@ls} 
\def\bigsymbol{\textstyle}
\symbols{\bigsymbol\ast,\dagger,\ddagger,\sharp,\flat,\natural,\star}
\def\ftnumbers{\ftn@mberstrue} \def\ftsymbols{\ftn@mbersfalse}
\def\paginal{\byp@ge} \def\resetftnumbers{\ftn@mber=1}
\def\ftnote#1{\defm@rk\expandafter\expandafter\expandafter\footnote
    \expandafter\b@@st\m@rk{#1}}

\long\def\jump#1\endjump{}
\def\ssum{\mathop{\lower .1em\hbox{$\textstyle\Sigma$}}\nolimits}

\def\qed{\nobreak\kern 1em \vrule height .5em width .5em depth 0em}
\def\newneq{\hbox{\rlap{\hbox to 1\wd9{\hss$=$\hss}}\raise .1em 
   \hbox to 1\wd9{\hss$\scriptscriptstyle/$\hss}}}
\def\subsetne{\setbox9 = \hbox{$\subset$}\mathrel{\hbox{\rlap
   {\lower .4em \newneq}\raise .13em \hbox{$\subset$}}}}
\def\supsetne{\setbox9 = \hbox{$\subset$}\mathrel{\hbox{\rlap
   {\lower .4em \newneq}\raise .13em \hbox{$\supset$}}}}

\def\vbar{\mathchoice{\vrule height6.3ptdepth-.5ptwidth.8pt\kern-.8pt}
   {\vrule height6.3ptdepth-.5ptwidth.8pt\kern-.8pt}
   {\vrule height4.1ptdepth-.35ptwidth.6pt\kern-.6pt}
   {\vrule height3.1ptdepth-.25ptwidth.5pt\kern-.5pt}}
\def\f@dge{\mathchoice{}{}{\mkern.5mu}{\mkern.8mu}}
\def\b@c#1#2{{\rm \mkern#2mu\vbar\mkern-#2mu#1}}
\def\b@b#1{{\rm I\mkern-3.5mu #1}}
\def\b@a#1#2{{\rm #1\mkern-#2mu\f@dge #1}}
\def\bb#1{{\count4=`#1 \advance\count4by-64 \ifcase\count4\or\b@a A{11.5}\or
   \b@b B\or\b@c C{5}\or\b@b D\or\b@b E\or\b@b F \or\b@c G{5}\or\b@b H\or
   \b@b I\or\b@c J{3}\or\b@b K\or\b@b L \or\b@b M\or\b@b N\or\b@c O{5} \or
   \b@b P\or\b@c Q{5}\or\b@b R\or\b@a S{8}\or\b@a T{10.5}\or\b@c U{5}\or
   \b@a V{12}\or\b@a W{16.5}\or\b@a X{11}\or\b@a Y{11.7}\or\b@a Z{7.5}\fi}}

\catcode`\X=11 \catcode`\@=12

\expandafter\ifx\csname citeadd.tex\endcsname\relax
\expandafter\gdef\csname citeadd.tex\endcsname{}
\else \message{Hey!  Apparently you were trying to
\string twice.   This does not make sense.} 
\errmessage{Please edit your file (probably \jobname.tex) and remove
any duplicate ``\string\input'' lines} \fi

\sectno=-2   
\localtags
\ifx\shlhetal\undefinedcontrolsequenc\let\shlhetal\relax\fi
\NoBlackBoxes
\define\mr{\medskip\roster}
\define\nl{\newline}
\define\sn{\smallskip\noindent}
\define\mn{\medskip\noindent}
\define\bn{\bigskip\noindent}
\define\ub{\underbar}
\define\wilog{\text{without loss of generality}}
\define\ermn{\endroster\medskip\noindent}
\define\dbca{\dsize\bigcap}
\define\dbcu{\dsize\bigcup}
\documentstyle {amsppt}
\topmatter
\title{More Constructions for Boolean Algebras \\
 Sh652} \endtitle
\rightheadtext{More Constructions}
\author {Saharon Shelah \thanks {\null\newline This research was partially
supported by the Israel Foundation for Science \null\newline
I would like to thank Alice Leonhardt for the beautiful typing. \null\newline
First Typed - 97/May/12 \null\newline 
Latest Revision - 98/May/13} \endthanks} \endauthor
\affil{Institute of Mathematics\\
 The Hebrew University\\
 Jerusalem, Israel
 \medskip
 Rutgers University\\
 Mathematics Department\\
 New Brunswick, NJ USA} \endaffil
\endtopmatter
\document  

\expandafter\ifx\csname alice2jlem.tex\endcsname\relax
  \expandafter\gdef\csname alice2jlem.tex\endcsname{}
\else \message{Hey!  Apparently you were trying to
\string  twice.   This does not make sense.}
\errmessage{Please edit your file (probably \jobname.tex) and remove
any duplicate ``\string\input'' lines} \fi

\expandafter\ifx\csname bib4plain.tex\endcsname\relax
  \expandafter\gdef\csname bib4plain.tex\endcsname{}
\else \message{Hey!  Apparently you were trying to \string twice.   This does not make sense.}
\errmessage{Please edit your file (probably \jobname.tex) and remove
any duplicate ``\string\input'' lines} \fi

\def\renewcommand{\newcommand}	       
\edef\cite{\the\catcode`@}%
\catcode`@ = 11
\let\@oldatcatcode = \cite
\chardef\@letter = 11
\chardef\@other = 12
%
%
%
%
\def\@innerdef#1#2{\edef#1{\expandafter\noexpand\csname #2\endcsname}}%
%
%
\@innerdef\@innernewcount{newcount}%
\@innerdef\@innernewdimen{newdimen}%
\@innerdef\@innernewif{newif}%
\@innerdef\@innernewwrite{newwrite}%
%
%
%
\def\@gobble#1{}%
%
%
%
\ifx\inputlineno\@undefined
   \let\@linenumber = \empty 
\else
   \def\@linenumber{\the\inputlineno:\space}%
\fi
%
%
%
\def\@futurenonspacelet#1{\def\cs{#1}%
   \afterassignment\@stepone\let\@nexttoken=
}%
\begingroup 
\def\\{\global\let\@stoken= }%
\\ 
\endgroup
\def\@stepone{\expandafter\futurelet\cs\@steptwo}%
\def\@steptwo{\expandafter\ifx\cs\@stoken\let\@@next=\@stepthree
   \else\let\@@next=\@nexttoken\fi \@@next}%
\def\@stepthree{\afterassignment\@stepone\let\@@next= }%
%
%
%
\def\@getoptionalarg#1{%
   \let\@optionaltemp = #1%
   \let\@optionalnext = \relax
   \@futurenonspacelet\@optionalnext\@bracketcheck
}%
%
%
\def\@bracketcheck{%
   \ifx [\@optionalnext
      \expandafter\@@getoptionalarg
   \else
      \let\@optionalarg = \empty
      \expandafter\@optionaltemp
   \fi
}%
\def\@@getoptionalarg[#1]{%
   \def\@optionalarg{#1}%
   \@optionaltemp
}%
%
%
%
\def\@nnil{\@nil}%
\def\@fornoop#1\@@#2#3{}%
\def\@for#1:=#2\do#3{%
   \edef\@fortmp{#2}%
   \ifx\@fortmp\empty \else
      \expandafter\@forloop#2,\@nil,\@nil\@@#1{#3}%
   \fi
}%
\def\@forloop#1,#2,#3\@@#4#5{\def#4{#1}\ifx #4\@nnil \else
       #5\def#4{#2}\ifx #4\@nnil \else#5\@iforloop #3\@@#4{#5}\fi\fi
}%
\def\@iforloop#1,#2\@@#3#4{\def#3{#1}\ifx #3\@nnil
       \let\@nextwhile=\@fornoop \else
      #4\relax\let\@nextwhile=\@iforloop\fi\@nextwhile#2\@@#3{#4}%
}%
%
%
%
\@innernewif\if@fileexists
\def\@testfileexistence{\@getoptionalarg\@finishtestfileexistence}%
\def\@finishtestfileexistence#1{%
   \begingroup
      \def\extension{#1}%
      \immediate\openin0 =
         \ifx\@optionalarg\empty\jobname\else\@optionalarg\fi
         \ifx\extension\empty \else .#1\fi
         \space
      \ifeof 0
         \global\@fileexistsfalse
      \else
         \global\@fileexiststrue
      \fi
      \immediate\closein0
   \endgroup
}%
%
%
%
%
\def\bibliographystyle#1{%
   \@readauxfile
   \@writeaux{\string\bibstyle{#1}}%
}%
\let\bibstyle = \@gobble
%
%
\let\bblfilebasename = \jobname
\def\bibliography#1{%
   \@readauxfile
   \@writeaux{\string\bibdata{#1}}%
   \@testfileexistence[\bblfilebasename]{bbl}%
   \if@fileexists
      \nobreak
      \@readbblfile
   \fi
}%
\let\bibdata = \@gobble
%
%
\def\nocite#1{%
   \@readauxfile
   \@writeaux{\string\citation{#1}}%
}%
\@innernewif\if@notfirstcitation
%
%
\def\cite{\@getoptionalarg\@cite}%
%
%
\def\@cite#1{%
   \let\@citenotetext = \@optionalarg
   \printcitestart
   \nocite{#1}%
   \@notfirstcitationfalse
   \@for \@citation :=#1\do
   {%
      \expandafter\@onecitation\@citation\@@
   }%
   \ifx\empty\@citenotetext\else
      \printcitenote{\@citenotetext}%
   \fi
   \printcitefinish
}%
\def\@onecitation#1\@@{%
   \if@notfirstcitation
      \printbetweencitations
   \fi
   \expandafter \ifx \csname\@citelabel{#1}\endcsname \relax
      \if@citewarning
         \message{\@linenumber Undefined citation `#1'.}%
      \fi
      \expandafter\gdef\csname\@citelabel{#1}\endcsname{%
\strut
\vadjust{\vskip-\dp\strutbox
\vbox to 0pt{\vss\parindent0cm \leftskip=\hsize 
\advance\leftskip3mm
\advance\hsize 4cm\strut\openup-4pt 
\rightskip 0cm plus 1cm minus 0.5cm ?  #1 ?\strut}}
         {\tt
            \escapechar = -1
            \nobreak\hskip0pt
            \expandafter\string\csname#1\endcsname
            \nobreak\hskip0pt
         }%
      }%
   \fi
   \csname\@citelabel{#1}\endcsname
   \@notfirstcitationtrue
}%
%
%
\def\@citelabel#1{b@#1}%
%
%
\def\@citedef#1#2{\expandafter\gdef\csname\@citelabel{#1}\endcsname{#2}}%
%
%
%
\def\@readbblfile{%
   \ifx\@itemnum\@undefined
      \@innernewcount\@itemnum
   \fi
   \begingroup
      \def\begin##1##2{%
         \setbox0 = \hbox{\biblabelcontents{##2}}%
         \biblabelwidth = \wd0
      }%
      \def\end##1{}
      %
      %
      \@itemnum = 0
      \def\bibitem{\@getoptionalarg\@bibitem}%
      \def\@bibitem{%
         \ifx\@optionalarg\empty
            \expandafter\@numberedbibitem
         \else
            \expandafter\@alphabibitem
         \fi
      }%
      \def\@alphabibitem##1{%
         \expandafter \xdef\csname\@citelabel{##1}\endcsname {\@optionalarg}%
         \ifx\biblabelprecontents\@undefined
            \let\biblabelprecontents = \relax
         \fi
         \ifx\biblabelpostcontents\@undefined
            \let\biblabelpostcontents = \hss
         \fi
         \@finishbibitem{##1}%
      }%
      \def\@numberedbibitem##1{%
         \advance\@itemnum by 1
         \expandafter \xdef\csname\@citelabel{##1}\endcsname{\number\@itemnum}%
         \ifx\biblabelprecontents\@undefined
            \let\biblabelprecontents = \hss
         \fi
         \ifx\biblabelpostcontents\@undefined
            \let\biblabelpostcontents = \relax
         \fi
         \@finishbibitem{##1}%
      }%
      \def\@finishbibitem##1{%
         \biblabelprint{\csname\@citelabel{##1}\endcsname}%
         \@writeaux{\string\@citedef{##1}{\csname\@citelabel{##1}\endcsname}}%
         \ignorespaces
      }%
      %
      %
      \let\em = \bblem
      \let\newblock = \bblnewblock
      \let\sc = \bblsc
      \frenchspacing
      \clubpenalty = 4000 \widowpenalty = 4000
      \tolerance = 10000 \hfuzz = .5pt
      \everypar = {\hangindent = \biblabelwidth
                      \advance\hangindent by \biblabelextraspace}%
      \bblrm
      \parskip = 1.5ex plus .5ex minus .5ex
      \biblabelextraspace = .5em
      \bblhook
      \input \bblfilebasename.bbl
   \endgroup
}%
%
%
\@innernewdimen\biblabelwidth
\@innernewdimen\biblabelextraspace
%
%
%
\def\biblabelprint#1{%
   \noindent
   \hbox to \biblabelwidth{%
      \biblabelprecontents
      \biblabelcontents{#1}%
      \biblabelpostcontents
   }%
   \kern\biblabelextraspace
}%
%
%
%
\def\biblabelcontents#1{{\bblrm [#1]}}%
%
%
\def\bblrm{\rm}%
%
%
\def\bblem{\it}%
%
%
\def\bblsc{\ifx\@scfont\@undefined
              \font\@scfont = cmcsc10
           \fi
           \@scfont
}%
%
%
\def\bblnewblock{\hskip .11em plus .33em minus .07em }%
%
%
\let\bblhook = \empty
%
%
%
\def\printcitestart{[}
\def\printcitefinish{]}
\def\printbetweencitations{, }
\def\printcitenote#1{, #1}
%
%
%
\let\citation = \@gobble
%
%
%
\@innernewcount\@numparams
%
%
\def\newcommand#1{%
   \def\@commandname{#1}%
   \@getoptionalarg\@continuenewcommand
}%
%
%
\def\@continuenewcommand{%
   \@numparams = \ifx\@optionalarg\empty 0\else\@optionalarg \fi \relax
   \@newcommand
}%
%
%
\def\@newcommand#1{%
   \def\@startdef{\expandafter\edef\@commandname}%
   \ifnum\@numparams=0
      \let\@paramdef = \empty
   \else
      \ifnum\@numparams>9
         \errmessage{\the\@numparams\space is too many parameters}%
      \else
         \ifnum\@numparams<0
            \errmessage{\the\@numparams\space is too few parameters}%
         \else
            \edef\@paramdef{%
               \ifcase\@numparams
                  \empty  No arguments.
               \or ####1%
               \or ####1####2%
               \or ####1####2####3%
               \or ####1####2####3####4%
               \or ####1####2####3####4####5%
               \or ####1####2####3####4####5####6%
               \or ####1####2####3####4####5####6####7%
               \or ####1####2####3####4####5####6####7####8%
               \or ####1####2####3####4####5####6####7####8####9%
               \fi
            }%
         \fi
      \fi
   \fi
   \expandafter\@startdef\@paramdef{#1}%
}%
%
%
%
%
\def\@readauxfile{%
   \if@auxfiledone \else 
      \global\@auxfiledonetrue
      \@testfileexistence{aux}%
      \if@fileexists
         \begingroup
            \endlinechar = -1
            \catcode`@ = 11
            \input \jobname.aux
         \endgroup
      \else
         \message{\@undefinedmessage}%
         \global\@citewarningfalse
      \fi
      \immediate\openout\@auxfile = \jobname.aux
   \fi
}%
%
%
\newif\if@auxfiledone
\ifx\noauxfile\@undefined \else \@auxfiledonetrue\fi
%
%
%
%
\@innernewwrite\@auxfile
\def\@writeaux#1{\ifx\noauxfile\@undefined \write\@auxfile{#1}\fi}%
%
%
%
\ifx\@undefinedmessage\@undefined
   \def\@undefinedmessage{No .aux file; I won't give you warnings about
                          undefined citations.}%
\fi
%
%
\@innernewif\if@citewarning
\ifx\noauxfile\@undefined \@citewarningtrue\fi
%
%
%
\catcode`@ = \@oldatcatcode


\def\widestnumber#1#2{}

\def\rm{\fam0 \tenrm}

\def\fakesubhead#1\endsubhead{\bigskip\noindent{\bf#1}\par}


%
%
%

%

\font\textrsfs=rsfs10
\font\scriptrsfs=rsfs7
\font\scriptscriptrsfs=rsfs5

\newfam\rsfsfam
\textfont\rsfsfam=\textrsfs
\scriptfont\rsfsfam=\scriptrsfs
\scriptscriptfont\rsfsfam=\scriptscriptrsfs

\edef\oldcatcodeofat{\the\catcode`\@}
\catcode`\@11

\def\Cal@@#1{\noaccents@ \fam \rsfsfam #1}

\catcode`\@\oldcatcodeofat

\newpage

\head {Annotated Content} \endhead  \resetall
\bn
\S1 $\quad$  The depth of free product may be bigger than the depths of those
multiplied
\mr
\item "{{}}" [We construct, in ZFC, for any Boolean Algebra $B$, and cardinal
$\kappa$ Boolean Algebras $B_1,B_2$ extending $B$ such that the depth of the 
free product of $B_1,B_2$ over $B$ is strictly larger than the depths of $B_1$
and of $B_2$ than $\kappa$.  Thus, we answer problem 10 of Monk 
\cite{M}.  We give a condition $\boxtimes_{\lambda,\mu,\theta}$ which implies
that for some Boolean Algebra $A = A_\theta$ there are 
$B_1 = B^1_{\lambda,\mu,\theta},
B_2 = B_2 = B^2_{\lambda,\mu,\theta}$, Depth$(B_t) \le \mu$ and Depth$(B_1
\underset A {}\to \oplus B_1) \ge \lambda$.  We then start to investigate for
a fixed $A$, the existence of such $B_1,B_2$; gives sufficient and 
necessary conditions, involving consistency results.]
\endroster
\bn
\S2 $\quad$  On the family of homomorphic images of a Boolean Algebra
\mr
\item "{{}}"  [We prove that e.g. if $B$ is a Boolean Algebra of cardinality
$\lambda,\lambda \ge \mu$ and $\lambda,\mu$ are strong limit singular of the
same cofinality, then $B$ has a homomorphic image of cardinality $\mu$
(and with $\mu$ ultrafilters). \nl
More generally if $\lambda \ge \mu > \text{ cf}(\mu) = \text{ cf}(\lambda)$
and $B$ a Boolean Algebra of cardinality $\lambda$, then for some homomorphic
image $B'$ of $B$ we have $\mu \le |B'| \le 2^{< \mu}$]
\endroster 
\bn
\S3 $\quad$ If $d(B)$ is small, then depth or independence are not small
\mr
\item "{{}}"  [We prove for a Boolean Algebra $B$, \ub{that} if $d(B)^\kappa <
|B|$, then ind$(B) > \kappa$ or Depth$(B) \ge \text{ log}(|B|)$.]
\endroster
\bn
\S4 $\quad$  On omitting cardinals by compact spaces
\mr
\item "{{}}"   $\mu \le |B'|$
\endroster
\bn
\S5 $\quad$ Depth of ultraproducts of Boolean Algebras
\mr
\item "{{}}"  [We show that if $\square_\lambda$ and $\lambda =
\lambda^{\aleph_0}$ then for some Boolean Algebras $B_n$, Depth$(B_n) \le
\lambda$ but for any uniform ultrafilter $D$ on $\omega,\dsize \prod_{n <
\omega} B_n/D$ has depth $\ge \lambda^+$.]
\endroster
\newpage

\head {\S0} \endhead  \resetall
\bn
We consider some problems \nl
\ub{\stag{0.A} Problem}:  Is there a class of cardinals $\lambda$ (or just two)
such that there is a $\lambda^+$-thin tall superatomic Boolean Algebra $B$
(i.e. $|B| = \lambda^+,B$ is superatomic and for every $\alpha < \lambda^+,B$
has $\le \lambda$ atoms of order $\alpha$), or even a $\lambda^+$-tree;
provably in ZFC.
\smallskip

Note that if $\lambda = \lambda^{< \lambda}$ the answer is yes, so for
$\lambda = \aleph_0$ there is one.  Also note that if there is a
$\lambda^+$-tree, then there is such $\lambda^+$-thin tall superatomic
Boolean Algebra.  The point is that for several problems in Monk \cite{M},
72, 74, 75 and ZFC version of 73, 77, 78, 79 (all solved in \cite{RoSh:599}
(in the original version asking for consistency) 
there is no point to try to get positive answers as long as we do not
know it for \scite{0.A}. \nl
Also for several problems of \cite{M}(49,57,58,61,63,87) there is no point 
to try to get consistency of 
non-existence as long as we have not proved the consistency of the GSH
(generalized Souslin hypothesis) which says there is no $\lambda^+$-Souslin
trees or there is no $\lambda$-Souslin tree for $\lambda = \text{ cf}(\lambda)
 > \aleph_0$ for some others this is not provable, but it seems that this is
very advisable.
\bn
\ub{\stag{0.B} Problem}:  1) [M, Problem 28].  Is there a class of 
(or just one) $\lambda$ such that some Boolean Algebra $B$ of cardinality 
$\lambda^+$ has irredundancy $\lambda^+$. \nl
2) On irr: can we build a Boolean Algebra $B$,
irr$(B) < |B| = \lambda^{+n},n$ large enough?

Note that there is no point to try to construct examples as in problems [M]
(25,26,65,66,67,73,82,83,85,89) of Monk before we construct one for 
\scite{0.B}(1).
\bn
\ub{\stag{0.C} Question}:  1) Is it consistent to have a 
Boolean Algebra $B$ such that \nl
$|B| \ge \text{ irr}(B)^{++}$? \nl
2)  More general for cardinal invariants
with small difference with $|B|$ that is \nl
$2^{\text{inv}(B)} \ge |B|$, we should ask e.g. if 
$|B| < \text{ inv}(B)^{+ \omega}$. \nl
3) Similarly for irr$_n(B)$.
\bn
\ub{\stag{0.D} Question}:  Ultraproduct of length/Depth not near singular. \nl
See \S1.
\bn
\ub{\stag{0.E} Question}: Investigate $SpDpFP(A)$ for 
a Boolean Algebra $A$ (see in \S1).
\bn
\ub{\stag{0.F} Problem}:  It is true that for any large enough Boolean Algebra 
$B$ we have
\mr
\item "{{}}"  id$(B) = \text{ id}(B)^{< \theta}$ when e.g. $\theta =
\log_2|B|$, or at least for some constant $n,\theta = \text{ Min}\{\mu:
\beth_n(\mu) \ge |B|\}$.
\ermn
Similarly for compact spaces. \nl
By \cite{Sh:233}, for every $B$ there is such $n$.  If you like to try 
consistency, you have to use the phenomena proved consistent in Gitik Shelah
\cite{GiSh:344} (a problem of Hajnal).
\bn
\ub{\stag{0.G} Question}:  1) For which pairs $(\lambda,\theta)$ of cardinals
$\lambda \ge \theta$ there is a
superatomic Boolean Algebra with $> \lambda$ elements, $\lambda$ atoms, and
every $f \in \text{ Aut}(B) \text{ mod } < \theta$ atoms? 
(That is $|\{x:B \models ``x \text{ an atom and } f(x) \ne x\}| < \theta$).
\nl
2) In particular, is it true that for some $\theta$ for a class of $\lambda$
there is such Boolean Algebra? \nl
3) Replace ``automorphism" by ``1-to-1 endomorphism". \nl
4) In particular in \cite[\S2]{Sh:641} for $\mu$ strong limit singular. \nl
See \cite[\S2]{Sh:641}. \nl
\cite[\S5]{Sh:641} larger difference? \nl
Try: with $n$ depending on arity of the term as in \cite[\S6]{Sh:641}.
\bn
Concerning attainment in ZFC: \nl
\ub{\stag{0.H} Question}:  1) Can we show the distinction made between 
the attainments of variant of $hL$ (and $hd$), in semi-ZFC way?  That is, 
in Roslanowski Shelah \cite{RoSh:599} such examples are forced.  
Can we prove such examples exist adding to ZFC only restrictions on 
cardinal arithmetic? \nl
2) Similarly for other consistency results.  
(Well, preferably of low consistency strength).
\bn
\ub{\stag{0.I} Question}:  Let $\theta \ge \aleph_0$ or any cardinal.  
Is there a class of cardinal $\lambda = \lambda^\theta$ such that there 
is an entangled linear order of cardinality $\lambda^+$ (see \cite[AP2]{Sh:e}).

If we omit $\lambda = \lambda^+$, the answer is yes, and we ``almost can
prove the answer is yes (in ZFC)", meaning that if the answer is no in $V$ 
then on cardinal arithmetic there are very severe limitations: on the one 
hand for no $\lambda,2^\lambda = \lambda^+ \and \lambda = \lambda^\theta$ 
(and more) but on the other hand $(\neg(\exists \mu)(\mu^\theta = 
\mu^\theta \and 2^\mu \ge \aleph_{(\mu^{+4})}))$ and more, see 
\cite[\S5]{Sh:462}. \nl
This is closed connected to
\mr
\item "{$(*)$}"  can we have Inc$(\dsize \prod_{i < \theta} B_i/F) >
\dsize \prod_{i < \theta} \text{ Inc}(B_i)/F$ where $F$ is an ultrafilter on
$\kappa$.
\ermn
A ``yes answer" for the question gives yes to $(*)$.
\bn
This paper can be translated to compact topologies.
\newpage

\head {\S1 On the depth of free products} \endhead  \resetall
\bigskip

Monk \cite{M}, Problem 10, 11 ask about the depth of $B
{\underset A {}\to \oplus} C$ (see there for the known results, \cite{M}).  
\nl
We shall define a spectrum $SpDpFP(A)$ for a Boolean Algebra $A$ (in
\scite{1.0}) and phrase Monk's question with it (\scite{1.0A}(2),(3)).  We
then phrase a combinatorial statement $\boxtimes_{\lambda,\mu,\theta}$ and
prove it gives examples of $B,C$ extending $A$ while 
$B {\underset A {}\to \oplus} C$ has depth larger than both (in \scite{1.1}),
and note that it (provably in ZFC) holds for many cardinals (with $\lambda =
\mu$ near a singular) (in \scite{1.2}).  Later we note some variants of
$\boxtimes_{\lambda,\mu,\theta}$ and investigate when the construction in
\scite{1.1} works for a Boolean Algebra $A$, in particular for any infinite
Boolean Algebra $A$, it holds for a class of $\lambda$'s.
\bigskip

\definition{\stag{1.0} Definition}  1) For a Boolean Algebra $A$, we define the
spectrum of depth of free products over $A, SpDpFP(A)$ as

$$
\align
\bigl\{\kappa:&\text{ there are Boolean Algebras } B,C \text{ extending } A
\text{ such that}: \\
  &\,\text{Depth}(B {\underset A {}\to \oplus} C) > \kappa \ge
\text{ Depth}(B) + \text{ Depth}(A) \bigr\}.
\endalign
$$
\mn
2) Similarly

$$
\align
SpDpFP^+(A) = \bigl\{\kappa:&\text{ there are Boolean Algebras } B,C 
\text{ extending } A \\
  &\,\text{ such that Depth}(B {\underset A {}\to \oplus} C) \ge \kappa >
\text{ Depth}(B) + \text{ Depth}(A) \bigr\}.
\endalign
$$
\enddefinition
\bigskip

\remark{\stag{1.0A} Remark}  1) Note that

$$
\kappa^+ \in SpDpFP^+(A) \text{ iff } \kappa \in SpDpFP(A)
$$
\mn
so we can deal with \scite{1.0}(2) only. \nl
2) So written in our terms, problem 10 of Monk \cite{M} is:
\mr
\item "{$(*)$}"  for every infinite Boolean Algebra 
$A,SpDpFP(A)$ is a set of cardinals, i.e. has an upper bound.
\ermn
3) Written in our terms, problem 11 of Monk is
\mr
\item "{$(*)$}"  for every infinite Boolean Algebra $A,SpDpFP^+(A)$ 
is non-empty.
\ermn
4) By \scite{1.1}, \scite{1.2} (see $(*)_3$) below, e.g. for some countable
Boolean Algebra $A$, for every strong limit cardinal $\mu$ of cofinality
$\aleph_0$, we have $\mu^+ \in SpDpFP^+(A)$ (hence $\mu \in SpDPFP(A)$), so
Monk's question 10 is answered. \nl
5) The combinatorial property in $\boxtimes_{\lambda,\mu,\theta}$ is close to
one considered for investigating the ``bad stationary set of a successor of
singulars," and more generally the ideal $I[\lambda]$ in \cite{Sh:108},
\cite{Sh:88a}, \cite[\S1]{Sh:420}.
\mn
We then may ask ourselves: \nl
\ub{\stag{1.0B} Question}:  What occurs to cardinals which are 
not ``near singular"
(e.g. $\lambda = \mu = \chi^+,\chi = \chi^{< \chi} > 2^{|A|}$).
\mn
\ub{\stag{1.0C} Question}:  Can you say more on $SpDpFP^+(A)$ when we are given
$A$?
\sn
We give some information concerning those problems. 
\endremark
\bigskip

\proclaim{\stag{1.1} Claim}  Assume
\mr
\item "{$\boxtimes_{\lambda,\mu,\theta}(a)$}"   $\theta < \text{{\rm cf\/}}
(\mu) \le \mu \le \lambda$
\sn
\item "{$(b)$}"  $\bold c:[\lambda]^2 \rightarrow \theta$ satisfies: if
$\zeta_1 < \zeta_2 < \zeta_3 < \lambda$, then \nl
$\bold c \{\zeta_1,\zeta_3\} \le 
\text{ Max}\{\bold c\{\zeta_1,\zeta_2\},\bold c\{\zeta_2,\zeta_3\}\}$,
\sn
\item "{$(c)$}"  if $n < \omega,w_\alpha \in [\lambda]^n$ for 
$\alpha < \mu$, \ub{then} we can find $\alpha,\beta,i,j$ such that \nl
$\alpha < \beta < \mu,w_\alpha =
\{\zeta_\ell:\ell < n\}$ increasing, $w_\beta = \{\xi_\ell:\ell < n\}$
increasing, $\zeta_\ell = \xi_k \Rightarrow \ell = k,\bold c
\{\zeta_\ell,\zeta_k\} = \bold c\{\xi_\ell,\xi_k\}$ for $\ell < k < n$ 
and for some \nl
$i < j$ satisfy $i \ge \sup\{ \bold c\{\zeta_\ell,\zeta_k\},
\bold c\{\xi_\ell,\xi_k\}:\ell,k < n\}$ we have: for
$\ell,k < n$ one of the following occurs: 
$\bold c\{\zeta_\ell,\xi_k\} \ge j \and \bold c\{\zeta_k,\xi_\ell\} \ge j$ 
or $\bold c\{\zeta_\ell,\xi_k\} = 
\bold c\{\zeta_k,\xi_\ell\} < j \and [\zeta_\ell < \xi_k \leftrightarrow 
\xi_\ell < \zeta_k]$
\sn
\item "{$(d)$}"  $\theta = \text{ cf}(\theta) = \sup(\text{Rang}(\bold c))$
\sn
\item "{$(e)$}"  in clause (c) we can demand $i \ge i(*)$ for any pregiven
$i(*) < \theta$.
\ermn
\ub{Then} for any $\kappa \in [\theta,\mu)$ we can find Boolean Algebras
$A,B_1,B_2$ such that:
\mr
\item "{$(\alpha)$}"  $|A| = \kappa,A$ depend just on $\theta$ and $\kappa$
\sn
\item "{$(\beta)$}"  $|B_1| = |B_2| = \lambda$
\sn
\item "{$(\gamma)$}"  Depth$^+(B_1 \underset A {}\to \oplus \,B_2) 
= \lambda^+$
\sn
\item "{$(\delta)$}"  Depth$^+(B_1) \le \mu$ and Depth$^+(B_2) \le \mu$
\sn
\item "{$(\delta)^+$}"  moreover, Length$^+(B_i) \le \mu$ for $i = 1,2$.
\endroster
\endproclaim
\bigskip

\remark{\stag{1.1A} Remark}  1) Let $\boxtimes^-_{\lambda,\mu,\theta}$ 
just means (a), (b), (c), (e), i.e. omitting clause (d).  So by \scite{1.5}(3)
below $\boxtimes^-_{\lambda,\mu,\theta}
\Rightarrow \dsize \bigvee_{\sigma \le \theta} \boxtimes_{\lambda,\mu,\sigma}$
and there are obvious monotonicity properties. \nl
2) Really we can omit $i$ in clause (c) of $\boxtimes_{\lambda,\mu,\theta}$, 
see proof. \nl
3) Let $\boxtimes^+_{\lambda,\mu,\theta}$ means that in the definition of
$\boxtimes_{\lambda,\mu,\theta}$ in clause (c), when $\bold c\{\zeta_\ell,
\xi_k\} \ge j \and \bold c\{\zeta_k,\xi_\ell\} \ge j$ we add
$\bold c\{\zeta_\ell,\xi_k\} = \bold c\{\zeta_k,\xi_\ell\}$. \nl
In \scite{1.2} below we can get this version by working a little more.
\endremark
\bn
We quote \cite{Sh:g}
\demo{\stag{1.2} Observation}  The demand on $\theta < \mu \le \lambda$ is not
hard, in fact we can find $\bold c$ such that clauses (b), (c), (d), (e) of
$\boxtimes_{\lambda,\mu,\theta}$ hold if the cardinals $\theta,\mu,\lambda$
satisfies at least one of the following statements:
\mr
\item "{$(*)_1$}"  for some $\chi \in (\theta,\mu)$ we have 
$\dsize \bigwedge_{\alpha < \chi} |\alpha|^{< \theta} < \chi$, 
cf$(\chi) = \theta$, pp$^+_{J^{bd}_\theta}(\chi) > \lambda$
\ermn
or
\mr
\item "{$(*)_2$}"  $\theta > \aleph_0,\theta = \text{ cf}(\chi),\dsize
\bigwedge_{\alpha < \chi} |\alpha|^\theta < \chi,\chi < \mu \le \lambda \le
\chi^\theta$
\sn
\item "{$(*)_3$}"  $\lambda = \text{ cf}(\lambda) \ge \mu > \chi >
\text{ cf}(\chi) = \theta,\langle \chi_i:i < \theta \rangle$ is strictly
increasing with limit $\chi$, max pcf$\{\chi_j:j < i\} < \chi$ (or just
$< \mu$) and $\dsize \prod_{i < \theta} \chi_i/J^{bd}_\theta$ has true
cofinality $= \lambda$ (or just there is a $<_{J^{bd}}$-increasing sequence
in the product for every regular $\lambda' \in (\chi,\lambda]$)
\sn
\item "{$(*)_4$}"  if $\aleph_0 < \text{ cf}(\sigma) = \sigma =
\text{ cf}(\chi) < \chi$, and $(\forall \alpha < \chi)(\text{cf}([\alpha]
^{\le \sigma},\subseteq) < \chi)$ and $\sigma^{\aleph_0} < \chi$, \ub{then}
for some club $C$ of $\chi,\chi_1 \in \{\chi\} \cup C \and \lambda = 
\text{ cf}([\chi_1]^{\le \sigma},\subseteq) \and \mu = \chi^+_1 \and \theta
= \text{ cf}(\chi_1) \Rightarrow \boxtimes_{\lambda,\mu,\theta}$.
\endroster
\enddemo
\bigskip

\demo{Proof}  \ub{For $(*)_1$}  Let $\chi = \dsize \sum_{i < \theta} \chi_i$,
with $\chi_i = \text{ cf}(\chi_i) < \chi,\theta = \text{ cf}(\chi)$ and in 
$\dsize \prod_{i < \theta} \chi_i$ we can find a 
$<_{J^{bd}_\theta}$-increasing sequence $\langle \eta_\alpha:\alpha < 
\lambda \rangle$ such that $|\{\eta_\alpha \restriction i:\alpha < \lambda\}|
\le \chi$ (really $< \mu$ suffice: if $\lambda$ is regular; by the definition
of $J^{bd}_\theta$; if $\lambda$ is singular by combining such examples for
regular (and \cite[Ch.II,3.5]{Sh:g}). \nl
For $\alpha < \beta$ let $\bold c\{\alpha,\beta\} = \text{ Min}\{i < \theta:
(\forall j \in [i,\theta))(\eta_\alpha(i) < \eta_\beta(i)\}$.
\enddemo
\bn
Let us check 
\mn
\ub{Clause $(a)$}:  By the assumptions on $\theta,\mu,\lambda$ clearly
$\theta < \chi < \mu \le \lambda$; note that we do not demand 
$\theta < \text{ cf}(\mu)$.
\bn
\ub{Clause $(b)$}:  Clearly $\bold c$ is a function from $[\lambda]^2$ to
$\theta$.  Now suppose $\alpha < \beta < \gamma < \lambda$ and $i = \bold c
\{\alpha,\beta\}$ and $j = \bold c\{\beta,\gamma\}$ hence: max$\{i,j\} \le
\zeta < \theta \Rightarrow \eta_\alpha(\zeta) < \eta_\beta(\zeta) < 
\eta_\gamma(\zeta)$ hence $\bold c\{\alpha,\gamma\} \le \text{ max}\{i,j\}$
as required.
\bn
\ub{Clause $(c),(e)$}:  So suppose $n < \omega$ and $w_\alpha \in
[\lambda]^n$ for $\alpha < \mu$.  Let $w_\alpha = \{\zeta_{\alpha,\ell}:\ell
< n_\alpha\}$ with $\zeta_{\alpha,\ell} < \zeta_{\alpha,\ell +1}$ and $i(*) < 
\theta$. \nl
Also \wilog \, for some $v \subseteq n$: (just shrink the set for each
$\ell < n$)
\mr
\item "{$\otimes_1$}"  if $\ell \in v$ then 
$\langle \zeta_{\alpha,\ell}:\alpha < \chi^+ \rangle$ is strictly 
increasing with limit $\zeta^*_\ell$
\sn
\item "{$\otimes_2$}"  if $\ell < n,\ell \notin v$ then 
$\langle \zeta_{\alpha,\ell}:\alpha < \chi^+ \rangle$ is constantly 
$\zeta^*_\ell$.
\ermn
We can also demand
\mr
\item "{$\otimes_3$}"  if $\ell \in v,m < n$ the truth value of
$\zeta_{\alpha,\ell} > \zeta^*_m,\zeta_{\alpha,\ell} < \zeta^*_m$ does \nl
not depend on $\alpha$
\sn
\item "{$\otimes_4$}"  if $\ell \ne m$ are from $v$ and $\zeta^*_\ell =
\zeta^*_m$ then $\alpha < \beta < \chi^+ \Rightarrow \zeta_{\alpha,\ell} <
\zeta_{\beta,m}$.
\ermn
For each $\alpha < \mu^{\chi^+}$ let 
$i_\alpha = \text{ Max}[\{\bold c\{\zeta,\xi\}:
\zeta \ne \xi \in w_\alpha \cup \{\zeta^*_\ell:\ell < n\}\} \cup \{i(*)\}]$,
now clearly $i_\alpha < \theta$ so as cf$(\chi^+) > \theta,\chi^+ \le \mu$
(by clause (a) or more exactly by our assumptions) \wilog \, $i_\alpha = 
i^*$ for $\alpha < \chi^+$.  As
$\langle \eta_{\zeta_{\alpha,\ell}}(i^* +1):\ell < n \rangle$ have only
$\le (\chi_{i^*+1})^n < \chi$ possible values, without loss of generality
$\alpha < \chi^+ \Rightarrow \eta_{\zeta_{\alpha,\ell}}(i^*+1) = \gamma_\ell$,
moreover $\alpha < \chi^+ \Rightarrow \eta_{\zeta_{\alpha,\ell}} \restriction
(i^* +2) = \nu_\ell$.
\mn
Also \wilog
\mr
\item "{$\otimes_5$}"  for every $\alpha < \chi^+$ and finite $u \subseteq
\theta$, for $\chi^+$ ordinals $\beta < \chi^+$ we have
$\dsize \bigwedge_{i \in u} \dsize \bigwedge_{\ell < n} f_{\zeta_{\beta,\ell}}
(i) = f_{\zeta_{\alpha,\ell}}(i)$.
\ermn
Let $\alpha_1 < \alpha_2 < \chi^+$, so we can find $j^* \in (i(*)+2,\theta)$ 
such that $\{\ell,m\} \subseteq v \and \zeta^*_\ell = \zeta^*_m 
\Rightarrow f_{\zeta_{\alpha_1,\ell}}(j^*) <
f_{\zeta_{\alpha_2,\ell}}(j^*)$.  Choose $\alpha_3 \in (\alpha_2,\chi^+)$ such
that $\dsize \bigwedge_\ell f_{\zeta_{\alpha_3,\ell}}(j^*) =
f_{\zeta_{\alpha_1,\ell}}(j^*)$.  Now let $\alpha =: \alpha_2,
\beta =: \alpha_3,i = i^*,j = i^* +1$, they are as required. \nl
Why?  
\mr
\item "{$\otimes_6$}"  assume $\ell \ne m$ are $< n$ and $\zeta^*_\ell \ne
\zeta^*_m$ then we have: $\bold c\{\zeta_{\alpha,\ell},\zeta_{\beta,m}\} =
\bold c\{\zeta_{\beta,m},\zeta_{\alpha,\ell}\}$ is determined by
$\nu_\ell,\nu_m$ and is $\le i^* = i$.  Also $\zeta_{\alpha,\ell} <
\zeta_{\beta,m} \leftrightarrow \zeta_{\beta,\ell} < \zeta_{\alpha,m}]$. \nl
[Why?  E.g. if $\zeta^*_\ell < \zeta^*_m$ for $i \in [i^*,\theta)$ we have
$f_{\zeta_{\alpha,\ell}}(i) \le f_{\zeta^*_\ell}(i) < f_{\zeta_{\beta,m}}(i)$,
so $\bold c\{\zeta_{\alpha,\ell},\zeta_{\beta,m}\} \le i^*$, noticing  
$i_\alpha = i_\beta \le i^*$, (and the choice of $i_\alpha$ and $\otimes_5$ 
above).]
\sn
\item "{$\otimes_7$}"  if $\zeta_{\alpha,\ell} \ne \zeta_{\beta,\ell}$ (i.e.
$\ell \in v$) then $\bold c(\zeta_{\alpha,\ell},\zeta_{\beta,\ell}) \ge j^*$
\nl
[why?  by the choice of $j^*$ we have $f_{\beta,\ell}(j^*) =
f_{\alpha_1,\ell}(j^*) < f_{\alpha,\ell}(j^*)$ and by the choice of $\beta$
we have $\alpha < \beta$.]
\sn
\item "{$\otimes_8$}"  if $\ell \ne m < n,\zeta^*_\ell = \zeta^*_n$ 
and $\{\ell,m\} \subseteq v$ then
$\bold c(\zeta_{\alpha,\ell},\zeta_{\beta,m}) \ge j^* > j$ \nl
[why?  clearly $\zeta_{\alpha,\ell} < \zeta_{\beta,m}$ as $\alpha < \beta$ 
(see $\otimes_4$) and \nl
$\eta_{\zeta_{\beta,m}}(j^*) = \eta_{\zeta_{\alpha_3,m}}
(j^*) = \eta_{\zeta_{\alpha_1,m}}(j^*) < \eta_{\zeta_{\alpha_2,\ell}}(j^*) =
\eta_{\zeta_{\beta,\ell}}(j^*)$.]
\sn
\item "{$\otimes_9$}"  if $\ell \ne m < n,\zeta^*_\ell = \zeta^*_m$ and
$\{\ell,m\} \nsubseteq v$, then $\bold c\{\zeta_{\alpha,\ell},
\zeta_{\beta,\ell}\} \le i^*$.  Also $[\zeta_{\alpha,\ell} < \zeta_{\beta,m}
\leftrightarrow \zeta_{\beta,\ell} < \zeta_{\alpha,m}]$. \nl
[why?   clearly $\{\ell,m\} \cap v = \emptyset$ is impossible, and so
$|\{\ell,m\} \cap v| = 1$; now the proof is similar to that of $\otimes_6$.]
\ermn
Together we are done.
\bn
\ub{For $(*)_2$}:

By \cite[Ch.II,5.4(2),Ch.VIII,\S1]{Sh:g} the assumptions in $(*)_1$ holds.
\bn
\ub{For $(*)_3$}:  

By \cite[Ch.II,3.5]{Sh:g} we can choose $\bar \eta = \langle \eta_\alpha:
\alpha < \lambda \rangle,\eta_\alpha \in {}^\theta \chi,\bar \eta$ is
$<_{J^{bd}_\theta}$-increasing, and $|\{\eta_\alpha \restriction i:\alpha <
\lambda\}| \le \chi$ for $i < \theta$  (again for singulars we combine
examples). \nl
For $\lambda$ singular combine such examples.
\bn
\ub{For $(*)_4$}:  

See \cite{Sh:E11}.
\bigskip

\demo{Proof of \scite{1.1}} 
\sn
\ub{Stage A}:  Let $A$ be the Boolean Algebra generated by 
$\{a^t_i:i < \kappa \text{ and } t \in \{1,2\}\}$ freely except the
equations:
\mr
\item "{$(*)_1$}"  $a^t_i \le a^t_j$ for $i \le j < \theta$ and
$t \in \{1,2\}$
\sn
\item "{$(*)_2$}"  $a^1_i \cap a^2_j = 0$ for $i,j < \theta$.
\ermn
Let $I^t$ be the ideal of $A$ generated by $\{a^t_i:i < \theta\}$ for 
$t = 1,2$.  Let $I$ be the ideal of $A$ generated by $I^1 \cup I^2$ so 
$A/I$ is the trivial (= two elements) Boolean Algebra. 
\nl
For $t = 1,2$ let $B_t$ be the extension of $A$ by $\{x^t_\alpha:\alpha <
\lambda\}$ freely except that:
\mr
\item "{$(*)_3$}"  $x^t_\alpha - x^t_\beta \le 
a^t_{\bold c\{\alpha,\beta\}} \qquad$ for $\alpha < \beta < \lambda$
\sn
\item "{$(*)_4$}"  $x^t_\alpha \cap a^{3-t}_i = 0 \qquad \qquad$ for $\alpha <
\lambda,i < \theta$.
\ermn
Let $B = B_1 \underset A {}\to \oplus B_2$ and let $x_\alpha = x^1_\alpha
\cap x^2_\alpha \in B$.  Let $J^1,J^2,J$ be the ideal of $B$ which $I^1,
I^2,I$ generates (resp.), \nl
Clearly $|B| \le \lambda$.
\bn
\ub{Stage B}:  Depth$^+(B) = \lambda^+$.

Clearly Depth$^+(B) \le |B|^+ \le \lambda^+$.  Hence it suffices to prove 
$\langle x_\alpha:\alpha < \lambda \rangle$ is strictly
increasing.  So let $\alpha < \beta$.  First we prove $x_\alpha \le x_\beta$.
Let $D$ be an ultrafilter on $B$ (so $D_t = D \cap B_t \in \text{ Ult}(B_t)$
and $D_1 \cap A = D_2 \cap A$), and we shall prove $x_\alpha \in D
\Rightarrow x_\beta \in D$.  This suffices for proving $x_\alpha \le x_\beta$ 
(as $D$ was any ultrafilter).
\enddemo
\bn
\ub{Case 1}:  $D \cap A$ is disjoint to $I$.

Now modulo $J^t$ in $B^t,x^t_\alpha \le x^t_\beta$ hence modulo $J$,
$x^t_\alpha \le x^t_\beta$ hence mod $J,x_\alpha \le x_\beta$ hence 
$x_\alpha \in D \Rightarrow x_\beta \in D$.
\bn
\ub{Case 2}:  $D \cap A$ is not disjoint to $I$.

So $D \cap I \ne \emptyset$ hence for some $t \in \{1,2\}$ and $i < \theta,
a^t_i \in D$, but $x^{3-t}_\alpha$ is disjoint to $a^t_i$ by $(*)_4$, so
$x^{3-t}_\alpha \notin D$ hence $x_\alpha = x^1_\alpha \cap x^2_\alpha 
\notin D$ so trivially $x_\alpha \in D \Rightarrow x_\beta \in D$.
\mn

We still have to prove $x_\alpha \ne x_\beta$.  Let $D^t_\beta$ be the
ultrafilter on $B^t$ generated by $x^t_\gamma \, (\gamma \ge \beta),
-x^t_\gamma \, (\gamma < \beta),-a^s_i \, (i < \theta,s \in \{1,2\})$ (check
that it is okay trivially).  Now $D^1_\beta \cap A = I = D^2_\beta \cap A$ 
(check, trivial).  So there is an ultrafilter $D$ on $B,D \cap B_t = 
D^t_\beta$.
So $-x^t_\alpha,x^t_\beta \in D^t_\beta$ so $-x_\alpha,x_\beta \in D$ as 
required.
\bn
\ub{Stage C}:  Length$^+(B_t) \le \mu$.

Assume not, so we can find $\langle c_\alpha:\alpha < \mu \rangle$ a chain
(so with no repetition).  Let \nl
$c_\alpha = \tau_\alpha(a^{s_{\alpha,1}}_{i_{\alpha,1}},\dotsc,
a^{s_{\alpha,k_\alpha}}_{i_{\alpha,k_\alpha}},k_\alpha,
x^t_{\zeta_{\alpha,1}} \ldots 
x^t_{\zeta_{\alpha,n_\alpha}})$ where $i_{\alpha,1} < \ldots <
i_{\alpha,k_\alpha} < \kappa$ and $\zeta_{\alpha,1} < \ldots < 
\zeta_{\alpha,n_\alpha} < \lambda$ and $s_{\alpha,1},\dotsc,
s_{\alpha,k_\alpha} \in \{1,2\}$ and
$\tau_\alpha$ a Boolean term. 
\sn
As cf$(\mu) > \theta \ge \aleph_0$, without loss of generality:
\mr
\item "{$(*)_5$}"  $k_\alpha = k^*,n_\alpha = n^*,\tau_\alpha = \tau^*$,
$s_{\alpha,\ell} = s_\ell$.
\sn
\item "{$(*)_6$}"  for some $m^* \le k^*$ we have:
{\roster
\itemitem{ $(i)$ }  $0 < \ell \le m^* \Rightarrow i_{\alpha,\ell} = 
i_\ell < \theta$
\sn
\itemitem{ $(ii)$ }  $\ell \in (m^*,k^*] \Rightarrow i_{\alpha,\ell} \ge
\theta$
\sn
\itemitem{ $(iii)$ }  for $\alpha < \beta < \mu$, for some $v =
v_{\alpha,\beta} \subseteq (m^*,k^*]$ we have \nl
\sn

$\qquad (a) \quad \ell \in v \Rightarrow i_{\alpha,\ell} =
i_{\beta,\ell}$ and \nl

$\qquad (b) \quad \ell \in (m^*,k^*] \and \ell \notin v$ implies
$i_{\alpha,\ell} \ne i_{\beta,\ell}$ \nl

$\qquad \qquad$ and are not in $\{i_{\alpha,k},i_{\beta,k}:k \ne \ell\}$.
\endroster}
\sn
[why $(*)_2$?  if $\mu$ is regular, by the $\triangle$-system lemma (so then
$w_{\alpha,\beta} = w$) and if singular, apply it twice.]
\ermn
Let $i(*) = \sup\{i_\ell:\ell < m^*\}$ so $i(*) < \theta$. 
\mn 
Let $w_\alpha = \{\zeta_{\alpha,1},\dotsc,\zeta_{\alpha,n^*}\}$.  Let
$\alpha \ne \beta < \mu$ and $i < j < \theta$ be as guaranteed by clause (c) 
of the assumption and $i > i(*)$ (see clause (e)).  
So $c_\alpha,c_\beta$ are distinct members of a chain of $B_t$. \nl
Now read Stage D below.

So let $w =: \{\zeta_{\alpha,\ell},\zeta_{\beta,\ell}:\ell = 1,\dotsc,n\}
\in [\lambda]^{< \aleph_0}$.
\sn
By Stage D, $B_{t,w}$ is a subalgebra of $B_t$ and $c_\alpha,c_\beta \in
B_{t,w}$, hence also in $B_{t,w},c_\alpha,c_\beta$ are distinct members of a
chain.  

By symmetry assume $B_{t,w} \models ``c_\alpha < c_\beta"$, hence there is
a homomorphism $f$ from $B_{t,w}$ to (the trivial Boolean algebra) $\{0,1\}$
such that $f(c_\alpha) = 0,f(c_\beta) = 1$.  Let $\gamma(f) = \text{ Min}
\{\gamma \le \kappa:\text{if } \gamma < \theta \text{ then } f(a^t_\gamma) 
=1\}$.
\nl
As in the Boolean Algebra $A \subseteq B^t,\langle a^t_\gamma:\gamma <
\theta \rangle$ is increasing, clearly
\mr
\item "{$(*)_7$}"  for $\xi < \theta,f(a^t_\xi) =1 \Leftrightarrow \xi \ge
\gamma(f)$
\sn
\item "{$(*)_8$}"  for $\xi < \theta,f(a^{3-t}_\xi) = 0$ \nl
[why?  otherwise $\alpha < \lambda \Rightarrow f(x^t_\alpha) =0$ by $(*)_4$,
hence $f(x^t_{\zeta_{\alpha,\ell}}) = 0 = f(x^t_{\zeta_{\beta,\ell}})$.  
Now let $g:\{a^1_i,a^2_i:i < \kappa\} \cup \{x^t_\alpha:\alpha < \lambda\} 
\rightarrow B_t$ be defined as follows: 
$g(a^{s_\ell}_{i_{\alpha,\ell}}) = a^{s_\ell}_{i_{\beta,\ell}},
g(a^{s_\ell}_{i_{\beta,\ell}}) = a^{s_\ell}_{i_{\alpha,\ell}}$ and otherwise
it is the identity so $i < \theta \and s \in \{1,2\} \Rightarrow
g(a^s_i) = a^s_i$.  By the assumption toward contradiction and $(*)$ of
Stage D, $g$ induces a homomorphism $\hat g$ from
$B_t$ to $B_t$, clearly it is an automorphism, so $f \circ \hat g$ is a
homomorphism from $B_t$ to $\{0,1\}$ and:
$(f \circ \hat g)(c_\alpha) = f(\hat g(c_\alpha)) = f(c_\beta) = 1,
(f \circ \hat g)(c_\beta) = f(\hat g(c_\beta)) = f(c_\beta) = 0$ contradicting
$B_t \models c_\alpha < c_\beta$]
\ermn
Now
\mr
\item "{$(*)_9$}"  the function $g:\{a^s_i:i < \kappa,s \in \{1,2\}\} 
\rightarrow \{0,1\}$ induce a homomorphism $\hat g$ from $A$ to $\{0,1\}$ 
where $g$ is defined by:
\sn
{\roster
\itemitem{ $(i)$ }  $g(a^t_\xi) = f(a^t_\xi)$ for $\xi < j$
\sn
\itemitem{ $(ii)$ }  $g(a^t_\xi) = 1$ if $\xi \ge j,\xi < \theta$
\sn
\itemitem{ $(iii)$ }  $g(a^{3-t}_\xi) = f(a^{3-t}_\xi) = 0$ for $\xi < \theta$
\sn
\itemitem{ $(iv)$ }  $g(a^{s_\ell}_{i_{\alpha,\ell}}) =
f(a^{s_\ell}_{i_{\beta,\ell}}),g(a^{s_\ell}_{i_{\beta,\ell}}) =
f(a^{s_\ell}_{i_{\alpha,\ell}})$ for $\ell = 1,\dotsc,k^*$ 
\sn
\itemitem{ $(v)$ }  $g(a^s_i) = f(a^s_i)$ if $a^s_i \notin \{
a^{s,\ell}_{i_{\alpha,\ell}},a^{s_\ell}_{i_{\beta,\ell}}:\ell = 1,\dotsc,
k^*\}$ and $i \in [\theta,\kappa)$ 
\endroster}
[why?  now $g$ is well defined as, e.g. for contradiction concerns (iii), two 
instances do not contradict by $(*)_6(iii)$ and they do not contradict 
others by $(*)_6(i)+(ii)$.  By $(*)$ of Stage D below we should check the 
equations appearing in the definition of $A$.  For those in $(*)_1$, i.e.
$a^s_\varepsilon \le a^s_\xi$ for $\varepsilon < \xi$, if $s = 3-t$
this is trivial by clause (iii), if $s=t,\xi \ge j$, this is trivial by
clause (ii) and if $s=t,\xi < j$, then $g(a^s_\varepsilon) = 
f(a^s_\varepsilon) \le f(a^s_\xi) = g(a^s_\xi)$.
\sn
As for the equations in $(*)_2$ that is $a^1_\varepsilon \cap a^2_\xi=0$ for
$\varepsilon,\xi < \theta$ they are preserved trivially by $(*)_7$ and
clause (iii).]
\ermn
Define a function $h$ from $A \cup \{x_{\zeta_{\alpha,\ell}},
x_{\zeta_{\beta,\ell}}:\ell = 1,\dotsc,n^*\}$ to $\{0,1\}$ as follows:
$h \restriction A$ is the homomorphism $\hat g$ to $\{0,1\}$.   Now define
$h(x_{\zeta_{\alpha,\ell}}) = f(x_{\zeta_{\beta,\ell}})$ and
$h(x_{\zeta_{\beta,\ell}}) = f(x_{\zeta_{\alpha,\ell}})$.  Now
\mr
\item "{$(*)_{10}$}"  $h$ induces a homomorphism $\hat h$ from
$B_{t,w}$ to $\{0,1\}$ \nl
[why?  by $(*)_6$ the function $h$ is well defined; we use $(*)$ of Stage D:  
now
{\roster
\itemitem{ $(a)$ }  the equations in $A$ are respected by the choice of
$h \restriction A$ as $g$ (and $g$ being a homomorphism)
\sn
\itemitem{ $(b)$ }  the equations $x^t_{\zeta_{\gamma,\ell}} \cap
a^{3-t}_i = 0$ for $\gamma \in \{\alpha,\beta\}$. \nl
This is respected as $h(a^{3-t}_i) = g(a^{3-t}_i) = 0$.
\sn
\itemitem{ $(c)$ }  for $\gamma \in \{\alpha,\beta\}$ the equation
$x^t_{\zeta_{\gamma,\ell}} - x^t_{\zeta_{\gamma,m}} \le 
a^t_{\bold c\{\zeta_{\gamma,\ell},\zeta_{\gamma,m}\}}$ for $\ell < m$.
\nl
Let $\delta$ be such that $\{\gamma,\delta\} = \{\alpha,\beta\}$ and remember
that $i \ge \bold c\{\zeta_{\gamma,\ell},\zeta_{\gamma,m}\} = \bold c
\{\zeta_{\delta,\ell},\zeta_{\delta,m}\}$ (by the choice of 
$\alpha,\beta,i,j$) so $h(a^t_{\bold c\{\zeta_{\gamma,\ell},
\zeta_{\gamma,m}\}}) = f(a^t_{\bold c\{\zeta_{\gamma,\ell},\zeta_{\gamma,m}
\}}) = f(a^t_{\bold c\{\zeta_{\delta,\ell},\zeta_{\gamma,m}\}})$ by the 
choice of $h \restriction A$ hence

$$
\align
h(x^t_{\zeta_{\gamma,\ell}}) - h(x^t_{\zeta_{\gamma,m}}) &=
f(x^t_{\zeta_{\delta,\ell}}) - f(x^t_{\zeta_{\delta,m}}) \\
  &= f(x^t_{\zeta_{\delta,\ell}} - x^t_{\zeta_{\delta,m}}) \le
f(a^t_{\bold c\{\zeta_{\delta,\ell},\zeta_{\delta,m}\}}) \\
  &= h(a^t_{\bold c\{\zeta_{\gamma,\ell},\zeta_{\gamma,m}\}})
\endalign
$$
\sn
\itemitem{ $(d)$ }  if $\ell,m \in \{1,\dotsc,m^*\}$ and $\gamma \ne \delta
\in \{\alpha,\beta\}$ and $\zeta_{\gamma,\ell} < \zeta_{\delta,m}$ the 
equation $x^t_{\zeta_{\gamma,\ell}} - x^t_{\zeta_{\delta,m}} \le
a^t_{\bold c\{\zeta_{\gamma,\ell},\zeta_{\delta,m}\}}$. \nl
Now we have to look at clause (c) of $\boxtimes_{\lambda,\mu,\theta}$ of
\scite{1.1}, there are two possibilities \nl
\block
\ub{possibility 1}:  $\bold c\{\zeta_{\gamma,\ell},\zeta_{\delta,m}\} \ge j$
(but necessarily $< \theta$), then \nl 
$g(a^t_{\bold c\{\zeta_{\gamma,\ell},\zeta_{\delta,m}\}}) = 1$ 
by clause (ii) of $(*)_9$ so the equation is respected. \nl
\sn
\ub{possibility 2}:  $\bold c\{\zeta_{\gamma,\ell},\zeta_{\delta,m}\} =
\bold c\{\zeta_{\gamma,\ell},\zeta_{\gamma,m}\} \le i$ hence \nl
$[\zeta_{\gamma,\ell} < \zeta_{\delta,m} \Leftrightarrow \zeta_{\delta,\ell}
< \zeta_{\gamma,m}]$ \nl
(read (c) of $\boxtimes_{\lambda,\mu,\theta}$ so both holds). \nl
\mn
So $h(x^t_{\zeta_{\gamma,\ell}}) - h(x^t_{\zeta_{\delta,m}}) =
f(x^t_{\zeta_{\delta,\ell}}) - f(x^t_{\zeta_{\gamma,m}}) =
f(x^t_{\zeta_{\delta,\ell}} - x^t_{\zeta_{\gamma,m}})$ \nl

$\qquad \qquad \qquad \le f(a^t_{\bold c
\{\zeta_{\gamma,\ell},\zeta_{\delta,m}\}}) =
g(a^t_{\bold c\{\zeta_{\gamma,\ell},\zeta_{\delta,m}\}}) =
h(a^t_{\bold c\{\zeta_{\gamma,\ell},\zeta_{\delta,m}\}})$.
So we have proved $(*)_{10}$.]
\endblock
\endroster}
\ermn
Now we have two homomorphisms $f,\hat h$ from $B_{t,w}$ to $\{0,1\}$ and
they satisfy:
\mr
\item "{$(*)_{11}(a)$}"  $f(a^{s_\ell}_{i_{\alpha,\ell}}) = \hat h
(a^{s_\ell}_{i_{\beta,\ell}})$
\sn
\item "{$(b)$}"  $f(a^{s_\ell}_{i_{\beta,\ell}}) = 
\hat h(a^{s_\ell}_{i_{\alpha,\ell}})$
\sn
\item "{$(c)$}"  $f(x^t_{\zeta_{\alpha,\ell}}) = \hat h
(x^t_{\zeta_{\beta,\ell}})$
\sn
\item "{$(d)$}"  $f(x^t_{\zeta_{\beta,\ell}}) = \hat h
(x^t_{\zeta_{\alpha,\ell}})$ \nl
[why?  for (a) + (b) not $\hat h \restriction A = g \restriction A$ and
$(*)_9(iv)$ for (c) + (d) just see the choice of $h$.]
\ermn
So clearly $f(c_\alpha) = \hat h(c_\beta),f(c_\beta) = \hat h(c_\alpha)$, but
$f(c_\alpha) = 0 < f(c_\beta) = 1$ so $\hat h(c_\beta) = 0 < \hat h
(c_\alpha) = 1$ whereas we assume $B_{w,t} \models c_\alpha < c_\beta$,
contradiction.

So we have finished the proof of Stage C, hence of \scite{1.1} except a debt:
Stage D.
\bn
\ub{Stage D}:  First recall
\mr
\item "{$(*)$}"  \ub{if} a Boolean algebra $B$ is defined by: generated freely
by $\{x_i:i < i^*\}$ except the set of equations $\Gamma$, and $B'$ another
Boolean Algebra and the function \nl
$h:\{x_i:i < i^*\} \rightarrow B'$ respect the equations
in $\Gamma$ \nl
(i.e. if $\tau'_i(x_{i_1},\ldots) = \tau''(x_{j_1},\ldots)
\in \Gamma$ then \nl
$B' \models \tau'(h(x_{i_1}),\ldots) = \tau''(h(x_{j_1}),\ldots)$) \nl
\ub{then} $h$ can be extended to a homomorphism from $B$ to $B'$ (and we
call it $\hat h$), similarly for ``extensions of a Boolean Algebra $A$".
\ermn
For $w \subseteq \lambda$ let $B_{t,w}$ be defined just like $B_t$
restricting ourselves to $\alpha \in w$ (so also in the set of equations we
consider involve only $\alpha,\beta \in w$). \nl
A priori it is not guaranteed that $w \subseteq u \subseteq \lambda
\Rightarrow B_{t,w} \subseteq B_{t,w}$.  Note $B_t = B_{t,\lambda}$.
\bigskip

\demo{Fact}  For $w \subseteq u \subseteq \lambda,B_{t,w} \subseteq B_{t,u}$
and $B_{t,u}$ is the direct limit of \nl
$\{B_{t,w}:w \subseteq u \text{ finite}\}$.
\enddemo

\demo{Proof}  It is enough to prove this for finite $u$, so we can ignore
the second phrase as it follows.  The first phrase we prove by induction 
on $|u \backslash w|$, so without loss of generality 
$|u \backslash w| = 1$, let $\zeta \in u \backslash w$. \nl
We define $h:A \cup \{x^t_i:i \in u\} \rightarrow B_{t,w}$,

$$
h \restriction A = \text{ identity}
$$

$$
h(x^t_i) = x^t_i \text{ if } i \in w
$$

$$
h(x^t_\zeta) = \cup\{x^t_\xi - a^t_{\bold c\{\zeta,\xi\}}:
\xi \in w \cap \zeta\}.
$$
\mn
Now $h$ is as in $(*)$ (see beginning of stage D, checked below) so there 
is a homomorphism from $B_{t,u}$
to $B_{t,w}$ which obviously extends the identity so we are done. \nl
Why is $h$ as required in $(*)$?  We check the ``new" equations, i.e.
the ones appearing in the definition of $B_{t,u}$ and not in the definition
of $B_{t,w}$ (which are satisfied as $h \restriction B_{t,w}$ is the
identity:
\mr
\widestnumber\item{$(iii)$}
\item "{$(i)$}"  $x^t_\zeta \cap a^{3-t}_j = 0$ \nl
[why?  obvious by the choice of $h(x^t_\zeta)$ as 
$h(x^t_\zeta) \cap h(a^{3-t}_j) = h(x^t_\zeta) \cap a^{3-t}_j \le
\bigl( \dbcu_{\xi \in w \cap \zeta} h(x^t_\xi) \bigr) \cap a^{3-t}_j =
\dbcu_{\xi \in w \cap \zeta} h(x^t_\xi) \cap a^{3-t}_j =
\dbcu_{\xi \in w \cap \zeta} 0 = 0$ so (i) holds]
\sn
\item "{$(ii)$}"  if $\varepsilon < \zeta,\varepsilon \in w$, then the
equation $x^t_\varepsilon - x^t_\zeta \le 
a^t_{\bold c\{\varepsilon,\zeta\}}$ \nl
[why?  by the choice of $h(x^t_\zeta)$ clearly $h(x^t_\varepsilon) -
h(x^t_\zeta) =$ \nl
$x^t_\varepsilon - h(x^t_\zeta) \le x^t_\varepsilon - 
(x^t_\varepsilon - a^t_{\bold c\{\zeta,\varepsilon\}}) \le
a^t_{\bold c\{\varepsilon,\zeta\}} = h(a^t_{\bold c\{\varepsilon,\zeta\}})$ 
so (ii) holds]
\sn
\item "{$(iii)$}"  if $\varepsilon > \zeta,\varepsilon \in w$, then 
the equation 
$x^t_\zeta - x^t_\varepsilon \le a^t_{\bold c\{\varepsilon,\zeta\}}$ \nl
[why?  the meaning of the demand is that we should check

$$
h(x^t_\zeta) - h(x^t_\varepsilon) \le h(a^t_{\bold c\{\varepsilon,\zeta\}})
$$
\mn
that is

$$
h(x^t_\zeta) - x^t_\varepsilon \le a^t_{\bold c\{\varepsilon,\zeta\}}
$$
\mn
that is

$$
\xi \in w \cap \zeta \Rightarrow (x^t_\xi - a^t_{\bold c\{\zeta,\xi\}}) -
x^t_\varepsilon \le a^t_{\bold c\{\varepsilon,\zeta\}}
$$
\ermn
for this it suffices to show

$$
\xi \in w \cap \zeta \Rightarrow B_{t,w} \models x^t_\xi - x^t_\varepsilon \le
a^t_{\bold c\{\xi,\zeta\}} \cup a^t_{\bold c\{\zeta,\varepsilon\}}
$$
\mn
but we know

$$
\xi < \zeta \and \xi \in w \and \zeta \in w \Rightarrow
B_{t,w} \models x^t_\xi - x^t_\varepsilon \le a^t_{\bold c\{\xi,\varepsilon\}}
$$
\mn
and $\langle a^t_i:i < \theta \rangle$ is increasing and

$$
\bold c\{\xi,\varepsilon\} \le \text{ Max}\{\bold c\{\xi,\zeta\},
\bold c\{\zeta,\varepsilon\}\}
$$
\mn
so we are done.]  \hfill$\square_{\scite{1.1}}$
\enddemo 
\bigskip

\proclaim{\stag{1.4} Observation}  1) Assume
\mr
\item "{$(*)$}"  $\lambda = \mu$ is weakly compact $> \kappa$.
\ermn
If $A$ is a Boolean Algebra, $A \subseteq B_1,A \subseteq B_2,
B = B_1 \underset A {}\to \oplus \,B_1,|A| \le \kappa$ and: \nl
Depth$^+(B) > \lambda$ or just Length$^+(B) > \lambda$, 
\ub{then} $\dsize \bigvee^2_{t=1} \text{ Depth}^+(B_t) > \mu$. \nl
2) Similarly if $\lambda \rightarrow (\mu)^2_\kappa$,cf$(\lambda) > 2^\kappa$.
\endproclaim
\bigskip

\demo{Proof}  We can find distinct $c_\alpha \in B$ (non-zero) for
$\alpha < \lambda$ such that $\alpha < \beta < \lambda \Rightarrow c_\alpha,
c_\beta$ comparable in $B$. \nl
For each $\alpha$ we can find $b^t_{\alpha,\ell} \in B_t$ ($\ell < n_\alpha,
t=1,2$) such that $c_\alpha = \dbcu^{n_\alpha-1}_{\ell=0}
(b^1_{\alpha,\ell} \cap b^2_{\alpha,\ell})$.  Without loss of generality
$\langle b^2_{\alpha,\ell}:\ell < n_\alpha \rangle$ are pairwise disjoint
(in $B_2$).  Without loss of generality $n_\alpha = n(*)$.
\mn
For each $\alpha < \lambda$ there is an ultrafilter $D_\alpha$ of $B$
such that

$$
c_\alpha \in D_\alpha, \dsize \bigwedge_{\beta < \lambda} [c_\beta <_B
c_\alpha \Rightarrow c_\beta \notin D_\alpha]
$$
\mn
(remember: $\{c_\alpha:\alpha < \lambda\}$ is a chain in $B$).

As $\lambda = \text{ cf}(\lambda) > 2^\kappa \ge 2^{|A|}$, without loss of
generality $\dsize \bigwedge_{\alpha < \lambda} D_\alpha \cap A = D^*$.
Let $I^* = \{a \in A:1_A - a \in D^*\}$, it is a maximal ideal of $A$.  
Let $I_t$ be the ideal which $I^*$ generates in $B_t$ and $I$ is the ideal 
which $I^*$ generates in $B$.
So easily $\langle c_\alpha/I:\alpha < \lambda \rangle$ is a chain with no
repetition.  Now easily $B/I = (B_1/I_1) \oplus (B_2/I_2)$; as in $B/I$ 
there is a chain of cardinality $\lambda$, by Monk McKenzie \cite{MoMc} this 
holds in $B_1/I_1$ or in $B_2/I_2$, so \wilog \, in $B_1/I_1$ say it is 
$\langle b_\alpha:\alpha < \lambda \rangle$.  Now for some 
$a_{\{\alpha,\beta\}} \in I^*$, if $b_\alpha/I_1 < b_\beta/I_1$ then

$$
B_1 \models b_\alpha - a_{\{\alpha,\beta\}} < b_\beta -
a_{\{\alpha,\beta\}}.
$$
\mn
So it is enough to find $X \in [\lambda]^\mu$ and $a \in A$ and truth value
$\bold t$ such that $a_{\{\alpha,\beta\}} = a \and [B \models c_\alpha <
c_\beta \Leftrightarrow \bold t = \text{ truth}]$ for $\alpha < \beta \in X$
which we can.  So $\langle b_\alpha - a:\alpha \in X \rangle$ or
$\langle 1-(b_\alpha -a): \alpha \in X \rangle$ is a strictly increasing
sequence of order type $\lambda$ .  (Of course, for the version with depth
not length, $\bold t$ is redundant.  In more details for each 
$\alpha < \beta$ choose if possible $\ell = \ell_{\alpha,\beta}
< n(*)$ and $a = a_{\alpha,\beta} \in A \backslash \{0\}$ such that
$b^1_{\alpha,\ell} \cap a,b^1_{\beta,\ell} \cap a$ are distinct but
comparable and $\bold t_{\alpha,\beta} \in \{0,1\}$ be such that
$b^1_{\alpha,\ell_\beta} < b^1_{\alpha,\ell_\beta} \equiv
t_{\alpha,\beta} = 1$. \nl
Define a colouring $c:[\lambda]^2 \rightarrow A \times n \times 3$ by: if
$a_{\alpha,\beta},\ell_{\alpha,\beta}$ are well defined,
$c\{\alpha,\beta\} = (a_{\{\alpha,\beta\}},
\ell_{\{\alpha,\beta\}},\bold t_{\alpha,\beta}) \in A \times n \times 2$, 
if not, $c\{\alpha,\beta\} = (0,0,2)$. \nl
It is enough to note
\mr
\item "{$(*)$}"  for no $X \in [\lambda]^{(2^{|A|})^+}$ is $c \restriction
[X]^2$ constantly $(0,0,2)$. \nl
[Why?  Otherwise, \wilog \, $X = (2^{|A|})^+$ and repeat the 
proof above for it.]   \hfill$\square_{\scite{1.4}}$
\endroster
\enddemo
\bigskip

\proclaim{\stag{1.5} Claim}  1) In \scite{1.1} we may replace in clauses
(b), (c) of the assumption the usual order of the ordinals by a linear 
order $<^*$, provided that we weaken clause $(\gamma)$ of the conclusion by
\mr
\item "{$(\gamma)^-$}"  Length$^+(B_1 \underset A {}\to \oplus B_2) = 
\lambda^+$
\ermn
2)  In $\boxtimes_{\lambda,\mu,\theta}$ of \scite{1.1}, we can omit clause
(e) as it follows. \nl
3)  If $\bold c,\lambda,\mu,\theta$ satisfies (a),(b),(c) of
$\boxtimes_{\lambda,\mu,\theta}$ and Rang$(\bold c)$ has no last element,
\ub{then} for some regular $\sigma \le \theta$ we have 
$\boxtimes_{\lambda,\mu,\sigma}$. \nl
4)  If $\boxtimes_{\lambda,\mu,\theta}$ and $\mu \le \mu_1 \le \lambda_1 \le
\lambda$, \ub{then} $\boxtimes_{\lambda_1,\mu_1,\theta}$.
\endproclaim
\bigskip

\demo{Proof}  1) Same proof as in \scite{1.1}. \nl
2) Add to $w_\alpha$ dummy members to increase $i$ (and included in the
next proof). \nl
3) Let $\delta^* = \text{ sup Rang}(\bold c)$, and let $\sigma =
\text{ cf}(\delta^*)$ and $\langle \gamma_\varepsilon:\varepsilon < \sigma
\rangle$ be increasing continuously with limit $\delta^*$, a limit ordinal. 
Define $\bold c':[\lambda]^2 \rightarrow \sigma$ by $\bold c'\{\alpha,\beta\}
= \text{Min}\{\varepsilon < \sigma:\bold c\{\alpha,\beta\} < 
\gamma_\varepsilon\}$ so $\bold c'\{\alpha,\beta\}$ is always 
a successor ordinal.

Let us prove (c) + (e).  So let $\varepsilon(*) < \sigma$ and let
$w_\alpha = \{\zeta_{\alpha,\ell}:\ell < n\},\zeta_{\alpha,0} < \zeta_{\alpha
,1} < \ldots < \zeta_{\alpha,n-1} < \lambda$ for $\alpha < \mu$.  Let
$\varepsilon_\alpha = \text{ max}\{\bold c'\{\zeta_{\alpha,\ell},
\zeta_{\alpha,m}\}:\ell < m < n\}$ and as cf$(\mu) > \sigma$ (so
cf$(\mu) > \theta$ is an overkill) \wilog \, $\varepsilon_\alpha$ is constant
so $\varepsilon^* = \text{ max}\{\varepsilon(*)+1,\varepsilon_\alpha +1:
\alpha < \lambda\} < \sigma$, so as $\delta^* = \text{ sup(Rang } \bold c)$, 
for some $\alpha^* < \beta^* < \lambda$ we have $\bold c\{\alpha^*,\beta^*\} >
\gamma_{\varepsilon^*}$, hence $\bold c'\{\alpha^*,\beta^*\} > \varepsilon^*$.
Without loss of generality for all $\alpha$'s the truth value of
$\zeta_{\alpha,\ell} < \alpha^*,\zeta_{\alpha,\ell} > \alpha^*,
\zeta_{\alpha,\ell} < \beta^*,\zeta_{\alpha,\ell} > \beta^*$ are the same.
Now apply the ``old clause (c)" to $w'_\alpha = w_\alpha \cup
\{\alpha^*,\beta^*\}$ and we can find $\alpha \ne \beta,i,j,a$ there.  Now
$\alpha,\beta,i' = \text{ Min}\{\varepsilon:j < \gamma_\varepsilon\},j' =:
i'+1$ are as required. \nl
4) Trivial.  \hfill$\square_{\scite{1.5}}$ 
\enddemo
\bigskip

\definition{\stag{1.6} Definition}  1)  $Qr_2(\lambda,\mu,\theta)$ means:
\mr
\item "{$(*)$}"  \ub{if} $f:[\lambda]^2 \rightarrow {\Cal P}(\theta)
\backslash \{\emptyset\}$ satisfies \nl
$\alpha < \beta < \gamma \Rightarrow \emptyset \ne
f\{\alpha,\beta\} \cap f\{\beta,\gamma\} \subseteq f\{\alpha,\gamma\}$, and
\nl
$a \in [\theta]^{< \aleph_0}$ \ub{then}
for some $X \in [\lambda]^\mu$ we have \nl
$\dbca_{\alpha \ne \beta \in X} f\{\alpha,\beta\} \nsubseteq a$.
\ermn
2) $NQs_2(\lambda,\mu,A,I)$ where $\lambda \ge \mu$ are (infinite)
cardinals and $A$ is a Boolean algebra and $I$ is an ideal of $A$ means that
there is a function $f:[\lambda]^2 \rightarrow {\Cal J}(I) =: \{J \subseteq
I:J \text{ is non-empty closed upward and closed under intersection but }
0_A \notin I\}$ such that
\mr
\item "{$(a)$}"  if $\alpha < \beta < \gamma$, then
$f\{\alpha,\gamma\} \supseteq f\{\alpha,\beta\} \cap f\{\beta,\gamma\}$
\sn
\item "{$(b)$}"  for some $a^* \in I$ for no $X \in [\lambda]^\mu$ and 
$b \in I$ do we have $a < b$ and:
$$
\alpha \ne \beta \in X \Rightarrow b \in f(\{\alpha,\beta\}).
$$
We say in this case ``$f$ witnesses $NQs_2(\lambda,\mu,A,I)$".
\ermn
3) $NQs^+_2(\lambda,\mu,A,I)$ means that some $f:[\lambda]^2 \rightarrow I$
witnesses it which means that \nl
$f':[\lambda]^2 \rightarrow {\Cal J}(I)$ which is defined by $f'\{\alpha,
\beta\} =: \{b \in I:f\{\alpha,\beta\} \le b\}$ witnesses $NQs_2(\lambda,\mu,
A,I)$. \nl
4) $NQs^*_2(\lambda,\mu,A,I)$ means $NQs_2(\lambda,\mu,A,I)$ is witnessed
by some $f$ which satisfies
\mr
\item "{$(c)$}"  if $n < \omega,\zeta_{\alpha,0} < \ldots < 
\zeta_{\alpha,n-1}$ for $\alpha < \mu$ and $a \in I$, \ub{then} 
for some $\alpha < \beta$ and $b$ we have: $a < b \in I$ we have:
\sn
{\roster
\itemitem{ $(i)$ }  $\zeta_{\alpha,\ell} \le \zeta_{\beta,\ell}$ and
$\zeta_{\alpha,\ell} = \zeta_{\beta,m} \rightarrow \ell = m$
\sn
\itemitem{ $(ii)$ }  if $\ell < m < n$ \ub{then} \nl
$b \in f\{\zeta_{\alpha,\ell},\zeta_{\alpha,m}\} = f\{\zeta_{\beta,\ell},
\zeta_{\beta,m}\}$
\sn
\itemitem{ $(iii)$ }  if $\ell,m < n$ \ub{then} \nl
$b \in f\{\zeta_{\alpha,\ell},\zeta_{\beta,m}\} = f(\zeta_{\beta,\ell},
\zeta_{\beta,m})$ \ub{or} $b$ does not belong to \nl
$f\{\zeta_{\alpha,\ell},\zeta_{\beta,m}\} \cup f\{\zeta_{\beta,\ell},
\zeta_{\alpha,m}\}$.
\endroster}
\ermn
5) In part (1) addition of the letter $N$ (that is 
$NQr_2(\lambda,\mu,\theta)$) means the negation; similarly in parts 
(2), (3) omitting $N$ means the negation.  We can 
replace $\mu$ by $D$, a filter on $\lambda$ meaning replacing ``there is 
$X \in [\lambda]^\mu$ such that ..." by 
$``\{\text{otp} X;X \subseteq \mu \text{ is such that } \ldots\} \in D"$. \nl
6)  If we omit $A$ we mean $I$ is a Boolean ring, $A$ the Boolean Algebra it
generates. \nl
If we omit $A$ and $I$ and write $\theta$ we mean in the $N$-version,
``for some $A,I,|A| \le \theta$" (so without $N$ for every such $A,I$).
\enddefinition
\bn
Among obvious implications are
\proclaim{\stag{1.6A} Claim}  1) $NQr_2(\lambda,\mu,\theta)$ \ub{implies}
$NQs_2(\lambda,\mu,{\Cal P}(\theta),[\theta]^{< \aleph_0})$. \nl
2) $NQs^+_2(\lambda,\mu,A,I)$ \ub{implies} $NQs_2(\lambda,\mu,A,I)$. \nl
3) If $\lambda = \mu$ is weakly compact $> \theta$ \ub{then} $Qr_2(\lambda,
\mu,\theta)$.
\endproclaim
\bigskip

\proclaim{\stag{1.7} Claim}  Assume cf$(\lambda) > 2^\theta$ and 
$Qr_2(\lambda,\mu,\theta)$ or just $Qs_2(\lambda,\mu,\theta)$. \nl
1) If $A,B_1,B_2$ are Boolean Algebras, $A \subseteq B_1,A \subseteq B_2,
|A| \le \theta$ and \nl
Depth$^+(B_1 \underset A {}\to \oplus B_2) > \lambda$ \ub{then}
$\dsize \bigvee^2_{t=1} \text{ Depth}^+(B_t) > \mu$. \nl
2) Similarly for length.
\endproclaim
\bigskip

\demo{Proof}  We start as in the proof of \scite{1.4}(1), getting
$I^*,I,I_1,I_2$ and find $t \in \{1,2\}$ and $b_\alpha \in B_t$ such that:
\mr
\item "{$(*)_1$}"  $\langle b_\alpha/I_t:\alpha < \lambda \rangle$ is strictly
increasing in $B_t/I_t$
\sn
\item "{$(*)_2$}"  $I^*$ is an ideal of $A$ (in fact, a maximal one), and
$I_t$ is the ideal of $B_t$ which $I^*$ generates.
\ermn
Let $f:[\lambda]^2 \rightarrow I^*$ be defined as
follows: for $\alpha < \beta < \lambda$ we let $f\{\alpha,\beta\} =$ \nl
$\{d \in I^*:b_\alpha - d < b_\beta - d\}$.  So $f\{\alpha,\beta\}$ is a
non-empty subset of $I^*$, upward closed and closed under intersection
(remember $\alpha < \beta < \lambda
\Rightarrow b_\alpha/I^* < b_\beta/I^*$).  Now the assumption of $(*)$ of
Definition \scite{1.6} holds as $d \in f\{\alpha,\beta\} \cap f\{\beta,
\gamma\},\alpha < \beta < \gamma$ implies $B_t \models b_\alpha - d <
b_\beta - d < b_\gamma - d$ so $B_t \models b_\alpha - d < b_\gamma - d$ 
hence $d \in f\{\alpha,\gamma\}$, thus proving the $\subseteq$.  As
for the $\ne \emptyset$, if $d_1 \in f\{\alpha,\beta\},d_2 \in d\{\beta,
\gamma\}$ then $d_1 \cup d_2 \in I^*$ belongs to $f\{\alpha,\gamma\}$ (note
then $b_\alpha/I^* < b_\gamma/I^*$ so the strict inequality is not a problem).
\nl
2) Similar only now $\{b_\alpha/I_t:\alpha < \lambda\}$ is a chain with no
repetitions.  \hfill$\square_{\scite{1.7}}$
\enddemo
\bigskip

\proclaim{\stag{1.8} Claim}  If $\lambda \rightarrow [\mu]^2
_{2^\theta,< \aleph_0}$ then $Qs_2(\lambda,\mu,\theta)$.
\endproclaim
\bigskip

\demo{Proof}  Straight.
\enddemo
\bigskip

\remark{Remark}  So we have consistency results by \cite{Sh:276},
\cite{Sh:546} in fact by the proofs $2^\theta$ can be reduced to $\theta$.
\endremark
\bigskip

\proclaim{\stag{1.8A} Claim}  1) Assume $\theta = \text{ cf}(\theta) < \chi
= \chi^{< \chi} < \lambda,\lambda$ measurable, then in $V_1 =
V^{\text{Levy}(\chi,< \lambda)}$ we have $\neg \boxtimes_{\lambda,\lambda,
\theta}$, moreover $Qr_2(\lambda,\lambda,\mu)$. \nl
2) If in $V_1,\lambda = \chi^+,\chi > \theta$ and
\mr
\item "{$(*)$}"  $D$ is a normal filter on $\lambda$ for which in the game
$\Game = \Game(\lambda,D,\theta)$ of length $\theta +1$ between the even and
odd players choosing $A_i \in D^+$ for $i \le \theta$ decreasing (of course,
for $i$ even, the even player chooses $A_i$, for $i$ odd, the odd player
chooses $A_i$), the even player can guarantee that for limit $\delta \le 
\theta,\dbca_{i < \delta} A_i \in D^+$.
\ermn
Then $V_2$ satisfies the conclusion in part (1).
\endproclaim
\bigskip

\demo{Proof}  1) By part (2) by \cite{JMMP} (see on the subject \cite{Sh:b}
or \cite{Sh:f}). \nl
2) Let $\bold c:[\lambda]^2 \rightarrow \theta$ exemplifies
$\boxtimes_{\lambda,\lambda,\theta}$.  Let $St$ be a winning strategy of even.
We chose by induction on $i < \theta,A_i,B_i$ such that
\mr
\item "{$(a)$}"  $\langle A_j:j \le i \rangle$ is a play of $\Game(\lambda,
D,\theta)$ in which even use his winning \nl
strategy $St$
\sn
\item "{$(b)$}"  $\langle B_j:j \le i \rangle$ is a play of $\Game(\lambda,
D,\theta)$ in which even use his winning \nl
strategy $St$
\sn
\item "{$(c)$}"  for $i$ odd for some $\gamma_i < \lambda$ and $j_i \in
(i,\theta)$ we have

$$
\gamma_i < \text{ Min }A_i
$$

$$
\gamma_i < \text{ Min }B_i
$$

$$
(\forall \alpha \in A_i)[\bold c\{\gamma_i,\alpha\} < j_i]
$$

$$
(\forall \beta \in B_i)[\bold c\{\gamma_i,\beta\} \ge j_i]
$$
\ermn
For $i$ even we have no free choice.
\sn
For $i$ odd we ask
\mr
\item "{$(*)$}"  is there $\gamma < \lambda$ such that \nl
$j < \theta \Rightarrow B_{i,j} = \{\beta \in B_{i-1}:\beta > \gamma$ and 
$\bold c\{\gamma,\beta\} \ge j\} \in D^+$?
\ermn
If yes, choose such $\gamma = \gamma_i$, so
$A_{i-1} \backslash (\gamma_i +1) = \dbcu_{j < t}\{\alpha \in A_{i-1}:\bold c
\{\gamma_i,\alpha\} = j\}$.  Hence for some $j$

$$
\{\alpha:\alpha \in A_{i-1},\alpha > \gamma_i,\bold c\{\gamma_i,\alpha\} < j
\} \in D^+.
$$
\mn
Choose this set as $A_i$ and $B_{i,j}$ as $B_i$.

If no, let for $\gamma < \lambda,j(\gamma) < \theta$ be a counterexample.  So
by the normality of the filter 
$B' =: \{ \beta \in B_{i-1}:\text{ for every } \gamma < \beta \text{ we have }
\bold c\{\gamma,\beta\} < j(\gamma)\} = \emptyset$ mod $D$, so
$B_{i-1} \backslash B' \in D^+$ hence for some $j^* < \theta$ we have
$B^* = \{\gamma \in B_{i-1}:\gamma \notin B' \text{ and } j(\gamma) = j^*\}
\in D^+$.

So $B^* \in D^+$ and for $\gamma_1 < \gamma_2$ in $B^*,\bold c\{\gamma_1,
\gamma_2\} < j(\gamma_1) = j^*$.
\mn
This contradicts the choice of $\bold c$. \nl
So we succeed to choose $\langle A_i:i \le \theta \rangle,\langle B_i:
i \le \theta \rangle$, so we can find $\alpha \in A_\theta$ and $\beta \in
B_\theta \backslash (\alpha +1)$ (as $A_\theta,B_\theta \in D^+$). \nl
So for each $i < \theta$
\mr
\item "{$(*)_i$}"  $\gamma_i < \alpha < \beta,\bold c\{\gamma_i,\alpha\}
< j_i,\bold c\{\gamma_i,\beta\} \ge j_i$ so $\bold c\{\gamma_i,\alpha\} <
\bold c\{\gamma_i,\beta\}$.
\ermn
But $\bold c\{\gamma_i,\beta\} \le \text{ Max}\{\bold c\{\gamma_i,\alpha\},
\bold c\{\alpha,\beta\}\}$.  Hence $\bold c\{\alpha,\beta\} \ge
\bold c\{\gamma_i,\beta\}$ but the latter is $\ge j_i$ and $j_i \ge i$, so
$\bold c\{\alpha,\beta\} \ge i$.  As this holds for any $i < \theta$ we have
gotten a contradiction. 
\enddemo
\bigskip

\remark{\stag{1.8B} Remark}  1) We can weaken the demand on the even player
in the game $\Game(\lambda,D,\theta)$ for $\delta = \theta$ to the demand
$\dbca_{i < \delta} A_i \ne \emptyset$ is enough, as we can:
\mr
\item "{$(\alpha)$}"   make $A_1 \cap B_1 = \emptyset$
\sn
\item "{$(\beta)$}"  if $i = 4j+1$ retain demand (c) in the proof but if
$i=4j+3$ uses a similar demand interchanging $A_i$ and $B_i$.
\ermn
2) Moreover, instead ``even has a winning strategy" it is enough, that
``odd has no winning strategy for winning at least one of two plays, 
played simultaneously".
\nl
Now we can deal with other variants.
\endremark
\bn
We may wonder what is required from $A$.
\proclaim{\stag{1.10} Claim}  Assume $\boxtimes_{\lambda,\mu,\theta}$ and $A$
is a Boolean Algebra of cardinality $< \text{ cf}(\mu)$ (for simplicity) and
\mr
\item "{$(*)$}"  there are $a^t_i \in A$ for $i < \theta,t \in \{1,2\}$
such that:
{\roster
\itemitem{ $(i)$ }  $i < j < \theta \Rightarrow A \models a^t_i \le a^t_j$
\sn
\itemitem{ $(ii)$ }  $a^1_i \cap a^2_i = 0$ for $i < \theta$
\sn
\itemitem{ $(iii)$ }   for every $b \in A,t \in \{1,2\}$ there is
$i_t(b) < \theta$ such that: \nl
if $D$ is an ultrafilter on $A$ such that $b \in D,a^{3-t}_i \notin D$ for
$i < \theta$ and $j \in [i_t(b),\theta)$ and $\dsize \bigwedge_{i < j}
a^t_i \notin D$, \ub{then} there is an ultrafilter $D'$ on $A$ such that
$a^t_j \in D',b \in D'$ and $\dsize \bigwedge_{i < j} a^t_i \notin D'$.
\endroster}
\ermn
\ub{Then} $(\lambda,\mu) \in SpDpFP(A)$.
\endproclaim
\bigskip

\demo{Proof}  Similar to the proof of \scite{1.1}, only now fixing 
$\tau_\alpha$ we also fix the
parameters from $A$, say $b_i,\dotsc,b_{k^*}$ and choose $i(*) < \theta$ 
above sup$\{j_t(b):t \in \{1,2\}\}$ and $b \in \langle b_1,\dotsc,b_{k^*}
\rangle_A$.  \hfill$\square_{\scite{1.10}}$
\enddemo
\bigskip

\proclaim{\stag{1.11} Claim}  Assume
\mr
\item "{$(a)$}"  $A$ is a Boolean Algebra and $I_1,I_2$ are ideal of $A$
and $I_1 \cap I_2 = \{0\}$
\sn
\item "{$(b)$}"  for $\ell = 1,2$ we have $NQs^*_2(\lambda,\mu,A,I_\ell)$ and
$|A| < \text{ cf}(\mu)$
\sn
\item "{$(c)$}"  $I_1 \cup I_2$ generates $A$ (or less).
\ermn
\ub{Then} there are Boolean Algebras $B_1,B_2$ extending $A$ such that
Depth$^+(B_1 {\underset A {}\to \oplus} B_2) = \lambda^+$, Depth$^+(B_1) \le
\mu$, Depth$^+(B_2) \le \mu$.
\endproclaim
\bigskip

\remark{Remark}  We can weaken clause (c).  We may wonder on using more
ideals.
\endremark
\bigskip

\demo{Proof}  Like \scite{1.1}.
\enddemo
\bn
Note
\proclaim{\stag{1.12} Claim}  1) Assume
\mr
\item "{$(a)$}"  $\theta = \text{ cf}(\theta) < \mu$, and
\sn
\item "{$(b)$}"  $\mu$ is strong limit singular, cf$(\mu) = \kappa < \theta$
\sn
\item "{$(c)$}"  for any $\bold c:[\theta]^2 \rightarrow \kappa$ satisfying
$\alpha < \beta < \gamma \Rightarrow \bold c\{\alpha,\gamma\} \subseteq
\text{ Max}\{\bold c\{\alpha,\beta\},\bold c\{\beta,\gamma\}\}$, \ub{there}
is $X \in [\theta]^\theta$ such that Rang$(\bold c \restriction [X]^2)$ is a
bounded subset of $\kappa$.
\ermn
\ub{then} $\{\delta < \mu^+:\text{cf}(\delta) = \theta\} \in I[\lambda]$. \nl
2) We can replace $(b)$ by
\mr
\item "{$(b)^-$}"  $\mu > \theta > \kappa = \text{ cf}(\mu),(\forall \sigma
< \mu)[\sigma^{< \theta} < \mu]$.
\endroster
\endproclaim
\bigskip

\demo{Proof}  See \cite{Sh:108}, \cite{Sh:88a}.
\enddemo
\bn
Toward solving problem 11 of \cite{M} we may consider:
\definition{\stag{1.13} Definition}  1) $NQt(\lambda,\mu,A,I)$ means:
\mr
\item "{$(a)$}"  $A$ is a Boolean Algebra
\sn
\item "{$(b)$}"  $I$ is an ideal of $A$
\sn
\item "{$(c)$}"  $\lambda \ge \mu \ge \aleph_0$
\sn
\item "{$(d)$}"  there is a witness $(\bold c_1,\bold c)$ which means
{\roster
\itemitem{ $(\alpha)$ }  $\bold c_\ell:[\lambda]^2 \rightarrow {\Cal J}
[I]$ (see \scite{1.6}(2))
\sn
\itemitem{ $(\beta)$ }  $(\bold c_1\{\alpha,\beta\}) \cap (\bold c_2
\{\alpha,\beta\}) = 0_A$
\sn
\itemitem{ $(\gamma)$ }  if $\alpha < \beta < \gamma < \lambda$ then
$(\bold c_\ell\{\alpha,\beta\}) \cap (\bold c_\ell\{\beta,\gamma\})
\subseteq \bold c_\ell\{\alpha,\gamma\}$
\sn
\itemitem{ $(\delta)$ }  for no $X \in [\lambda]^\mu$ and $\ell \in
\{1,2\}$ and $\alpha \in I \backslash \{0\}$ do we have:

\block
$(*) \qquad$ for $\alpha < \beta$ in $X$ we have $d \in \bold c_\ell
\{\alpha,\beta\}$.
\endblock
\endroster}
\ermn
2) $NQt(\lambda,\mu,\theta)$ means $NQt(\lambda,\mu,A,I)$ for some $A,I$
such that $|A| \le \theta$.
\enddefinition
\bn
\ub{\stag{1.14} Fact}:  1) If $NQt(\lambda,\mu,A,I)$ and $\lambda \ge \mu \ge
\text{ cf}(\mu) > \theta$, \ub{then} the conclusion of \scite{1.1} holds for
the Boolean Algebra $A$. \nl
2) Assume $\lambda \ge \mu,\text{cf}(\lambda) > 2^{|A|}$ and $A \subseteq
B_t$ (for $t = 1,2$), Depth$^+(B_t) \le \mu$ and Depth$^+ (B_1
\underset A {}\to \bigoplus B_2) > \lambda$.  \ub{Then} 
$NQt(\lambda,\mu,A,I)$ for some maximal ideal $I$ of $A$.
\bigskip

\demo{Proof}  Like the earlier ones.
\enddemo
\newpage

\head {\S2 On the family of homomorphic images of a Boolean Algebra} \endhead  \resetall
\bn
Our best result is \scite{2.6}(2) but we first deal with more specific
cases.
\proclaim{\stag{2.1} Lemma}  Assume $\kappa < \mu < \lambda$ and $\lambda$ is
a strong limit singular of cofinality $\kappa$ and cf$(\mu) = \kappa$.

\ub{Then} every Boolean Algebra of cardinality $\lambda$ has a homomorphic
image of cardinality $\in [\mu,2^{< \mu}]$ or $= \mu^{\aleph_0}$.
\endproclaim
\bn
First we prove
\proclaim{\stag{2.2} Claim}  Assume $\kappa = \text{ cf}(\mu) < \mu < 
2^{< \mu}$.  Any Boolean Algebra $B$ of cardinality $\ge \chi =:
\dsize \sum_{\theta < \mu} (2^\theta)^+$ has a homomorphic image of 
cardinality $\in [\mu,\chi)$; note that $(2^{< \mu})^+ \ge \chi$.
\endproclaim
\bigskip

\demo{Proof}  Let $\mu = \dsize \sum_{i < \kappa} \mu_i$ such that
$i < j \Rightarrow \mu_i < \mu_j$ and moreover $i < j \Rightarrow \mu <
2^{\mu_i} \le 2^{\mu_j} \le 2^{< \mu}$ and $\mu_i > \kappa$ and if
$\langle 2^{\mu_i}:i < \kappa \rangle$ is not eventually constant it is
strictly increasing.

If $B$ has an independent subset of cardinality $\mu$, then it has a
homomorphic image of cardinality $\in [\mu,\mu^{\aleph_0}]$ but
$\mu^{\aleph_0} < \chi$ so \wilog \,
there is no independent $X \subseteq B$ of cardinality $\mu$ hence for 
no $\theta < \mu$ does $B$ satisfy the $\theta^+$-c.c., hence $c(B) > \theta$.
So $c(B) \ge \mu$ which is singular hence by a theorem of Erd\"os and Tarski,
$B$ has an antichain $\{a_\alpha:\alpha < \mu\}$.  For each $i < \mu$, let
$B_i$ be the subalgebra of $B$ generated by $\{a_\alpha:\alpha < \mu_i\}$
let $B^c_i$ be the completion of $B_i$.  So id$_{B_i}$ is a homomorphism
from $B_i$ into $B^c_i$ hence it can be extended to a homomorphism $f_i$ from
$B$ into $B^c_i$. \nl
Clearly $B_i \subseteq \text{ Rang}(f_i) \subseteq B^c_i$ so $\mu_i \le
(\text{Rang}(f_i)) \le |B^c_i| \le 2^{\mu_i}$.
\sn

If for some $i,|\text{Rang}(f_i)| \ge \mu$ we are done: $B'_i =: \text{ Rang}
(f_i)$ is a homomorphic image of $B,\mu \le |B'_i| \le |B^c_i| \le 2^{\mu_i} 
< \chi$.  Otherwise, we have $i < \kappa \Rightarrow |B'_i| < \mu$, let
$B^* = \dsize \prod_{i < \kappa} B'_i$, and we define a homomorphism $f$ from
$B$ into $B^*,f(x) = \langle f_i(x):i < \kappa \rangle$.  Clearly $B^*$ is a
Boolean Algebra and $f$ a homomorphism from $B$ into $B^*$.  Now let
$B' = \text{ Rang}(f)$, clearly $B'$ is a homomorphic image of $B$ and \nl
$|B'| = |\text{Rang}(f)| \le |B^*| \le \dsize \prod_{i < \kappa} |B'_i|
\le \dsize \prod_{i < \kappa} \mu = \mu^\kappa \le 2^{\mu_0} < \chi$.
On the other hand $f$ is one to one on each $B_i$ (as $f_i$ is) hence is one
to one on $\dbcu_{i < \kappa} B_i$ which has cardinality $\mu$, so
$|B'| \ge \mu$. \nl
So we are done.  \hfill$\square_{\scite{2.2}}$
\enddemo
\bigskip

\remark{\stag{2.2A} Remark}  If $\mu$ is regular, $|B| \ge \mu$, then $B$ has
a homomorphic image of cardinality $\in [\mu,2^{< \mu}]$, this follows by
Juhasz \footnote{on cardinality and weight of subspace of compact spaces}
\cite{Ju1}.
\endremark
\bigskip

\proclaim{\stag{2.3} Claim}  Assume
\mr
\item "{$(a)$}"  for $\ell = 1,2,3$ we have $\langle B^\ell_i:i \le \delta
\rangle$ is an increasing continuous sequence of Boolean Algebras and for
simplicity $B^\ell_0$ is trivial
\sn
\item "{$(b)$}"  $B^2_i \subseteq B^3_i$ and $B^0_i \subseteq B^1_i$ for 
$i \le \delta$
\sn
\item "{$(c)$}"   for $i < \delta$ non-limit, $B_\delta^1$ is complete
\sn
\item "{$(d)$}"  $h_i$ is a homomorphism from $B^2_i$ into $B^0_i$,
increasing with $i$
\sn
\item "{$(e)$}"  if $x \in B^2_{i+1},y \in B^3_i$ and
$B^3_{i+1} \models ``x \cap z = 0"$, \ub{then} for some $z \in B^2_i$ we
have $B^3_i \models z \cap y = 0,B^2_{i+1} \models x \le z$.
\ermn
\ub{Then} we can find a homomorphism $h$ from $B^2_\delta$ into $B^1_\delta$
extending $h_\delta$.
\endproclaim
\bigskip

\demo{Proof}  We choose by induction on $i$ a homomorphism $f_i$ from $B^3_i$
into $B^1_i$ extending $h_i$.
\mn
\ub{For $i=0$}:  Trivial as the $B^\ell_i$ are trivial.
\bn
\ub{For $i$ limit}:  Let $f_i = \dbcu_{j < i} f_i$.
\bn
\ub{For $i=j+1$}:  Let $H_i = \{f:f \text{ a homomorphism from some 
subalgebra } B_f \text{ of } B^3_{i+1}$ \nl
$\text{into } B^1_{i+1}
\text{ extending } f_j \text{ and } h_i\}$.  As $B^1_{i+1}$ is complete, it
is enough to prove that $H_i \ne \emptyset$.  Let $B'_{i+1}$ be the
subalgebra of $B^3_{i+1}$, generated by $B^2_{i+1} \cup B^3_i$, so it is
enough to prove that $f_i \cup h_{i+1}$, induce a homomorphism from 
$B'_{i+1}$ into $B^1_{i+1}$.  Easily it suffices to prove:
\mr
\item "{$(*)$}"  \ub{if} $x \in B^2_{i+1},y \in B^3_i$ are disjoint (in
$B^3_{i+1}$), \ub{then} $h_{i+1}(x),f_i(y)$ are disjoint (in $B^1_{i+1}$).
\ermn
By the assumption $(e)$ there is $z \in B^2_i$ such that 
$B^2_{i+1} \models x \le z,B^3_i \models x \cap z = 0$ 
hence $B^0_{i+1} \models h_{i+1}(x) \le h_{i+1}(z)$ and 
$B^1_i \models f_i(x) \cap f_i(z) = 0$.  As $f_i(z) = h_i(z)
= h_{i+1}(z)$ we are done.  \hfill$\square_{\scite{2.3}}$
\enddemo
\bigskip

\demo{\ub{Proof of \scite{2.1}}}  So 
by \scite{2.2} \wilog \,$\mu$ is strong limit
(as $\mu = 2^{< \mu} \and \mu$ not strong limit $\Rightarrow \mu =
\text{ cf}(\mu)$), so let $\mu = \dsize \sum_{i < \kappa} \mu_i,
\kappa < \mu_i$ and $\dsize \prod_{i < j} 2^{\mu^+_i} < \mu_j$.  
Similarly, let $\lambda = \dsize \sum_{i < \kappa} \lambda_i,
i < j \Rightarrow \lambda_i < \lambda_j$,
moreover, $\lambda_i = \text{ cf}(\lambda_i) > \dsize \prod_{j < i}
2^{\lambda_j}$.  As in the proof of \scite{2.2}, it suffices to deal with the
following three cases.
\mn
\ub{Case A}:  There is an antichain $\{a_\alpha:\alpha < \lambda\}$ of $B$.
As we can replace $B$ by any homomorphic image of cardinality $\ge \lambda$
\wilog \, $\{a_\alpha:\alpha < \lambda\}$ is a maximal antichain of $B$ and 
each $a_\alpha$ an atom.  So \wilog \, $B$ is a subalgebra of ${\Cal P}
(\lambda)$ and 
$a_\alpha = \{\alpha\}$.  Let $B = \dbcu_{i < \kappa} B_i,B_i$ increasing
continuous, $|B_i| \le \dsize \sum_{j < i} \lambda_j$. \nl
We can find $X_i \subseteq \lambda_i$ of cardinality $\lambda_i$ such that:
\mr
\item "{$(*)_i$}"  for each $a \in B_i$ either $(\forall \alpha \in X_i)
(a_\alpha \le b)$ or $(\forall \alpha \in X_i)(a_\alpha \cap b = 0)$
\sn
\item "{$(*)$}"   $\{a_\alpha:\alpha \in X_i\} \subseteq B_{j_i}$ where
$i < j_i < \kappa$.
\ermn
Choose $Y_i \subseteq X_i$ of cardinality $\mu_i$.

Let $Y = \dbcu_{i < \kappa} Y_i$ and let $f:B \rightarrow {\Cal P}(Y)$ be 
the following homomorphism: \nl
$f(a) = a \cap Y$ and let $B'$ be Rang$(f)$, 
so $B' \subseteq {\Cal P}(Y)$ is a homomorphic image of $B$ and

$$
|B'| \ge |\{a_\alpha:\alpha \in Y\}| = |Y| = \dsize \sum_{i < \kappa} |Y_i| =
\dsize \sum_{i < \kappa} \mu_i = \mu
$$

$$
\align
|B'| &\le |\{a \cap Y:a \in B\}| \le \dsize \sum_{i < \kappa} 
|\{a \cap Y:a \in B_i\}| \\
  &= \dsize \sum_{i < \kappa}|\{a:a \subseteq Y \text{ and } (\forall j)
(i \le j < \kappa) a \cap Y_i \in \{\emptyset,Y_i\}| \\
  &\le \dsize \sum_{i < \kappa} (2^\kappa \times 2^{\mu_i}) =
\dsize \sum_{i < \kappa} 2^{\mu_i} = 2^{< \mu} = \mu.
\endalign
$$
\mn
\ub{Case B}:  $B$ satisfies the $\theta$-c.c., $\theta < \lambda,\kappa >
\aleph_0$.  Then $B$ has an independent subset of cardinality $\chi$ for each
$\chi < \lambda$, in particular $\chi = \mu$, say $\{a_\alpha:\alpha < 
\mu\}$.  Let $B'_0$ be the subalgebra of $B$ which $\{a_\alpha:\alpha < 
\mu\}$ generates, and $B^c_0$ is completion, so id$_{B_0}$ can be extended to
a homomorphism $f$ from $B$ into $B^c_0$, let $B_1 = \text{ Rang}(f)$ so
$\mu = |B_0| \le |B_1| \le |B^c_0| \le \mu^{\aleph_0} = \mu$ so $B_1$ is
as required.
\mn
\ub{Case C}:  $B$ satisfies the $\theta$-c.c., $\theta < \lambda,\kappa =
\aleph_0$. \nl
Let $\lambda = \dsize \sum_{n < \omega} \lambda_n,\theta <
\lambda_n,2^{\lambda_n} < \lambda_{n+1},\lambda_n = \text{ cf}(\lambda_n),
\lambda^\theta_n = \lambda_n$.  Trivially we 
can find pairwise disjoint $\langle b_n:n < \omega \rangle$ 
such that $|B \restriction b_n| \ge \lambda^+_n$
(why?  choose by induction on $n,B_n$ such that $\ell < n \Rightarrow b_\ell
\cap b_n = 0,(B \restriction b_n) \ge \lambda^+_n$ and $|B \restriction
(1_B - \dbcu_{\ell \le n} b_\ell)| = \lambda$); if we are stuck in $n$, then
$I =: \{x \in B:x \cap \dbcu_{f < n} b_\ell = 0$ and $|B \restriction x|
\le \lambda_n\}$ has cardinality and we get a contradiction.

We can find for each $n,\langle a^n_\alpha:\alpha < \lambda^+_n \rangle$
such that $a^n_\alpha \le b_n$ and $\langle a^n_\alpha:\alpha < \lambda^+_n
\rangle$ is independent in $B \restriction b_n$.
So some homomorphic image of $B$ has $\lambda$ atoms and we get Case A.
Alternatively, let $B = \dbcu_{n < \omega} B_n,B_n$ a subalgebra of $B$ 
of cardinality
$\lambda_n$ and $B_n \subseteq B_{n+1}$ and $\{b_n:n < \omega\} \subseteq
B_0,a^n_\alpha \in B_n$ for $\alpha < \lambda^+_n$ pairwise distinct.  
We can find $X_n \subseteq \lambda^+_n$ of cardinality $\lambda^+_n$ such
that $\langle a^n_\alpha:\alpha \in X^+_n \rangle$ is an
indiscernible sequence over $B_n$ (\cite{Sh:92}), even more specifically 
for some $B'_n,B_n \subseteq B'_n \subseteq B_{n+1},|B'_n| = \lambda_n$, and 
some disjoint ideals $I_n,J_n$ of $B'_n$, we have:
\mr
\item "{$(*)_1$}"  $x \in I_n \and y \in J_n \and \alpha \in X_n \Rightarrow 
x \cap a^n_\alpha = x \and y \cap a^n_\alpha = 0$
\sn
\item "{$(*)_2$}"  $x \in B'_n \backslash I_n \backslash J_n \and 
\alpha_1 < \ldots < \alpha_n \in X_n \and \dsize \bigwedge^n_{\ell=1} 
b_i \in \{a^n_{\alpha_\ell},-a^n_{\alpha_\ell}\}
\Rightarrow \dbca^n_{\ell =1} b_\ell \cap x \ne 0 \and (\dbca^m_{\ell=1} b_i) 
- x \ne 0$.
\ermn
Let $A$ be the subalgebra of $B$ generated by $\{b_n:n < \omega\} \cup
\{a^n_\alpha:\alpha \in X_n\}$ and $A^c$, its completion and 
$A^* =: \{c \in A^c:\text{ for every } n \text{ large enough}, c \cap b_n \in
\{b_n,0\}\}$. \nl
Clearly $\{b_n:n < \omega\}$ is a maximal antichain of $A$ hence of $A^c$ and 
of $A^*$ and $B_n$ is a free Boolean Algebra generated by $\lambda^+_n$
elements hence $|A^c \restriction b_n| = (\lambda^+_n)^{\aleph_0}$ hence

$$
|A^*| = 2^{\aleph_0} \times (\dsize \sum_{n < \omega} \lambda^{\aleph_0}_n) =
\dsize \sum_{\chi < \mu} \chi^{\aleph_0}.
$$
\mn
Let $D_n$ be an ultrafilter on $B_n \restriction b_n$ such that 
$D_n \cap (I_n \cup J_n) = \emptyset$.  Let $h_n$ be a homomorphism from
$B_n \restriction b_n$ into $A^c \restriction b_n$ such that
$\dsize \bigwedge_{\alpha \in X_n} h_n(a^n_\alpha) = a^n_\alpha$, \nl
$x \in D_n \Rightarrow h_n(x) = b_n = 1_{A^c \restriction b_n}$ and $x \in 
B_n \backslash D_n \Rightarrow h_n(x) = 0$.  Possible by the choice of $X_n$.
Now we define $h:B \rightarrow A^c$ by
$h(x) = \dsize \sum_n h_n(x \cap b_n)$. \nl
Clearly $x \in B_n \Rightarrow h_n(x \cap b_n) \in \{0,b_n\}$ hence
Rang$(h) \subseteq A^*$.  So $h$ is a homomorphism $f$ from $B$ into $A^*$
such that $f(b_n) = b_n$ and $\alpha \in X_n \Rightarrow f(a^n_\alpha) =
a^n_\alpha$.

Hence $\mu = |A| \le |\text{Rang}(f)| \le |A^*| = \dsize \sum_{\chi < \mu}
\chi^{\aleph_0}$ and we are done. \hfill$\square_{\scite{2.1}}$
\enddemo
\bigskip

\demo{\stag{2.5} Observation}  If $\mu$ is strong limit of cofinality
$\aleph_0$ and $\mu \le \lambda < 2^\mu,B$ a Boolean Algebra of cardinality
$\lambda$, \ub{then} ult$(B) > \lambda \Rightarrow$ ult$(B) \ge 2^\mu$.
\enddemo
\bigskip

\demo{Proof}  Straight, or by 
\cite{Sh:454a} (for even a more general setting: a topology).
\enddemo
\bn
We can add
\proclaim{\stag{2.6} Claim}  1) If $\lambda$ is a strong limit singular and
$B$ a Boolean Algebra of cardinality $\lambda$, \ub{then} for some homomorphic
image $B'$ of $B$ we have $|B'| = |\text{Ult}(B')| = \lambda$. \nl
2) If $\lambda_i \ge \mu > \text{ cf}(\mu) = \text{ cf}(\lambda)$,
\ub{then} for some homomorphic image $B'$ of $B$ we have:

$$
\mu \le |B'| \le 2^{< \mu},\mu \le |\text{Ult}(B')| \le
\dsize \sum_{\theta < \mu} 2^{2^\theta}
$$
\mn
(remember: the free Boolean Algebra generated by $\{x_\alpha:\alpha <
\lambda\},B_\lambda$ (or its completion) has a homomorphic image any Boolean
Algebra (any $\sigma$-complete Boolean Algebra) of cardinality $\le \lambda$.)
\endproclaim
\bigskip

\demo{Proof}  1) Without loss of generality $B$ is a Boolean Algebra of
subsets of $\lambda$ and let $\lambda = \dsize \sum_{i < \theta}
\mu_\zeta,\theta=\text{ cf}(\lambda),\mu_\zeta > \theta + 
\dsize \prod_{\xi < \zeta} 2^{\mu_\xi}$ and let $B = \dbcu_{i < \theta}
B_i,B_i$ increasing continuously $|B_\zeta| \le \dsize \sum_{\xi < \zeta}
\mu^+_\xi$.  We know that there are $(a^\zeta_i,\alpha^\zeta_i)$ for
$i < \mu^+_\zeta$ such that $a^\zeta_i \in B,\alpha^\zeta_i < \lambda,
\alpha^\zeta_j \in a^\zeta_i \Leftrightarrow j=i$ (e.g. by the well known
$s(B) \ge \ell g_2(B)$.  So \wilog \, 
$\alpha^\zeta_i < \mu^+_\zeta,a^\zeta_i \in B_{\zeta +1}$.  As we can replace
$B$ by a homomorphic image \wilog \, $\alpha^\zeta_i = \mu_\zeta +i$.  We
can find $X_\zeta \in [\{\alpha:\mu_\zeta \le \alpha < \mu^+_\zeta\}]
^{\mu^+_\zeta}$ for $\zeta < \theta$ such that \nl
$(\forall a \in B_\zeta)((\forall \alpha \in X_i)
(\alpha \in a) \vee (\forall \alpha \in X_i)(\alpha \notin a))$ (like case A 
of the proof of \scite{2.1}, so actually there it suffices. 
Let $X = \dbcu_{\zeta < \theta} X_\zeta$ and let $Y_\zeta = 
\dbcu_{\varepsilon < \zeta} X_\varepsilon$. \nl
Now $X = \dbcu_{i < \theta} X_i,h:B \rightarrow {\Cal P}(X)$ is $h(a) =
a \cap X$.  Easily $B' =: h(B)$ is as required. \nl
Let $B'' = \{a \in B':\text{ for every } i < \theta,X_i \subseteq a
\text{ or } X_i \cap a = \emptyset\}$. \nl
Let $B^+$ be the Boolean Algebra of subsets of $X$ generated by
$B'' \cup B'$.  Now clearly

$$
\align
|\text{Ult}(B')| &\le |\text{Ult}(B^+)| \\
  &\le |\{p \in \text{ Ult}(B''):\text{for every } i < \theta, 
\text{ we have } Y_i \notin p\}| + 
\dsize \sum_{i < \theta} |\text{Ult}(B'' \restriction Y_i)| \\
  &\le 2^{2^\theta} + \dsize \sum_{i < \theta} 2^{2^{|Y_i|}} \le \lambda.
\endalign
$$
\enddemo
\mn
2) Easy (by parts (1), \scite{2.7} below and earlier proofs).
\hfill$\square_{\scite{2.6}}$
\bn
\ub{\stag{2.7} Observation}:  1) Let $B_\sigma$ be the Boolean Algebra 
generated freely by $\{x_\alpha:\alpha < \sigma\}$ and $B^c_\sigma$ its 
completion.  \ub{Then} there is a homomorphic image 
$B'_\sigma$ of $B^c_\sigma$ say by $f,[\sigma]^{< \aleph_0} \subseteq
f(B_\nu) \subseteq B'_\sigma \subseteq B^*_\sigma = 
\{X \subseteq \sigma:|X| \le \aleph_0$ or
$|\sigma \backslash X| \le \aleph_0\}$ so $B'_\sigma$ has 
$\le \sigma^{\aleph_0}$ elements and $\le \sigma^{+ \beth_2}$ ultrafilters 
and $\ge \sigma$ elements. \nl
2) We can make above use $B^*_{\sigma^{\aleph_0}}$ and get $|B'_\sigma| =
\sigma^{\aleph_0}$.
\bigskip

\demo{Proof}  1) Let the homomorphism $f_0:B \rightarrow B^*_\sigma$ 
be induced by $f_0(x_\alpha) = \{\alpha\}$. \nl
Let for $\alpha < \sigma,B^c_{\sigma,\alpha}$ be the complete subalgebra of
$B^c_\sigma$ which $\{x_\beta:\beta < \lambda$ and $\beta \ne \alpha\}$
generated, so $(\forall b)(b \in B^c_{\sigma,\alpha} \and 0 < b < 1
\rightarrow c \cap a_\alpha > 0 \and c - a_\alpha > 0)$.  Let $D$ be an
ultrafilter on $B^c_\sigma$ such that $\{-x_\alpha:\alpha < \sigma\} \subseteq
D$ and for $\alpha < \sigma$ let $D_\alpha$ be an ultrafilter on $B^c_\sigma$
such that: $D_\alpha \cap B^c_{\sigma,\alpha} = D \cap B^c_{\sigma,\alpha}$
and $a_\alpha \in D_\alpha$.  We define a homomorphism $f:B^c_\sigma
\rightarrow {\Cal P}(\sigma)$ and $f(b) = \{\alpha < \sigma:b \in D_\alpha\}$.
Clearly $f$ is a homomorphism from $B^c_\sigma$ into $B'_\sigma$ and if
$c \in B^c_\delta$ for some $\alpha_n < \sigma$ (for $n < \omega$) and
(infinite) Boolean term $\tau,b = \tau(\ldots,x_{\alpha_n},\dotsc)
_{n < \omega}$, so $\alpha \in \sigma \backslash \{\alpha_n:n < \omega\}
\Rightarrow [\alpha \in f(b) \leftrightarrow b \in D_\alpha \leftrightarrow
b \in D]$, so $f(b)$ contains $\sigma \backslash \{\alpha_n:n < \omega\}$ or
is disjoint to it hence $f(b) \in B'_\sigma$, so we are done. \nl
2) Let $\{Y_\gamma:\gamma < \sigma^{\aleph_0}\} \subseteq 
[\sigma]^{\aleph_0}$ be such that 
$\gamma < \beta \Rightarrow Y_\gamma \cap Y_\beta$
finite (clearly possible) let $Y_\gamma = \{\alpha_{\gamma,\ell}:\ell < 
\omega\}$ with no repetitions; moreover, \wilog \nl
$\alpha_{\gamma_1,\ell_2} =
\alpha_{\gamma_2,\ell_2} \Rightarrow \ell_1 = \ell_2 \and \dsize
\bigwedge_{\ell < \ell_1} \alpha_{\gamma_1,\ell} = \alpha_{\gamma_2,\ell}$.
\nl
Let $y_{\gamma,m} = \dbcu_{k < m} \bigl( x_{\alpha_{\gamma,3k}} -
\dbcu_{\ell < 3k} x_{\alpha_{\gamma,\ell}} \bigr) \in B_\sigma,
y_\gamma = \dbcu_m y_{\gamma,m} \in B^c_\sigma$.  
Let $A$ be the subalgebra of $B^c_\sigma$ which $\{y_\gamma:\gamma <
\sigma^{\aleph_0}\}$ generates and $f_0:A \rightarrow 
B^*_{(\sigma^{\aleph_0})}$ be the homormophism induced by 
$f_0(y_\gamma) = \{\gamma\}$, and continue as above.
\hfill$\square_{\scite{2.7}}$
\enddemo
\newpage

\head {\S3 If $d(B)$ is small, then depth or ind are not tiny} \endhead  \resetall
\bigskip

\definition{\stag{3.1} Definition}  1) We say $\langle a_\beta:\beta <
\beta^* \rangle$ is semi-independent if: it is a sequence of distinct 
elements in a Boolean Algebra $B$ and 
for some ideal $I$ on $B$ \nl
for any $\alpha < \gamma < \beta^*$ and $b \in \langle a_\beta:\beta <
\alpha \rangle_B$ we have
\mr
\item "{$(*)_1$}"  $b \in I \Rightarrow b \cap a_\gamma = b \cap a_\alpha$
\sn
\item "{$(*)_2$}"  $b \notin I \Rightarrow \{b \cap a_\alpha,
b \cap a_\gamma\}$ is an independent set in $B \restriction b$, (so e.g. 
$b \cap a_\alpha > 0$)
\sn
\item "{$(*)_3$}"  $b \notin I,b \cap a_\alpha \in I \Rightarrow 
b \cap a_\gamma \in I$
\sn
\item "{$(*)_4$}"  $b \notin I,
b - a_\alpha \in I \Rightarrow b-a_\gamma \in I$.
\ermn
2) $si^+(B) = \text{ Min}\{\lambda:\text{there is no } \langle a_\beta:
\beta < \lambda \rangle \text{ in } B \text{ which is semi-independent}\}$ and
we say $\bar a = \langle a_\beta:\beta < \lambda \rangle$ and $I$
witness $\lambda < si^+(B)$.
Let $si(B) = \sup\{\lambda:$ there is a semi-independent sequence
$\langle a_\beta:\beta < \lambda \rangle$ in $B\}$. \nl
3)  $si^{1+}(B)$ is defined similarly to $si^+(B)$ with 
$\langle a_\beta:\beta < \lambda + 1 \rangle$.  We say 
$I$, \nl
$\langle a_\beta:\beta \le \lambda \rangle$ witness $\lambda < si^{1+}(B)$.
\enddefinition
\bigskip

\demo{\stag{3.2} Fact}  0) $si(B),si^{1+}(B) \le si^+(B) \le (si(B))^+$. \nl
1) ind$^+(B) \le si^+(B)$. \nl
2) $si^+(B) \ge t^+(B)$. \nl
3) $si^+(B) \le \text{ ind}^+(B) + \text{ Depth}^+(B)$. 
\enddemo
\bigskip

\demo{Proof}  0) Read the definition. \nl
1) ind$^+(B) \le si^+(B)$ holds as independent implies semi-independent
for the ideal $\{O_B\}$. \nl
2) Let $\lambda < si^+(B)$ and $\langle a_\beta:\beta < \lambda \rangle,I$ 
witness it. 

Let $D$ be an ultrafilter on $B$ disjoint to $I$.  As we can replace
$\langle a_\beta:\beta < \lambda \rangle$ by $\langle a'_\beta:\beta < \lambda
\rangle$ when $a'_\beta \in \{a_\beta,-a_\beta\}$, \wilog \,
$\dsize \bigwedge_{\beta < \lambda} a_\beta \in D$.  Let $\beta_0 < \ldots
< \beta_{m-1} < \beta_m < \ldots < \beta_{n-1}$.  Let for $k \in [m,n],
b_k = \dbca_{\ell < m} a_{\beta_\ell} \cap \dbca_{\ell \in [m,k)}
(-a_{\beta_\ell})$.  So $b_m \in D$ as $b_{\beta_0},\dotsc,b_{\beta_{m-1}}
\in D$.  So $b_m \notin I$.
\sn
We should prove $b_n > 0$. \nl
Let $k \in [m,n]$ be maximal such that $b_k \notin I$.

So $k$ is well defined, if $k=n$ we are done as then $b_n \notin I \Rightarrow
b_n > 0$, so we can assume $k < n$.  
So by the Definition \scite{3.1}(1)$(*)_2$ with $(b_k,\beta_k,\beta_{k+1})$
here standing for $(b,\alpha,\gamma)$ there we have
$b_{k+1} \cap (-a_{\beta_{k+1}}) = b_k \cap (-a_{\beta_k}) \cap
(-a_{\beta_{k+1}}) = (b_k \cap (-a_{\beta_k})) \cap
(b_k \cap (-a_{\beta_{k+1}}))$ is $> 0$ and by the maximality of $k,
b_{k+1} \in I$, so by \scite{3.1}(1)$(*)_1$ we have 
$b_{k+1} \cap a_{\beta_{k+1}} = 
b_{k+1} \cap a_{\beta_{k+2}} = \ldots$ hence $b_{k+1} \cap (-a_{\beta_k}) =
b_{k+2} = \ldots = b_n$ so $b_n \ne 0$. \nl
3) Let $\langle a_\beta:\beta < \lambda \rangle,I$ witness $\lambda < si^+
(B)$.

If $\langle a_\beta:\beta < \lambda \rangle$ is independent we are done, so
assume not, so let $\beta^* < \lambda$ be minimal such that $\langle
a_\beta:\beta \le \beta^* \rangle$ is not independent modulo $I$; so $\langle
a_\beta:\beta < \beta^* \rangle$ is independent modulo $I$ and for some
$b \in \langle a_\beta:\beta < \beta^* \rangle_B$ satisfying $b > 0$, 
(so $b \notin I$ by the assumption on $\beta^*$), we have 
$b \cap a_{\beta^*} \in I$ or $b - a_{\beta^*} \in I$.  Now by symmetry
\wilog \, the former holds. \nl

So assume $\beta^* \le \gamma_1 < \gamma_2 < \gamma_3 < \lambda$, by
$(*)_3$ of \scite{3.1} we have
$b \cap a_{\gamma_\ell} \in I$ (for $\ell = 1,2,3$) and so by $(*)_1$ of
\scite{3.1} we have 
$b \cap a_{\gamma_1} \cap a_{\gamma_2} = b \cap a_{\gamma_1} 
\cap a_{\gamma_3}$ but $b \notin I$ hence by $(*)_3$ of \scite{3.1} we have
$b-a_{\gamma_3} \notin I$ hence $(b - a_{\gamma_1}) \cap a_{\gamma_2} \ne 0$.

Now $b \cap a_{\gamma_1} \cap a_{\gamma_2} = b \cap a_{\gamma_1} \cap
a_{\gamma_2} \cap a_{\gamma_3} \ge b \cap a_{\gamma_2} \cap a_{\gamma_3}$.
So $\gamma_1 < \gamma_2 \le \gamma_3 \le \gamma_4$ implies $b \cap
a_{\gamma_1} \cap a_{\gamma_2} \ge b \cap a_{\gamma_3} \cap a_{\gamma_4}$.
Let $c_\gamma = b \cap a_{\beta^* + 2\gamma} \cap a_{\beta^* +2 \gamma+1}$,
so $\gamma_1 < \gamma_2 < \lambda \Rightarrow c_{\gamma_1} \ge
c_{\gamma_2}$.  Can equality hold?  Clearly $c_{\gamma_1} \in I$ so 
$b-c_{\gamma_1} \notin I$ hence by $(*)_2$ of \scite{3.1}(1) we have 
$(b-c_{\gamma_1}) \cap a_{\beta^* + 2 \gamma_2} \cap 
a_{\beta^* + 2 \gamma_2+1} > 0$ so necessarily $(b - c_{\gamma_1}) \cap
c_{\gamma_2} > 0$ hence $c_{\gamma_1} \ne c_{\gamma_2}$.  Together
$c_{\gamma_1} > c_{\gamma_2}$.  So $\langle c_\gamma:\gamma < \lambda \rangle$
exemplify $\kappa < \text{ Depth}^+(B)$. \nl
${{}}$  \hfill$\square_{\scite{3.2}}$
\enddemo
\bigskip

\proclaim{\stag{3.3} Claim}  Assume 
$B$ is a Boolean Algebra $\kappa \ge \text{ ind}^+(B),
d(B) \le \lambda =$ \nl
$\lambda^{< \kappa} < |B|$. \nl
1) Depth$^+(B) > \kappa$. \nl
2) If in addition $\lambda^+ \rightarrow (\mu +1)^3_\sigma$ for 
$\sigma < \kappa$, \ub{then} Depth$^+(B) > \mu,si^+(B) > \mu$ [Saharon?]
\endproclaim
\bigskip

\demo{Proof}  1) Let 
$\langle a_\alpha:\alpha < \lambda^+ \rangle$ be a list of 
pairwise distinct elements of $B$.  As $\lambda \ge d(B)$ \wilog \, $B$ is 
a subalgebra of ${\Cal P}(\lambda)$.  Let $B_\alpha = \langle \{a_\beta:
\beta < \alpha\} \rangle_B$ and

$$
E =: \bigl\{\delta < \lambda^+:B_\delta \cap \{a_\alpha:\alpha < \lambda^+\} =
\{a_\alpha:\alpha < \delta\} \bigr\}
$$
\mn
clearly it is a club of $\lambda^+$.

For every $\delta \in S_0 =: \{\delta \in E:\text{cf}(\delta) \ge \kappa\}$
we let $I_\delta =: \{b \in B_\delta:a_\delta \cap b \in B_\delta\}$, so
$I_\delta$ is an ideal of the subalgebra $B_\delta$ of $B$.  Let $\delta \in
S_0$, and $J$ be an ideal on $B_\delta$ and now we try to choose by induction 
on $i < \kappa$, an ordinal $\alpha_{\delta,J,i}$ such that:
\mr
\item "{$(a)$}"  $\alpha_{\delta,J,j} < \alpha_{\delta,J,i} < \delta$ 
for $j<i$
\sn
\item "{$(b)$}"  $\langle a_{\alpha_{\delta,J,j}}/J:j < i \rangle$ is 
independent in the Boolean Algebra $B_\delta/J$.
\ermn
If we succeed, then $\langle a_{\alpha_{\delta,J,i}}:i < \kappa \rangle$
contradict the assumption $\kappa \ge \text{ ind}^+(B)$, so for some
$i(\delta,J) < \kappa$ we have: $\alpha_{\delta,J,i}$ is defined iff 
$i < i(\delta,J)$.
So for some stationary $S_1 \subseteq S_0$ and $i(*) < \delta$ and $\langle
\alpha_i:i < i(*) \rangle$, an increasing sequence of ordinals $< \lambda^+$,
we have $S_1 = \{\delta \in S:i(\delta) = i(*)$ and $i < i(*) \Rightarrow
\alpha_{\delta,I_\delta,i} = \alpha_i\}$. \nl
Let $\langle b_\gamma:\gamma < \gamma(*) \rangle$ list the non-zero Boolean
combination of $\{a_{\alpha_i}:i < i(*)\}$ so $\gamma(*) < \kappa$.  
As $B$ is a subalgebra of ${\Cal P}(\lambda)$ we can choose a function 
$H$ such that Dom$(H) = B \backslash \{\emptyset\},H(c) \in c$.  
Choose a function
$F_\delta$, Dom$(F_\delta) = I_\delta$ and $c \in I_\delta \Rightarrow
F_\delta(c) = c \cap a_\delta \in B_\delta$.
Again for some $Y \subseteq \gamma(*)$ and $\langle x_\gamma,y_\gamma:
\gamma < \gamma(*) \rangle$ we have

$$
S_2 = \{\delta \in S_1:b_\gamma \in I_\delta \Leftrightarrow \gamma \in Y
\text{ and for } \gamma \in \gamma(*) \backslash Y,H(b_\gamma \cap
a_\delta) = x_\gamma,H(b_\gamma - a_\delta) = y_\gamma\}
$$
\mn
is a stationary subset of $\lambda^+$.

For each  $\delta \in S_2$ and $\bold t \in \{0,1\}$ and $\gamma < \gamma(*)$
we try to choose by induction on $i < \kappa$, an ordinal
$\beta_{\delta,\gamma,i}$ such that:
\mr
\item "{$(a)'$}"  $\delta > \beta_{\delta,\gamma,i} > \beta_{\delta,\gamma,j}$
for $j<i$
\sn
\item "{$(b)'$}"  $\beta_{\delta,\gamma,i} > \alpha_j$ for $j < i(*)$
\sn
\item "{$(c)'$}"  $a^{\bold t}_{\beta_{\delta,\gamma,i}} \cap b_\gamma \in
I_\delta$ \nl
(remember $a^{\bold t}_{\beta_{\delta,\gamma,i}}$ is
$a_{\beta_{\delta,\gamma,i}}$ if $\bold t = 1$ and is
$-a_{\beta_{\delta,\gamma,i}}$ if $t=0$)
\sn
\item "{$(d)'$}"  $B_{\delta,\gamma,i}$ is the smallest subalgebra of
$B_\delta$ containing $\{b_\gamma\} \cup \{a_{\alpha_j}:j < i(*)\} \cup
\{a_{\beta_{\delta,\gamma,j}}:j < i\}$
\sn
\item "{$(e)'$}"  if $c \in I_\delta \cap B_{\delta,\gamma,i}$ then
$a_{\beta_{\delta,\gamma,i}} \cap c = F_\delta(c) = a_\delta \cap c$ (in
fact just $c \in \{b_\gamma \cap a^{\bold t}_{\beta_{\delta,\gamma,j}}:
j < i\}$ suffice)
\sn
\item "{$(f)'$}"  if $c \in B_{\delta,\gamma,i}$ and $\bold s \in \{0,1\}$
then $c \cap a^{\bold s}_\delta \ne 0 \Rightarrow H(c \cap a^{\bold t}
_\delta) \in a^{\bold s}_{\beta_{\delta,\gamma,i}}$ (in fact just
$c \in \{b_\gamma \cap a^{\bold t}_{\beta_{\delta,\gamma,j}}:j < i\}$
suffice).
\ermn
Now
\mr
\item "{$(*)$}"  If for some $\delta \in S_2$ and 
$\gamma \in \gamma(*) \backslash Y$ and $\bold t \in
\{0,1\}$ we succeed, we can prove Depth$^+(B) > \kappa$.
\ermn
[Why?  We just prove that $\langle a^{\bold t}_\delta \cap b_\gamma \cap
a^{\bold t}_{\beta_{\delta,\gamma,i}}:i < \kappa \rangle$ is strictly
increasing.  Let
\sn
$j < i < \kappa$, so by clause $(c)'$ we know that
$b_\gamma \cap a^{\bold t}_{\beta_{\delta,\gamma,j}} \in I_\delta$ but
$b_\gamma \cap a^{\bold t}_{\beta_{\delta,\gamma,j}} \in B_{\delta,\gamma,i}$
by clause $(d)'$ hence $a_{\beta_{\delta,\gamma,i}} \cap b_\gamma \cap
a^{\bold t}_{\beta_{\delta,\gamma,j}} = F_\delta(b_\gamma \cap
a^{\bold t}_{\beta_{\delta,\gamma,j}}) = a_\delta \cap b_\gamma \cap
a^{\bold t}_{\beta_{\delta,\gamma,j}}$ by clause $(e)'$ and similarly
$(-a_{\beta_{\delta,\gamma,i}}) \cap b_\gamma \cap a^{\bold t}
_{\beta_{\delta,\gamma,i}} = (-a_\delta) \cap b_\gamma \cap
a^{\bold t}_{\beta_{\delta,\gamma,j}}$.
\sn
So in any case $a^{\bold t}_{\beta_{\delta,\gamma,i}} \cap b_\gamma \cap 
a^{\bold t}_{\beta_{\delta,\gamma,j}} = a^{\bold t}_\delta \cap b_\gamma
\cap a^{\bold t}_{\beta_{\delta,\gamma,j}}$. \nl
So

$$
x \in a^{\bold t}_{\beta_{\delta,\gamma,i}} \cap b_\gamma \cap 
a^{\bold t}_{\beta_{\delta,\gamma,j}} \Leftrightarrow x \in
a^{\bold t}_\delta \cap b_\gamma \cap
a^{\bold t}_{\beta_{\delta,\gamma,j}}.
$$
\mn
So if $x \in a^{\bold t}_\delta \cap b_\gamma$ then

$$
x \in a^{\bold t}_{\beta_{\delta,\gamma,i}} \cap  
a^{\bold t}_{\beta_{\delta,\gamma,j}} \Leftrightarrow x \in
a^{\bold t}_{\beta_{\delta,\gamma,j}}.
$$
\mn
so

$$
x \in a^{\bold t}_{\beta_{\delta,\gamma,j}} \Rightarrow x \in
a^{\bold t}_{\beta_{\delta,\gamma,i}}.
$$
\mn
The above statement means 

$$
x \in a^{\bold t}_\delta \cap b_\gamma \Rightarrow [x \in a^{\bold t}
_{\beta_{\delta,\gamma,j}} \Rightarrow x \in a^{\bold t}_{\beta,\delta,i}]
$$
\mn
hence $\langle a^{\bold t}_\delta \cap b_\gamma \cap a^{\bold t}
_{\beta_{\delta,\gamma,i}}:i < \kappa \rangle$ is $\le$-increasing.  But
$b_\gamma \notin I_\delta$ hence $a^{\bold t}_\delta \cap b_\gamma \notin
B_\delta$, hence for $j < i,a^{\bold t}_\delta \cap b_\gamma \ne 
a^{\bold t}_\delta \cap b_\gamma \cap
a^{\bold t}_{\beta_{\delta,\gamma,j}}$ (as by clause $(c)'$ we know that
$b_\gamma \cap a^{\bold t}_{\beta_{\delta,\gamma,j}} \in I_\delta$ so 
$a^{\bold t}_\delta \cap b_\gamma \cap a^{\bold t}_{\beta_{\delta,\gamma,j}} 
\in B_\delta$).
\sn
So $0 < a^{\bold t}_\delta \cap b_\gamma - a^{\bold t}_\delta \cap b_\gamma
\cap a^{\bold t}_{\beta_{\delta,\gamma,j}} = (b_\gamma - a^{\bold t}
_{\beta_{\delta,\gamma,j}}) \cap a^{\bold t}_\delta$ hence
$x = H((b_\gamma - a^{\bold t}_{\beta_{\delta,\gamma,j}}) \cap
a^{\bold t}_\delta) \in a^{\bold t}_\delta$ is well defined and belongs to
$a_{\beta_{\delta,\gamma,i}}$ by clause $(f)'$.  
So $x$ belongs to $a^{\bold t}_\delta,
b_\gamma - a^{\bold t}_{\beta_{\delta,\gamma,j}},
a^{\bold t}_{\beta_{\delta,\gamma,i}}$ so it exemplifies
$a^{\bold t}_\delta \cap b_\gamma \cap a^{\bold t}_{\beta_{\delta,\gamma,j}}
\ne a^{\bold t}_\delta \cap b_\gamma \cap a^{\bold t}
_{\beta_{\delta,\gamma,i}}$.  So $\langle a^{\bold t}_\delta \cap b_\gamma 
\cap a^{\bold t}_{\beta_{\delta,\gamma,i}}:i < \kappa \rangle$ is strictly
increasing.] \nl
\sn
We have proved $(*)$, so assume toward contradiction
that for $\delta \in S_2,\gamma \in \gamma(*) \backslash Y,
a_{\beta_{\delta,\gamma,i}}$ is well defined iff $i < j(\delta,\gamma)$ 
where $j(\delta,\gamma) < \kappa$.  
So again for some
$\beta_{\gamma,i}(\gamma \in \gamma(*) \backslash Y,i < j(\gamma))$ we have

$$
\align
S_3 = \biggl\{ \delta \in S_2:&\text{for every } \gamma \in \gamma^*
\backslash Y \text{ we have} \\
  &j(\delta,\gamma) = j(\gamma),i < j(\gamma) \Rightarrow 
\beta_{\delta,\gamma,i} = \beta_{\gamma,i} \\
  &\text{and } i < j(\gamma) \Rightarrow B_{\delta,\gamma,i} =
B_{\gamma,i} \biggr\}
\endalign 
$$
\mn
and is stationary.
\sn
Let $\langle b_\gamma:\gamma < \gamma(**) \rangle$ list
$\langle \{a_{\alpha_i}:i < i(*)\} \cup \{a_{\beta_{\gamma,i}}:\gamma \in
\gamma(*) \backslash Y,i < j(\gamma)\} \rangle$, so for some stationary
$S_4 \subseteq S_3$ we have:
$\delta_1,\delta_2 \in S_4 \Rightarrow H_{\delta_1} \restriction
\dbcu_{\gamma,i} B_{\gamma,i} = H_{\delta_2} \restriction \dbcu_{\gamma,i}
B_{\gamma,i}$ and $F_{\delta_1}(y \cap a^{\bold t}_{\delta_1}) =
F_{\delta_2}(y \cap a^{\bold t}_{\delta_1})$, for $\bold t = 0,1,y \in
\dbcu_{\gamma,i} B_{\gamma,i}$ and $I_{\delta_1} \cap \dbcu_{\gamma,i}
B_{\gamma,i} = I_{\delta_2} \cap \dbcu_{\gamma,i} B_{\gamma,i}$.  
Let $\delta_1 < \delta_2$ be in $S_4$ and we get a contradiction. \nl
2) In the proof of part (1) define, for $\gamma < \gamma(*)$ the colouring
$\bold c_\gamma:[S_4]^3 \rightarrow \{0,1,2\}$ by: \nl
for $\delta_0 < \delta_1 < \delta_2$ from $S_1$ we have

$$
\align
\bold c_\gamma\{\delta_0,\delta_1,\delta_2\} \text{ is}: &0 \text{ if }
a^{\bold t}_{\delta_2} \cap a^{\bold t}_{\delta_0} < a^{\bold t}_{\delta_2}
\cap a^{\bold t}_{\delta_1}, \text{ for } \bold t =0 \\
  &1 \text{ if } a^{\bold t}_{\delta_2} \cap a^{\bold t}_{\delta_0} <
a^{\bold t}_{\delta_2} \cap a^{\bold t}_{\delta_1} \text{ for } \bold t=1 \\
  &\qquad \qquad \text{ and it is not } 0 \\
  &2 \text{ if otherwise}.
\endalign
$$
\mn
Lastly, let $\bold c:[S_1]^3 \rightarrow \gamma(*) +1$ be: \nl
for $\delta_0 < \delta_1 < \delta_2$ from $S_1$

$$
\align
\bold c \{\delta_0,\delta_1,\delta_2\} \text{ is } (\gamma,\bold t) 
\text{ if } &\bold c_\gamma \{\delta_0,\delta_1,\delta_2\} = \bold t \in
\{0,1\} \and \\
  &(\forall \gamma' < \gamma)[\bold c_\gamma \{\delta_0,\delta_1,\delta_2\}
=2] \\
  &\text{and is } (\gamma(*),2) \text{ otherwise}.
\endalign
$$
\mn
Now apply the partition property on $S_4$ (noting the ``just" clauses in
$(e)', (f)'$ in proof of part (1)).  It is simpler to apply $\lambda^+
\rightarrow (\lambda^+,(\mu +1)_\sigma)^3$.  \hfill$\square_{\scite{3.3}}$
\enddemo
\bigskip

\demo{\stag{3.4} Conclusion}  1) If $\kappa = \text{ Depth}^+(B) + \text{ ind}
^+(B)$, \ub{then} $|B| \le d(B)^{< \kappa}$. \nl
2)?  If $|B| > |d(B)|^{< \text{ind}^+(B)}$, then for some $\sigma <
\text{ ind}^+(B)$ we have $\lambda < |B| \Rightarrow \lambda^+ \nrightarrow
(\text{Depth}^+(B))^2_\sigma$.
\enddemo
\bigskip

\remark{\stag{3.5} Remark}  1) Instead $\lambda = d(B)^{< \kappa}$ we can use
$\lambda = \lambda^{< \kappa}$ such that:
\mr
\item "{$(*)$}"  there are no $\gamma^* < \kappa$ and $b_\gamma > 0$
(for $\gamma < \gamma^*$) and $\langle a_\alpha:\alpha < \lambda^+ \rangle$
such that for every $\alpha < \beta < \lambda^+$, for some $\gamma$ in
$B \restriction b_\gamma$ we have $a_\gamma \cap a_\alpha \ne b_\gamma 
\cap a_\beta$ are disjoint or have disjoint compliment. 
\ermn
2) Can think of parallel replacing in $(*)_2$ of \scite{3.1}(1), 2 by $n$,
that is: let $2 \le n < \omega$.  We say $\langle a_\beta:\beta < \beta^*
\rangle$ is $n$-semi-independent (in the Boolean Algebra $B$) if some $I$
witnesses it which means $a_\beta \in B$ are pairwise distinct and if for
$\alpha \le \gamma_0 < \ldots < \gamma_{n-1}$ and $b \in \langle a_\beta:
\beta < \alpha \rangle_B$ we have $(*)_1, (*)_3, (*)_4$ and
\mr
\item "{$(*)_{2,n}$}"  in $\{a_{\gamma_\ell} \cap b:\ell < n\}$ is
independent in $B \restriction b$.
\ermn
There are some variants and I have not tried if this gives something
interesting.
\endremark
\newpage

\head {\S4 On omitting cardinals by compact spaces} \endhead  \resetall
\bn
We continue Juhasz Shelah \cite{JuSh:612}.  We investigate what homomorphic
images some Boolean Algebras may have, and (in \scite{4.15}) prove the
topological analog of \S2, showing the existence of some subspaces for
Hausdorff spaces (not, necessarily compact).
\definition{\stag{4.1} Definition}  1) $\bold U_\theta(\mu) = \text{ Min}\{
|{\Cal P}|:{\Cal P} \subseteq [\mu]^{\le \theta}$ and $(\forall X \in 
[\mu]^\theta)(\exists a \in {\Cal P})$ \nl

$\qquad \qquad \qquad \qquad \qquad \qquad \qquad \qquad (|a \cap x) 
= \theta)\}$. \nl
2) Let ${\frak a}_\theta(\mu) = \text{ Min}\{|{\Cal A}|:{\Cal A} \subseteq
[\mu]^\theta$ is $\theta$-MAD$\}$ where ${\Cal A} \subseteq [\mu]^\theta$
is called $\theta$-AD if $A \ne B \in {\Cal A} \Rightarrow |A \cap B| <
\theta$ and we say ${\Cal A}$ is $\theta$-MAD means that in addition 
${\Cal A}$ is maximal under those restrictions.
\enddefinition
\bigskip

\remark{\stag{4.2} Remark}  In the case $\mu \ge 2^\theta$, in which we are
interested, $\bold U_\theta(\mu) = {\frak a}_\theta(\mu) = |{\Cal A}|$ 
whenever ${\Cal A}$ is $\theta$-MAD for $\mu$.
(See more and connection of pcf theory \cite{Sh:506}, \nl
\cite{Sh:589}, but we do not use any non-trivial fact.)
\endremark
\bigskip

\definition{\stag{4.3} Definition}  1) For $J_1 \subseteq J_2$ ideals of a
Boolean Algebra ${\Cal B}$, we say $J_2$ is $\theta$-full over $J_1$ inside
${\Cal B}$, \ub{if} for every $X \in [J_1]^\theta$ there is $b \in J_2$ such
that \nl
$|\{x \in X:{\Cal B} \models x \le b\}| = \theta$. \nl
2) The ideal $J$ of ${\Cal B}$ is $\theta$-full if: $J$ is $\theta$-full over
$J$ inside ${\Cal B}$. \nl
3) $J$ is $\tau$-local inside ${\Cal B}$, \ub{if} $|\{b \in {\Cal B}:{\Cal B}
\models b \le x\}| \le \tau$ for $x \in J$. \nl
4) In part (1), $J_2$ is strongly $\theta$-full over $J_1$ inside
${\Cal B}$ if for every $X \in [J_1]^\theta$ for some $Y \in [X]^\theta$,
for every $Z \subseteq Y$ we have $\dbcu_{b \in Z} a \in {\Cal B}$ in the 
sense of ${\Cal B}$ exist.  Similarly in part (2).
\enddefinition
\bn
\ub{\stag{4.4} Fact}:  1) If $J_2$ is $\theta$-full over $J_1$ inside
${\Cal B},h$ is a homomorphism from ${\Cal B}$ onto ${\Cal B}^*$ and
$J^*_\ell =
h(J_\ell)$ for $\ell=1,2$, \ub{then} $J^*_2$ is $\theta$-full over $J^*_1$
inside ${\Cal B}^*$. \nl
2) If $J_2$ is $\theta$-full over $J_1$ inside ${\Cal B}$ and $J_2$ is
$\tau$-local inside ${\Cal B}$, \ub{then} $\bold U_\theta(|J_1|) \le
\tau^\theta + |J_2|$. \nl
3) If the ideal $J$ of ${\Cal B}$ is $\theta$-full and $\tau$-local inside
${\Cal B}$, then $\bold U_\theta(|J|) \le \tau^\theta + |J|$. \nl
4) In parts 2) and 3), if $\tau = \theta$ we can replace $\tau^\theta$ by
$\tau$.
\bigskip

\demo{Proof}  1) Trivial. \nl
2) Let $J_1 = \{a_i:i < |J_1|\}$, let ${\Cal P}_y = \{X \subseteq |J_1|:|X|
= \theta$ and $(\forall i \in X)({\Cal B} \models a_i \le y)\}$ for each
$y \in J_2$, so $|{\Cal P}_y| \le \tau^\theta$.  Lastly, let ${\Cal P} = 
\cup\{{\Cal P}_y:y\in J_2\}$ so $|{\Cal P}| \le |J_2| \times
\underset {y \in J_2} {}\to \sup |{\Cal P}_y| \le \tau^\theta + |J_2|$.
Easily ${\Cal P}$ is as required in Definition \scite{4.1}. \nl
3) Follows by part (2).  \nl
4) Similar to the proof of part (2) only now
${\Cal P}_y = \{\{i < |J_1|:{\Cal B} \models a_i \le y\}\}$.
\hfill$\square_{\scite{4.4}}$
\enddemo
\bn
\ub{\stag{4.6} Fact}:  1) Assume $\lambda < \kappa \le \mu \le \kappa^\lambda$
and $\Theta \subseteq \Theta_{\mu,\lambda} =: 
\{ \theta \le \lambda:\bold U_\theta(\mu) = \mu\}$ and let $\sigma \in
\{\sigma:\sigma = \text{ cf}(\sigma) \le \lambda^+$ and for every
$\theta \in \Theta$ we have cf$(\theta) \ne \sigma\}$, (clearly there is one:
$\sigma = \lambda^+$).  \ub{Then} there is a Boolean algebra ${\Cal B}$ such 
that
\mr
\item "{$(a)$}"  $|{\Cal B}| = \mu$
\sn
\item "{$(b)$}"  ${\Cal B}$ is atomic with exactly $\kappa$ atoms
\sn
\item "{$(c)$}"  ${\Cal B}$ has a maximal ideal $J$ which is 
$2^\lambda$-local, moreover $x \in J \Rightarrow |\{y:{\Cal B} \models ``y \le
x \and y$ an atom (of ${\Cal B})"\}| \le \lambda$
\sn
\item "{$(d)$}"  for every $\theta \in \Theta,J$ is $\theta$-full (inside
${\Cal B}$)
\sn
\item "{$(e)$}"  if $2^\lambda \le \mu$ then 
${\Cal P}(\lambda)$ is isomorphic to some ${\Cal B} 
\restriction \{x:x \le a\}$ \nl
(so can demand $\bigl\{\{\alpha\}:\alpha <
\kappa \bigr\} \subseteq J \subseteq {\Cal B} \subseteq {\Cal P}(\kappa),
J \subseteq [\kappa]^{\le \lambda}$). \nl
If $2^\lambda > \mu,{\Cal B}_0 \subseteq {\Cal P}(\lambda)$ has cardinality
$\le \mu$ we can replace above ${\Cal P}(\lambda)$ by ${\Cal B}_0$.
\ermn
2)  We can add in part (2)
\mr
\item "{$(f)$}"  if $x \in J$ and $a_i \le x$ is an atom of ${\Cal B}$ for
$i < i^* \le \lambda$, \ub{then} some 
$y \in {\Cal B}$ is $\dbcu_{i < i^*} a_i$ in ${\Cal B}$'s sense
\sn
\item "{$(g)$}"  $2^\lambda \le \mu \and \theta \in \Theta \Rightarrow J$ 
is strongly $\theta$-full inside ${\Cal B}$, in fact $2^\lambda \le \mu
\and x \in J \Rightarrow {\Cal P}(x) \subseteq J$.
\endroster
\bigskip

\demo{Proof}  1)  Below if $2^\lambda > \mu$ we should replace everywhere
``$J_\zeta$ is the ideal of ${\Cal P}(\kappa)$ such that $\ldots$" by 
``$J_\zeta$ is the Boolean subring of ${\Cal P}(\kappa)$ such that"; 
if $2^\lambda \le \mu$ we may or may not.
As $\mu \le \kappa^\lambda$ we can find pairwise distinct
$x_i \in [\kappa]^\lambda$ for $i < \mu,x_0 = \{i:i < \lambda\}$.  Let $J_0$
be the ideal of ${\Cal P}(\kappa)$ generated by $\{x_i:i < \mu\} \cup
\bigl\{ \{\alpha\}:\alpha < \kappa \bigr\} \cup {\Cal B}_0$, so 
$J_0 \subseteq [\kappa]^{\le \lambda},|J_0| = \mu$.  We choose by 
induction on $\zeta \le \sigma$,
an ideal $J_\zeta$ of ${\Cal P}(\kappa),J_\zeta \subseteq 
[\kappa]^{\le \lambda},|J_\zeta| = \mu,J_\zeta$ is increasing continuous in
$\zeta$.  For $\zeta = 0,J_0$ was defined, for $\zeta$ limit let $J_\zeta =
\dbcu_{\varepsilon < \zeta} J_\varepsilon$.  For $\zeta = \varepsilon +1$,
let $J_\varepsilon = \{a^\varepsilon_i:i < \mu\}$, for $\theta \in \Theta$
let ${\Cal P}^\varepsilon_\theta \subseteq [\mu]^{\le \theta}$ exemplifies
$\bold U_\theta(\mu) = \mu$ (which follows from $\theta \in \Theta$), so
$|{\Cal P}^\varepsilon_\theta| \le \mu$, and let $J_\zeta$ be the ideal of
${\Cal P}(\kappa)$ generated by $J_\varepsilon \cup \{\dbcu_{i \in x}
a^\varepsilon_i:x \in {\Cal P}^\varepsilon_\theta$ for some $\theta \in
\Theta\}$, easy to check the inductive demand.  It is also easy to check
that then $J_{\zeta +1}$ is $\theta$-full over $J_\zeta$ (inside 
${\Cal P}(\kappa)$), when $\theta \in \Theta$.  Let $J = J_\sigma,
B = J \cup \{\kappa \backslash x:x \in J\} =$ the Boolean subalgebra of
${\Cal P}(\kappa)$ which $J$ generates.  Clearly $J$ is $\theta$-full 
over $J$ inside ${\Cal B}$.  (If $\sigma = \lambda^+$, if $X \subseteq J,
|X| \le \lambda$, then for some $\zeta < \lambda^+$ we have $X \subseteq
J_\zeta$, so $|X| \in \Theta \Rightarrow (\exists y \in J_{\zeta +1})
[|X| = |\{x \in X:B \models x \in y\}|]$.
If $\sigma < \lambda^+$ and $|X| \in \Theta$, then for some $\zeta < \sigma,
|X \cap J_\zeta| = \theta$ and proved as above.)  
So easily ${\Cal B},J$ are as required. \nl
2) Should be clear.  \hfill$\square_{\scite{4.6}}$
\enddemo
\bigskip

\proclaim{\stag{4.7} Claim}  1) If ${\Cal B},J$ (and $\kappa,\mu,\Theta)$ 
are as in fact \scite{4.6}(2), and ${\Cal B}^*$ is a homomorphic image of 
${\Cal B}$ and $\|{\Cal B}^*\| \ge 2^\lambda$, \ub{then} 
$\theta \in \Theta \Rightarrow \bold U_\theta(\|{\Cal B}^*\|) = 
\|{\Cal B}^*\|$. \nl
2) Hence if $\theta \in \Theta,2^\lambda \le \chi < \kappa,(\forall \alpha
< \chi)(|\alpha|^{< \theta} < \chi \and \text{{\rm cf\/}}(\chi) = 
\theta)$, then $\|{\Cal B}^*\| \notin [\chi,\chi^\theta)$. \nl
3) Also it follows that the number of ultrafilters of ${\Cal B}^*$ is
$\le 2^{2^\lambda} + \|{\Cal B}^*\|$, if $\|{\Cal B}^*\| > 2^\lambda \and
\lambda \in \Theta$ equality holds (in any case $\ge \|{\Cal B}^*\|$).
\endproclaim
\bigskip

\demo{Proof}  Let $h^*:{\Cal B} \rightarrow {\Cal B}^*$ be a homomorphism from
${\Cal B}$ onto ${\Cal B}^*$, let $J^* = \{h^*(x):x \in J\}$.  First assume
$1_{{\Cal B}^*} \in J^*$, this means that for some $x \in J,h^*(x) = 
1_{{\Cal B}^*}$, so ${\Cal B} \restriction x =: {\Cal B} \restriction \{y:
{\Cal B} \models y \le x \}$ has ${\Cal B}^*$ as a homomorphic image so
$\|{\Cal B}^*\| \le |{\Cal P}(x)| \le 2^\lambda$, and the number of 
ultrafilters of ${\Cal B}^*$ is $\le 2^{2^\lambda}$; also if
$\|{\Cal B}^*\| \ge 2^\lambda$ we get $\|{\Cal B}^*\| = 2^\lambda$ hence
${\frak a}_\theta(\|{\Cal B}^*\|) = {\frak a}_\theta(2^\lambda) = 2^\lambda$
as $(2^\lambda)^\theta = 2^\lambda$; this finishes except the other 
inequalities in part (3).  So assume $1_{{\Cal B}^*} \notin J^*$,
hence $J^*$ is a maximal ideal of ${\Cal B}^*$, also clearly $J^*$ is
$\theta$-full inside ${\Cal B}^*$ for every $\theta \in \Theta$ (by Fact
\scite{4.4}(1)).  As $J$ is $2^\lambda$-local (see Fact \scite{4.6}(c)), 
clearly $J^*$ is $2^\lambda$-local, hence by \scite{4.4}(2), (letting 
$\chi =: |J^*|$ we have), ${\frak a}_\theta(\chi) \le 
(2^\lambda)^\theta + \chi$, but ${\frak a}_\theta
(\chi) \ge \chi \ge 2^\lambda = (2^\lambda)^\theta$ so ${\frak a}_\theta
(\chi) = \chi$.  Also $\|{\Cal B}^*\| = \chi + \chi = \chi$.  So we have
gotten the conclusion of \scite{4.7}(1).  Now \scite{4.7}(2) 
follows easily as ${\frak a}_\theta
(\chi) \ge \chi^\theta$ as $\bigl\{ \{\eta \restriction \alpha:\alpha <
\theta\}:\eta \in \chi^\theta \bigr\}$ is a $\theta$-AD family of subsets of
$\{\eta:\eta \in {}^{\theta >} \chi\}$ which has cardinality $\chi^\theta$.  
Lastly, for \scite{4.7}(3), if $D$ is an ultrafilter of 
${\Cal B}^*$ then either $D = {\Cal B}^* \backslash J^*$ or for some
$x \in J^*,x \in D$ but for each $x \in J^*,{\Cal B}^* \restriction
\{y \in {\Cal B}^*:y \le x\}$ has $\le 2^\lambda$ members so the number of
ultrafilters of ${\Cal B}$ to which $x$ belongs is $\le 2^{2^\lambda}$,
that means \nl
\mn
\ub{Fact}:  If the Boolean Algebra ${\Cal B}$ has a $\tau$-local maximal 
ideal, then \nl
$|(\text{set of ultrafilters of } {\Cal B}^*)| \le 2^\tau +
\|{\Cal B}^*\|$. 
\sn
Last point is the second inequality in part (3); 
assume $\mu = \|{\Cal B}^*\| > 2^\lambda$.  Let $x_i \in J$ for
$i < \mu$ be such that $i < j < \tau \Rightarrow h^*(x_\ell) \ne h^*(x_j)$
(possible as $\|{\Cal B}^*\| = |J^*|$).  So by the $\triangle$-system
argument without loss of generality for some 
$x^*,i < j < (2^\lambda)^+ \Rightarrow x_i \cap x_j = x^*$.  But 
$\langle h^*(x_i) \cap h^*(x^*):i < (2^\lambda)^+ \rangle$ is constant
hence $\langle h^*(x_i - x^*):i < (2^\lambda)^+ \rangle$ is a sequence of
$(2^\lambda)^+$ pairwise disjoint non-zero (in the sense of ${\Cal B}^*$)
members of ${\Cal B}^*$.  We can find $y \in {\Cal B}$ such that $w =
\{i < (2^\lambda)^+:x_i - x^* \le y\}$ has cardinality $\lambda$ (remember
$\lambda \in \Theta$ by assumption of \scite{4.7}(3)), hence
$(\forall u \subseteq w)(\exists z \in {\Cal B})([i \in u \rightarrow x_i -
x^* \le z] \and [i \in w \backslash u \Rightarrow (x_i,x^*) \cap z = 0])$
hence ${\Cal B}^*$ has a homomorphic image isomorphic to ${\Cal P}(\lambda)$
hence ${\Cal B}$ has $\ge 2^{2^\lambda}$ ultrafilters. \nl
Together we finish. \hfill$\square_{\scite{4.7}}$
\bigskip

\centerline {$* \qquad * \qquad *$}
\bigskip

By claims \scite{4.6}, \scite{4.7} we really finish.  
Let me point some specific conclusions: conclusion \scite{4.8} is the 
theorem of Juhasz Shelah \cite{JuSh:612}.
\bigskip

\demo{\stag{4.8} Conclusion}  For every $\kappa > 2^\lambda$, there is a
Boolean algebra ${\Cal B}_\kappa$ such that: ${\Cal B}_\kappa$ is atomic with
$\kappa$ atoms, $\|{\Cal B}\| = \kappa^\lambda$ and for every homomorphic
image ${\Cal B}^*$ of ${\Cal B}$ of cardinality $\chi > 2^\lambda$ we have
$\chi = \chi^\lambda$ and the number of ultrafilters of ${\Cal B}^*$ is
$2^{2^\lambda} + \chi$ in particular $\chi \in [\kappa,\kappa^\lambda)$ is
impossible.  (Check).
\enddemo
\bigskip

\demo{Proof}  We apply fact \scite{4.6} + claim \scite{4.7} to 
$\mu = \kappa^\lambda$ and our $\kappa$, so $\Theta = \{\theta:\theta \le
\lambda\}$.  So in $\Theta$ we get ${\Cal B}$.  Let ${\Cal B}^*$ be a
homomorphic image of ${\Cal B}$ (equivalently, a quotient of ${\Cal B}$) and
$\chi = \|{\Cal B}^*\| > 2^\lambda$.  So $\theta \in \Theta \Rightarrow
\bold U_\theta(\chi) = \chi$, now if $\chi < \chi^\lambda$ let $\sigma =
\text{ Min}\{\theta:\chi^\theta = \chi\}$, so $\sigma \in \Theta$ and
$\chi^{< \sigma} = \chi < \chi^\sigma$ and we get a contradiction to
$\bold U_\sigma(\chi) = \chi$.  For the number of ultrafilters use 
\scite{4.7}(3).  \hfill$\square_{\scite{4.8}}$
\enddemo
\bigskip

\demo{\stag{4.9} Conclusion}  If $\lambda$ is strong limit singular (e.g.
$\beth_\omega$) and $2^\lambda < \kappa \le \mu \le \kappa^\lambda$,
\ub{then} there is an (atomic) Boolean Algebra ${\Cal B}$ with $\kappa$ atoms,
$|{\Cal B}| = \mu$ such that: for every large enough regular
$\theta < \lambda$ we have:
\mr
\item "{{}}"  every homomorphic image ${\Cal B}^*$ of ${\Cal B}$ of
cardinality $> (2^\lambda)^+$ satisfies ${\frak a}_\theta[\|{\Cal B}^*\|] =
\|{\Cal B}^*\|$ \nl
(so for any cardinality $\tau$, we have
$\|{\Cal B}^*\| \in [\tau^{< \theta},\tau^\theta)$ is impossible).
\endroster
\enddemo
\bigskip

\demo{Proof}  By conclusion \scite{4.6} + \scite{4.7} using \cite{Sh:460}.
\enddemo
\bigskip

\remark{\stag{4.10} Remark}  In addition to $\bold U_\theta(-)$ we can use
other functions (e.g. as in \cite[\S1]{Sh:589}, even if their number is
$> \mu$ it does not matter as for each $\chi \in (2^\lambda,\mu)$ we can
choose one) but does not seem worth elaborating.
\endremark
\bigskip

\remark{\stag{4.11} Remark}  1) Assume $\kappa$ is strong limit singular of
cofinality $\theta^* < \lambda < \kappa$ and $2^\kappa = \kappa^\lambda >
\kappa$.  There are many $\mu \in [\kappa,2^\kappa)$ such that
$\Theta_{\mu,\lambda} = \{\theta \le \lambda:\text{cf}(\theta) \ne 
\theta^*\}$.  E.g. $\mu \in \{\kappa^{+ n}:n <\omega\}$, also (see
\cite[Ch.IX]{Sh:g}) for a club $E$ of $\lambda^{+4},\delta \in E \and
\text{cf}(\delta) \ge \lambda^+ \and n < \omega \Rightarrow 
\Theta_{\kappa^{+ \delta +n},\lambda} = \{\theta \le \lambda:\text{cf}(\theta)
\ne \theta^*\}$. 
\endremark
\bigskip

\proclaim{\stag{4.12} Claim}  In Claim \scite{4.7} we can add (see
\cite{Sh:460})
\mr
\item "{$(a)$}"  if $\theta \le \lambda$ and $\mu = \mu^{[\theta]}$ then
$\|{\Cal B}^*\|^{[\theta]} = \|{\Cal B}^*\|$
\sn
\item "{$(b)$}"  similar more general condition (see \cite[\S1,end]{Sh:589}).
\endroster
\endproclaim
\bigskip

\demo{\stag{4.13} Conclusion}  If $\mu$ is strongly limit singular and
cf$(\mu) \le \theta < \mu$, \ub{then} for some atomic Boolean Algebra
${\Cal B}$ we have:
\mr
\item "{$(a)$}"  ${\Cal B}$ has cardinality $\mu$
\sn
\item "{$(b)$}"  ${\Cal B}$ has $\mu$ ultrafilters
\sn
\item "{$(c)$}"  if ${\Cal B}^*$ is a homomorphic image of ${\Cal B}$ of
singular strong limit cardinality $\chi > \theta$, then cf$(\chi) =
\text{ cf}(\mu)$ and ${\Cal B}^*$ has $\chi$ ultrafilters
\sn
\item "{$(d)$}"  if ${\Cal B}^*$ is a homomorphic image of ${\Cal B},
\chi = \|{\Cal B}^*\| > 2^\lambda$, \ub{then} \nl
$\chi^{< \text{cf}(\mu)} = \chi$ and $\chi = \text{ cov}(\chi,\theta^+,
\theta^+,(\text{cf}(\mu))^+)$.
\endroster
\enddemo 
\bigskip

\demo{Proof}  By \scite{4.7} and \scite{4.12}.
\enddemo
\bn
Another example
\demo{\stag{4.14} Conclusion}  If $\mu$ is strong limit singular,
$\theta = \text{ cf}(\delta) < \text{ cf}(\mu),\mu^{+ \delta} < 2^\mu,
(\mu^{+ \delta})^{< \theta} = \mu^{+ \delta}$, \ub{then} for some Boolean
Algebra ${\Cal B}$ of cardinality $\mu^{+ \delta}$ it has $\mu^{+ \delta}$
ultrafilters, and for every homomorphic image ${\Cal B}^*$ of ${\Cal B}$ of
cardinality $\chi,2^{\text{cf}(\mu)} < \chi < \mu$ we have:
\mr
\item "{$(*)$}"  if $\chi$ is strong limit then cf$(\chi) \in
\{\text{cf}(\mu),\text{cf}(\delta)\}$.
\endroster
\enddemo
\bigskip

\proclaim{\stag{4.15} Claim}  Assume
\mr
\item "{$(a)$}"  $\lambda$ is strong limit singular, $\kappa = \lambda,
\kappa = \text{{\rm cf\/}}(\mu) < \mu \le \lambda$
\sn
\item "{$(b)$}"  $X$ is a Hausdorf topological space with $w(X) = \lambda$.
\ermn
\ub{Then} $X$ has a closed subset $Y$ such that:
\mr
\item "{$(\alpha)$}"  $\mu \le w(B) \le 2^{< \mu}$
\sn
\item "{$(\beta)$}"  $\mu \le |Y| \le \dsize \sum_{\theta < \mu}
2^{2^\theta}$.
\endroster
\endproclaim
\bigskip

\remark{Remark}  1) If we speak on Boolean Algebras ${\Cal B}$, let $X =
\text{ Ult}({\Cal B})$, so $w(X) = |{\Cal B}|$ and $\{Y:Y \subseteq X
\text{ closed}\} = \{\text{Ult}({\Cal B}/I):I \text{ and ideal of }{\Cal B}\}$
essentially. \nl
2) So this result reasonably compliments Juhasz \cite{Ju1}, Juhasz Shelah 
\cite{JuSh:612}.
\endremark
\bigskip

\demo{Proof}  \ub{Case 1}:  $\mu = \lambda$.

Let $W = \{U_i:i < \lambda\}$ be a basis of $X$.  Choose $\langle
\lambda_i:i < \kappa \rangle$ be increasing continuously with limit
$\lambda,\lambda_0 = 0,(\forall \sigma < \lambda_{i+1})[\sigma^{\lambda_i}
< \lambda_{i+1} = \text{ cf}(\lambda_{i+1})]$.  As $|X| \ge \lambda$ (as
$w(X) = \lambda$) necessarily $s^+(X) > \lambda$, see Juhasz \cite{Ju} so
there is $\{y_\alpha:\alpha < \lambda\} \subseteq X$, (with no repetitions)
which is discrete.  We can choose $Z_i \in [\lambda_{i+1}]^{\lambda_{i+1}}$
such that

$$
(\forall \alpha < \lambda_i)(\dsize \bigwedge_{\zeta \in Z_i} 
y_\zeta \in U_\alpha \vee \dsize \bigwedge_{\zeta \in Z_i} y_\zeta 
\notin U_\alpha).
$$
\mn
By renaming \wilog \, $Z_i = [\lambda_i,\lambda_{i+1})$.  Let $Y = c \ell
\{y_\alpha:\alpha < \lambda\}$.  It suffices to prove that $|Y^*| = \lambda$,
for this it suffices to prove
\mr
\item "{$(*)$}"  if $x \in Y^*$, then for some $i < \kappa$ we have
$$
x \in c \ell(\{y_\alpha:\alpha < \lambda_i\} \cup \{y_{\lambda_j}:j <
\kappa\}).
$$
\ermn
If $x$ contradicts $(*)$, then for every $i < \lambda$ there is $\alpha_i
< \lambda$ such that $x \in U_{\alpha_i},U_{\alpha_i} \cap (\{y_\alpha:
\alpha < \lambda_i\} \cup \{y_{\lambda_j}:j < \kappa\}) = \emptyset$. \nl
Now $\alpha_i < \lambda = \dbcu_{j < \kappa} \lambda_j$ so for some $j,
\alpha_i < \lambda_j$, so $U_{\alpha_i} \cap U_{\alpha_j}$ is an open 
neighborhood of $x$ (as 
$U_{\alpha_i},U_{\alpha_j}$ are) disjoint to $\{y_\alpha:\alpha <
\lambda_j\}$ as $U_{\alpha_j}$ is, and for each
$\beta \in [\lambda_j,\lambda)$ we have
\mr
\item "{$(*)$}"  for some $\zeta \in [j,\kappa),\lambda_\zeta \le \beta
< \lambda_{\zeta +1}$ and so (as $\alpha_i < \lambda_j,j \le \zeta$) we have
$y_\beta \in U_{\alpha_i} \Leftrightarrow y_{\lambda_\zeta} \in U_{\alpha_i}$;
but $y_{\lambda_\zeta} \notin U_{\alpha_i}$ by the choice of $U_{\alpha_i}$,
hence
$$
\beta \in [\lambda_j,\lambda] \Rightarrow y_\beta \notin U_{\alpha_i} 
\Rightarrow y_\beta \notin U_{\alpha_i} \cap U_{\alpha_j}.
$$
\ermn
So $U_{\alpha_i} \cap U_{\alpha_j}$ is an open neighborhood of $x$ disjoint to
$\{y_\alpha:\alpha < \lambda\}$ so $x \notin Y^*$, contradiction so $(*)$
holds, hence we are done.
\bn
\ub{Case 2}:  $\mu < \lambda$.

Let $\mu = \dsize \sum_{i < \kappa} \mu_i,\mu_i$ strictly increasing with
$i$.  Repeat the above and let \nl
$Y' = c \ell\{y_\alpha:\alpha \in [\lambda_i,\lambda_i + \mu_i)$ 
for some $i\}$.
\nl
So

$$
\align
y \in Y' &\Rightarrow y \in c \ell\{y_\alpha:\alpha \in \dbcu_{i < \kappa}
[\lambda_i,\lambda_i + \mu_i)\} \\
  &\Rightarrow \dsize \bigvee_{j < \kappa} x \in c \ell(\{y_\alpha:\alpha \in 
\dbcu_{i < j}[\lambda_i,\lambda_i + \mu_i)\} \cup \{y_{\lambda_i}:
i < \kappa\}).
\endalign
$$
\mn
So clearly $\mu = \dsize \sum_{i < \kappa} \mu_i \le |Y'| \le \dsize \sum
_{i < \kappa} 2^{2^{\mu_i}}$. \nl
Now a basis of $Y'$ is $\dbcu_{\Sb j < \kappa \\ A \subseteq \kappa \endSb}
W_{j,A}$ where

$$
\align
W_{j,A} = \biggl\{ U_\alpha \cap Y':&\alpha < \lambda_j \text{ and } 
[i \in [j,\kappa) \cap A \Rightarrow \dsize \bigwedge_{\alpha \in 
[\lambda_i,\lambda_i + \mu_i)} y_\alpha \in U_\alpha] \\
  &\text{and } [i \in [j,\kappa) \backslash A \Rightarrow
\dsize \bigvee_{\alpha \in [\lambda_i,\lambda_i + \mu_i)} y_\alpha \notin
U_\alpha] \biggr\}.
\endalign
$$
\mn
So considering $W_{j,\kappa}$ we can look at the space
$Y_j = c \ell(\{y_\alpha:\alpha < \dbcu_{i < j} [\lambda_i,\lambda_i +
\mu_i)\} \cup \{y_{\lambda_i}:i < \kappa\})$, which has density
$\le \mu_j + \kappa_i = \mu_i$ hence witness weight $\le 2^{\mu_j}$, so we
can finish easily. \hfill$\square_{\scite{4.15}}$
\enddemo
\bn
We may consider this in the framework
\definition{\stag{4.16} Definition}  1) For a Boolean Algebra ${\Cal B}$
let

$$
\text{wcSp}(B) = \{(\|{\Cal B}'\|,\text{ult}({\Cal B}')):{\Cal B}'
\text{ is an infinite homomorphic image of } {\Cal B}\}
$$
\mn
(remember

$$
\text{ult}({\Cal B}') = |\text{Ult}({\Cal B}')|,\text{Ult}(B') =
\{D:D \text{ an ultrafilter of } {\Cal B}'\}).
$$
\mn
2) For a topological space $X$ let

$$
\text{wcSp}(X) = \{(|Y|,w(Y)):Y \text{ is a closed subspace of } X\}.
$$
\mn
(remember $w(X)$ is the weight of the topological space $X$).
\enddefinition
\bigskip

\remark{\stag{4.17} Remark}  Of course, we can use disjoint 
sums of Boolean Algebras to get more examples (and similarly for 
topological spaces) as

$$
\align
\text{wcSp}(\dsize \sum_{i < \alpha} B_i) = \biggl\{ (\dsize \sum_{i < \alpha}
\lambda_i,1 + \dsize \sum_{i < \alpha} \mu_i):&(\lambda_i,\mu_i) \in
\text{wcSp}(B_i) \\
  &\cup \{(2^n,2^{2^n}):n < \omega \cup \{(0,0) \\
  &\text{for } i < \alpha \text{ and } \dsize \sum_{i < \alpha}
\lambda_i \text{ is infinite} \biggr\}.
\endalign
$$
\mn
In this way we can get more examples from the ones from \cite{JuSh:612}, but
this does not cover all the above. 
\endremark
\enddemo
\newpage

\head {\S5 Depth of Ultraproducts of Boolean Algebras} \endhead  \resetall \footnote""{Done 7-8/97}
\bigskip

\proclaim{\stag{5.1} Claim}  Assume $\square_\lambda$ (i.e. there is
$\langle C_\delta:\delta < \lambda^+$ limit$\rangle$ such that $C_\delta$ is
a club of $\delta$ of order type $< \delta$ if $\text{{\rm cf\/}}(\delta) <
\delta$ and $\delta_1 \in \text{{\rm acc\/}}(C_{\delta_2}) \Rightarrow
C_{\delta_1} = C_{\delta_2} \cap \delta_1$). \nl
Let $\kappa = \text{{\rm cf\/}}(\kappa) < \lambda$.  \ub{Then} there are 
Boolean Algebras $B_\varepsilon$ for $\varepsilon < \kappa$ such that:
\mr
\item "{$(a)$}"  Depth$(B_\varepsilon) \le \lambda$
\sn
\item "{$(b)$}"  for any uniform ultrafilter $D$ on $\kappa,\lambda^+ \le
\text{ Depth}(\dsize \prod_{\varepsilon < \kappa} B_\varepsilon/D)$.
\endroster
\endproclaim
\bigskip

\remark{\stag{5.2} Remark}  This can be expressed through \S1, see later.
\endremark
\bigskip

\demo{Proof}  Let $\langle C_\delta:\delta < \lambda^+ \text{ limit}
\rangle$ exemplify $\square_\lambda$ so there are an ordinal $\gamma^*$ and a 
stationary $S \subseteq \lambda^+$ such that $(\forall \alpha \in S)
[\text{otp}(C_\alpha) =
\gamma^*]$ (so $\alpha \in S \Rightarrow \alpha$ limit and), cf$(\gamma^*) =
\kappa$ and $\gamma^*$ divisible by $\omega^2$.  So without loss of generality
for every $\delta,C_\delta \cap S = \emptyset$ (by deleting the first
$\gamma^* +1$ elements from any $C_\delta$ of greater order type) also
\wilog \, $[\alpha \in C_\alpha \backslash
\text{ acc}(C_\alpha) \Rightarrow \alpha \text{ non-limit}]$.  
\enddemo
\bn
\ub{\stag{5.1A} Fact}:  Under the assumptions of \scite{5.1} there are sets
$A_{\alpha,\varepsilon}(\alpha < \lambda^+,\varepsilon < \kappa)$ such that:
\mr
\widestnumber\item{$(iii)$}
\item "{$(i)$}"  $\varepsilon < \zeta \Rightarrow A_{\alpha,\varepsilon}
\subseteq A_{\alpha,\zeta}$
\sn
\item "{$(ii)$}"  $\dbcu_{\varepsilon < \kappa} A_{\alpha,\varepsilon} =
\alpha$
\sn
\item "{$(iii)$}"  $\beta \in A_{\alpha,\varepsilon} \Rightarrow
A_{\beta,\varepsilon} = A_{\alpha,\varepsilon} \cap \beta$
\sn
\item "{$(iv)$}"  $\alpha \in S \and \varepsilon < \kappa \Rightarrow \sup
(A_{\alpha,\varepsilon}) < \alpha$
\sn
\item "{$(v)$}"  $A_{\alpha,\varepsilon}$ is a closed subset of $\alpha$
and disjoint to $S$
\sn
\item "{$(vi)$}"  if $\beta \in \text{ acc}(C_\alpha)$ then $\beta \in
A_{\alpha,0}$
\sn
\item "{$(vii)$}"  $\beta \in A_{\beta +1,0}$.
\endroster
\bigskip

\demo{Proof of \scite{5.1A}}  We choose by induction $\alpha < \lambda^+$ 
a sequence $\langle A_{\alpha,\varepsilon}:
\varepsilon < \kappa \rangle$ such that clauses (i)-(vii) holds. \nl
Let $\langle \gamma_\varepsilon:\varepsilon < \kappa \rangle$ be increasing
continuous in $\varepsilon$ sequence of ordinals with limit $\gamma^*$, 
each $\gamma_\varepsilon$ a limit ordinal. \nl
How do we carry the definition?

\bn
\ub{Case 1}:  $\alpha = 0$.

Let $A_{\alpha,\varepsilon} = \emptyset$.
\bn
\ub{Case 2}:  $\alpha = \beta +1$.

Let $A_{\alpha,\varepsilon} = A_{\beta,\varepsilon} \cup \{\beta\}$.
\bn
\ub{Case 3}:  $\alpha$ limit, $\alpha > \sup(\text{acc}(C_\alpha))$.

So necessarily cf$(\alpha) = \aleph_0$.  Let $\beta_0 \in C_\alpha$ and
$\beta_0 = \text{ max}(\text{acc } C_\alpha)$ if acc$(C_\alpha) \ne 
\emptyset$.  Choose $\beta_n$ (for $n \in [1,\omega))$ such that 
$n \ge 0 \Rightarrow \beta_n < \beta_{n+1},\beta_n \in C_\alpha,
\alpha = \dbcu_{n < \omega} \beta_n$.  Choose $\varepsilon_n < \kappa$ 
such that $\varepsilon_0 = 0,\varepsilon_n \le \varepsilon_{n+1},
\beta_n \in A_{\beta_{n+1},\varepsilon_{n+1}}$.  Lastly, for $\varepsilon
< \kappa$ we let $A_{\alpha,\varepsilon}$ be: $\bigcup\{\{\beta_n\} \cup
A_{\beta_n,\varepsilon}:n \text{ satisfies } \varepsilon_n \le \varepsilon\}$.
\nl
Now check.
\bn
\ub{Case 4}:   $\alpha$ limit, $\alpha = \sup(\text{acc }C_\alpha),\alpha
\notin S$. \nl
Let $A_{\alpha,\varepsilon} = \cup\{A_{\beta,\varepsilon}:\beta
\in \text{ acc}(C_\alpha)\}$. \nl
Remember that $C_\delta \cap S = \emptyset$ for every limit $\delta <
\lambda^+$.
\bn
\ub{Case 5}:  $\alpha$ limit, $\alpha = \sup(\text{acc } C_\alpha),
\alpha \in S$.

Let $A_{\alpha,\varepsilon} = A_{\beta_\varepsilon,\varepsilon}$ where
$\beta_\varepsilon$ is the $\gamma_\varepsilon$-th member of $C_\alpha$
(so necessarily $\beta_\varepsilon \in \text{acc}(C_\alpha)$ and $\xi <
\zeta \Rightarrow C_{\beta_\xi} = C_{\beta_\zeta} \cap \beta_\varepsilon
\Rightarrow (\forall \varepsilon)[(A_{\beta_\xi,\varepsilon} = 
A_{\beta_\zeta,\varepsilon} \cap \beta_\varepsilon)]$. \nl
Check. \hfill$\square_{\scite{5.1A}}$
\enddemo
\bigskip

\remark{\stag{5.3} Remark}  This is relevant 
to a problem from \cite{Sh:108}. \nl
\ub{Continuation of the proof of \scite{5.1}}. Let 
$<_\varepsilon$ be the following two place relation on $\lambda^+:
\alpha <_\varepsilon \beta \Rightarrow \alpha \in A_{\beta,\varepsilon}$.
It is a partial order (by clause (iii) of \scite{5.1A}).  
Also $\alpha < \beta \Rightarrow \dsize \bigvee_{\zeta < \kappa} \,\,
\dsize \bigwedge_{\varepsilon \in [\zeta,\kappa)} \alpha
<_\varepsilon \beta$ by clauses (i) + (ii) of \scite{5.1A}.  
Let $B_\varepsilon$ be $BA[(\lambda^+,<_\varepsilon)]$,
i.e. it is a Boolean Algebra generated by $\langle x_\alpha:\alpha <
\lambda^+ \rangle$ freely except
\mr
\item "{$\bigotimes$}"  $x_\alpha \le x_\beta$ when $\alpha <_\varepsilon
\beta$.
\ermn
Clearly if $D$ is a filter on $\kappa$ containing the co-bounded subsets of
$\kappa$ then Depth$^+(\dsize \prod_{i < \kappa} B_i/D) > \lambda^+$ as
$\left< <x_\alpha:\varepsilon < \kappa>/D:\alpha < \lambda^+ \right>$
exemplifies this.  
Assume toward contradiction Depth$^+(B_\varepsilon) > \lambda^+$ so assume
$\bar b = \langle b_\gamma:\gamma < \lambda^+ \rangle$ is (strictly) 
increasing in $B_\varepsilon$.  
Choose by induction on $\gamma < \lambda^+$ a model
$M_\gamma \prec ({\Cal H}(\lambda^{++}),\in,<^*_{\lambda^{++}})$ 
of cardinality $\lambda$ increasing continuous in $\gamma$ 
such that $\{\bar C,<_\varepsilon,\bar b\} \in M_0$ and 
$\langle M_\beta:\beta \le \gamma 
\rangle \in M_{\gamma +1}$.  So $C^* = \{\delta < \lambda^+:M_\delta \cap
\lambda^+ = \delta\}$ is a club of $\lambda^+$ so choose $\delta(*) \in S
\cap \text{ acc}(C^*)$.

Let $b_\gamma = \tau(x_{\alpha(\gamma,0)},\dotsc,
x_{\alpha(\gamma,n_\gamma-1)})$
with $\alpha(\gamma,0) < \alpha(\gamma,1) < \ldots \alpha(\gamma,n_\gamma-1)
< \lambda^+$, and for $Y \subseteq \lambda^+$ let $B_{\varepsilon,Y}$ be the
subalgebra of $B_\varepsilon$ generated by $\{x_\alpha:\alpha \in Y\}$.
Easily
\mr
\item "{$(*)$}"  $B_{\varepsilon,Y}$ is the algebra generated by $\{x_\alpha:
\alpha \in Y\}$ freely except $x_\alpha \le x_\beta$ when $\alpha
<_\varepsilon \beta \and \alpha \in Y \and \beta \in Y$, i.e. 
$B_{\varepsilon,Y} = BA[(Y,<_\varepsilon \restriction Y)]$.
\ermn
Clearly by clause (ii) for some $\zeta_\ell < \kappa$ we have
$\alpha(\delta(*),\varepsilon) > \delta(*) \Rightarrow \delta(*) \in
A_{\alpha(\delta(*),\ell),\zeta_\ell}$ and 
$\alpha(\delta(*),\varepsilon) < \delta(*)
\Rightarrow \alpha(\delta(*),\ell) \in A_{\delta(*),\zeta_\ell}$ and let
$\xi = \max\{\varepsilon\} \cup \{\zeta_\ell:\ell < n_{\delta(*)}\}$, so
$\delta(*) \cap \dbcu_{\ell < n_{\delta(*)}} A_{\alpha(\delta(*),\ell),
\varepsilon} \subseteq A_{\delta(*),\xi}$. \nl
Let $\alpha_0(*) = \sup(\delta(*) \cap
\dbcu_{\ell < n_{\delta(*)}} A_{\alpha(\delta(*),\ell),\varepsilon})$, 
now as $\delta(*) \in S,n_{\delta(*)} < \omega,A_{\delta(*),\xi}$ is 
closed and clause (iv) of \scite{5.1A} we have $\delta(*) > \sup
(A_{\delta(*),\xi})$ hence clearly 
$\alpha_0(*) < \delta(*)$.  Let $\alpha(*) =
\text{ Min}(C^* \backslash \alpha_0(*))$ so $\alpha(*) < \delta(*)$ (as
$\delta(*) \in \text{ acc } C^*)$.  Let $Y_0 = \alpha(*),Y_1 = \delta(*),
Y_2 = \alpha(*) \cup \{\alpha(\delta(*),0),\dotsc,\alpha(\delta(*),
n_{\delta(*)}-1)\}$.
\mn
Easily: $Y_0 = Y_1 \cap Y_2,Y_1 \cup Y_2 \subseteq \lambda^+$ and

$$
\beta_1 \in Y_1 \backslash Y_2 \and \beta_2 \in Y_2 \backslash Y_1
\Rightarrow (\beta_1,\beta_2 \text{ are } <_\varepsilon \text{-incomparable)}.
$$
\mn
By $(*)$ above and the definition of $C^*$ clearly:
\mr
\item "{$(a)$}"  $B_{\varepsilon,Y_1 \cup Y_2}$ is the free product of
$B_{\varepsilon,Y_1}$ and $B_{\varepsilon,Y_2}$ over
$B_{\varepsilon,Y_0}$
\sn
\item "{$(b)$}"  $b_{\delta(*)} \notin M_{\delta(*)},B_\varepsilon 
\restriction M_{\delta(*)} = B_{\varepsilon,Y_1}$ so $b_{\delta(*)} \notin
B_{\varepsilon,Y_1}$ hence $b_{\delta(*)} \notin B_{\varepsilon,Y_0}$
\sn
\item "{$(c)$}"  $b_{\alpha(*)} \notin M_{\alpha(*)},B_\varepsilon 
\restriction M_{\alpha(*)} = B_{\varepsilon,Y_0}$ so $b_{\alpha(*)} \notin
B_{\varepsilon,Y_0}$.
\ermn
But by the choice of $\bar b,B_\varepsilon \models 
b_{\alpha(*)} < b_{\delta(*)}$ hence 
$B_{\varepsilon,Y_1 \cup Y_2} \models b_{\alpha(*)} < b_{\delta(*)}$ hence
by clause (a) for some 
$c \in B_{\varepsilon,Y_0}$ we have $B_{\varepsilon,Y_1} \models
b_{\alpha(*)} \le c$ and 
$B_{\varepsilon,Y_2} \models c \le b_{\delta(*)}$.  So
$c \in M_{\alpha(*)}$ and $\alpha(*) \in Z = \{\alpha:b_\alpha \le c\}$,
this last set $Z$ is an initial segment of $\lambda^+$, it belongs to
$M_{\alpha(*)}$, hence its supremum belongs to $M_{\alpha(*)}$, but the
supremum $\in \{\alpha:\alpha \le \lambda^+\}$, and is not in $\lambda^+$
(as $M_{\alpha(*)} \cap \lambda^+ = \alpha(*)$ and $\alpha(*)$ belongs to
the set), so $Z = \lambda^+$.  So
$B_\varepsilon \models b_{\delta(*)} \le c$ but (by its choice)
$B_\varepsilon \models c \le b_{\delta(*)}$ so $c = b_{\delta(*)}$, but
$c \in M_{\alpha(*)} \prec M_{\delta(*)},b_{\delta(*)} \notin
M_{\delta(*)}$, contradiction.  \hfill$\square_{\scite{5.1}}$
\bn
\ub{\stag{5.4} Conclusion}:  Under the assumption of claim \scite{5.1}
\mr
\item "{$(a)$}"  we can find $\bold c:[\lambda]^2 \rightarrow \kappa$ such
that $\boxtimes_{\lambda,\lambda,\theta}$ (from \scite{1.1})
\sn
\item "{$(b)$}"  $NQs_2(\lambda,\mu,\kappa)$ (see Definition \scite{1.6}(2),
in fact $NQs_2(\lambda,\mu,A,I)$ where $A$ is the interval Boolean Algebra
of $\kappa,I = \{a \in A:\sup(a) < \kappa\}$.
\endroster
\bigskip

\demo{Proof}  Easy, e.g.
\mr
\item "{$(a)$}"  let $\left< \langle A_{\alpha,\varepsilon}:
\varepsilon < \kappa
\rangle:\alpha < \lambda^+ \right>$ be as in the proof of \scite{5.1}.
Let for $\alpha < \beta$
$$
\bold c\{\alpha,\beta\} =: \text{ Min}\{\varepsilon:\alpha \in
A_{\beta,\varepsilon}\}.
$$
\endroster
\enddemo
\bigskip

\remark{\stag{5.5} Remark}  For $\kappa$ singular we can deduce \scite{5.1}
straightforwardly.
\endremark
\endremark
\newpage

\shlhetal
\newpage
    
REFERENCES.  
\bibliographystyle{lit-plain}
\bibliography{lista,listb,listx,listf,liste}

\enddocument

\bye